\documentclass[a4paper,12pt,reqno]{amsart}
\usepackage{amssymb}
\usepackage[dvips]{graphicx}

\input epsf.tex
\newdimen\xsize
\newdimen\oldbaselineskip
\newdimen\oldlineskiplimit
\def\restorelineskip{\baselineskip=\oldbaselineskip%
\lineskiplimit=\oldlineskiplimit}
\def\putm[#1][#2]#3{
\hbox{\vbox to 0pt{\parindent=0pt%
\vskip#2\xsize\hbox to0pt{\hskip#1\xsize $#3$\hss}\vss}}}%
\long\def\Line#1{\hbox to \hsize{#1}}
\def\putt[#1][#2]#3{
\vbox to 0pt{\noindent\hskip#1\xsize\lower#2\xsize%
\vtop{\restorelineskip#3}\vss}}

\makeatletter
\def\xbig[#1]#2{{\hbox{$\m@th\left#2\vbox to#1\xsize{}%
\right.\n@space$}}}
\makeatother
\def\xlar[#1]#2{%
\smash{\mathop{ \hbox to #1\xsize{\leftarrowfill}}\limits^{#2}}}
\def\xrar[#1]#2{%
\smash{\mathop{ \hbox to #1\xsize{\rightarrowfill}}\limits^{#2}}}
\def\xline[#1]{\hbox to #1\xsize{\leaders\hrule\hfill}}

\DeclareFontFamily{U}{rsf}{\skewchar\font'177}%
\DeclareFontShape{U}{rsf}{m}{n}{<-6>rsfs5<6-8>rsfs7<8->rsfs10}{}%
\DeclareFontShape{U}{rsf}{b}{n}{<-6>rsfs5<6-8>rsfs7<8->rsfs10}{}%
\DeclareMathAlphabet\RSFS{U}{rsf}{m}{n}
\SetMathAlphabet\RSFS{bold}{U}{rsf}{b}{n}
\def\sf#1{{\mathsf{#1}}}

\def\slsf{\slshape \sffamily }

\def\msmall#1{\mathchoice{\hbox{\small$\displaystyle {#1}$}}{#1}{#1}{#1}}

\def\bb{{\mathbb B}}

\def\cc{{\mathbb C}}

\def\rr{{\mathbb R}}
\def\sss{{\mathbb S}}

\def\pp{{\mathbb P}}
  
\def\zz{{\mathbb Z}}
\def\ttt{{\mathbb T}}

\def\adyn{\sf{1}}

\def\st{_{\mathsf{st}}}
\def\area{\sf{area}}

\def\const{\sf{const}}

\def\dim{\sf{dim}\,}

\def\ev{\sf{ev}}

\def\hol{\sf{Hol}}

\def\ind{\sf{ind}}

\def\hol{\sf{Hol}}

\def\id{\sf{Id}}
\def\im{\sf{Im}\,}

\def\re{\sf{Re}\,}

\def\ker{\sf{Ker}\,}
\def\lim{\mathop{\sf{lim}}}

\def\max{\sf{max}}
\def\min{\sf{min}}

\def\v{{\mathrm{v}}}
\def\w{{\mathrm{w}}}
\def\vol{\sf{Vol}}

\def\eps{\varepsilon}

\def\<{\langle}\let\la=\<
\def\>{\rangle}\let\ra=\>
  
\def\comp{\Subset}
\def\d{\partial}
\def\dbar{{\barr\partial}}

\def\ddef{\mathrel{{=}\raise0.3pt\hbox{:}}}
\def\deff{\mathrel{\raise0.3pt\hbox{\rm:}{=}}}

\def\fraction#1/#2{\mathchoice{{\msmall{ #1\over#2}}}%
{{ #1\over #2 }}{{#1/#2}}{{#1/#2}}}
\def\norm#1{\left\Vert{#1}\right\Vert}

\def\le{\leqslant}

\def\emptyset{\varnothing}

\def\longpoints{\leaders\hbox to 0.5em{\hss.\hss}\hfill \hskip0pt}
\def\stateskip{\smallskip}
\def\state#1. {\stateskip\noindent{\bf#1. }} 
\def\statep#1. {\stateskip\noindent{\bf#1 }} 
\def\proof{\state Proof. \2}

\def\Chi{\raise 2pt\hbox{$\chi$}}

\def\ie{\hskip1pt plus1pt{\sl i.e.\/,\ \hskip1pt plus1pt}}

\def\sli{{\sl i)} } 
\def\slii{{\sl i$\!$i)} } 
\def\sliii{{\sl i$\!$i$\!$i)} }
\def\sliv{{\sl i$\!$v)} }

\def\sing{^\sf{sing}}

\def\barr#1{\mskip1mu\overline{\mskip-1mu{#1}\mskip-1mu}\mskip1mu}
\def\Chi{\raise 2pt\hbox{$\chi$}}

\let\phI=\phi\let\phi=\varphi\let\varphi=\phI
\let\cal=\mathcal

\def\calb{{\cal B}}
\def\calc{{\cal C}}

\def\cale{{\cal E}}
\def\calf{{\cal F}}

\def\calh{{\cal H}}
\def\cali{{\cal I}}

\def\calk{{\cal K}}
\def\call{{\cal L}}

\def\caln{{\cal N}}
\def\calo{{\cal O}}

\def\calr{{\cal R}}
\def\cals{{\cal S}}

\def\calu{{\cal U}}
\def\calv{{\cal V}}
\def\calw{{\cal W}}
\def\calx{{\cal X}}

\def\calz{{\cal Z}}

\def\eps{\varepsilon}

\def\comp{\Subset}
\def\d{\partial}
\def\dbar{{\barr\partial}}
\def\1{{1\mkern-5mu{\rom l}}}

\def\ge{\geqslant}

\def\fraction#1/#2{\mathchoice{{\msmall{ #1\over#2}}}%
{{ #1\over #2 }}{{#1/#2}}{{#1/#2}}}

\def\le{\leqslant}

\def\emptyset{\varnothing}
\newcommand{\2}{\thinspace}

\def\qed{\ \ \hfill\hbox to .1pt{}\hfill\hbox to .1pt{}\hfill $\square$\par}

\def\comment#1\endcomment{}

 \belowdisplayskip=\abovedisplayskip

\def\lineeqqno(#1){\hfill\llap{\vbox to 10pt%
{\vss\begin{align} \eqqno(#1)\end{align}\vss}}\vskip1pt}

\advance\headheight 1.2pt

\def\ShowwLLabel#1{}

\def\thechpt{\Roman{chpt}}
\def\newchapt[#1]#2{%
\refstepcounter{chpt}\setcounter{subsection}{0}%
\setcounter{thm}{0}\setcounter{defi}{0}%
\setcounter{rema}{0}\setcounter{exrc}{0}%
\renewcommand{\thesubsection}{\thechpt.\arabic{subsection}}%
\newpage
\section*{\begin{center}\huge \bf Chapter \thechpt\\
#2 \end{center}}\label{#1}%
\ \smallskip%
\addcontentsline{toc}{chpt}{Chapter \thechpt. #2}%
\markboth{Chapter \thepart}{#2}%
}
\def\newsect[#1]#2{\refstepcounter{section}\setcounter{equation}{0}%
\renewcommand{\thesubsection}{\arabic{section}.\arabic{subsection}}%
\section*{\arabic{section}.
#2}\vspace{-20pt}\label{#1}\vspace{20pt}%
\markboth{Section \arabic{section}}{#2}}

\def\newlect[#1]#2{\refstepcounter{section}%
\renewcommand{\thesubsection}{\arabic{section}.\arabic{subsection}}%
\section*{Lecture \arabic{section}\\
#2}\label{#1}%
\markboth{Lecture \arabic{section}}{#2}}

\def\newprg[#1]#2{\refstepcounter{subsection}%
\subsection*{{\thesubsection.\ #2}} \label{#1}%
}

\def\newappx[#1]#2{%
\refstepcounter{appx}\setcounter{section}{0}%
\renewcommand{\thesubsection}{A\arabic{appx}.\arabic{subsection}}%
\section*{Appendix \arabic{appx}\\ #2}
\label{#1}%
\markboth{Appendix A\arabic{appx}}{#2}
}

\newtheorem{thm}{Theorem}[section]
   \def\newthm#1{\begin{thm}\label{#1}}

\newtheorem{nnthm}{Theorem}
   \def\newthm#1{\begin{nnthm}\label{#1}}

\newtheorem{lem}{Lemma}[section]
   \def\newlemma#1{\begin{lem} \label{#1}}

\newtheorem{prop}{Proposition}[section]
   \def\newprop#1{\begin{prop}\label{#1}}

\newtheorem{nnprop}{Proposition}
   \def\newprop#1{\begin{nnprop}\label{#1}}

\newtheorem{corol}{Corollary}[section]
   \def\newcorol#1{\begin{corol} \label{#1}}

\newtheorem{nncorol}{Corollary}
   \def\newcorol#1{\begin{nncorol} \label{#1}}

\newtheorem{defi}{Definition}[section]
   \def\newdefi#1{\begin{defi} \label{#1}\rm }

\newtheorem{nndefi}{Definition}
   \def\newthm#1{\begin{nndefi}\label{#1}}

\newtheorem{exmp}{Example}[section]
   \def\newexmp#1{\begin{exmp} \label{#1}\rm }

\newtheorem{nnexmp}{Example}
   \def\newthm#1{\begin{nnexmp}\label{#1}}

\newtheorem{exrc}{Exercise}
   \def\newexrc#1{\begin{exrc} \label{#1}\rm }

\newtheorem{rema}{Remark}[section]
   \def\newrema#1{\begin{rema} \label{#1}\rm }

\newtheorem{nnrema}{Remark}
   \def\newthm#1{\begin{nnrema}\label{#1}}

\newtheorem{quest}{Question}

\def\eqqno(#1){\label{(#1)}}
\def\eqqref(#1){(\ref{(#1)})}

\title{Vanishing Cycles in Holomorphic Foliations by Curves and  Foliated Shells}
\author{S. Ivashkovich}
\date{\today}

\address{
Universit\'e de Lille-1, UFR de Math\'ematiques, 59655 Villeneuve
d'Ascq, France} \email{ivachkov@math.univ-lille1.fr}
\address{IAPMM Nat. Acad. Sci. Ukraine
Lviv, Naukova 3b,
79601 Ukraine}


\begin{document}
\begin{abstract}
The purpose of this paper is the study of vanishing cycles in
holomorphic foliations by complex curves on compact complex
manifolds. The main result consists in showing that a vanishing
cycle comes together with a much richer complex geometric object -
we call this object a {\slsf  foliated shell}.
\end{abstract}

\maketitle


\setcounter{tocdepth}{1}

\tableofcontents

\newsect[sect.INT]{Introduction.}

\newprg[prgINT.van-cyc]{Vanishing cycles, compact leaves and simultaneous uniformization}

\smallskip  Let $\call$ be a holomorphic foliation by complex curves on a compact
complex manifold $X$. For the sake of clarity and simplicity of
exposition we describe our results in this Introduction assuming
$\call$ to be smooth. In the main body of the paper this assumption
will be removed (as well, the assumption of compactivity of $X$ will
be replaced by the {\slsf disc-convexity}).

\smallskip Take a point $z\in X$ and denote by $\call_z$ the leaf of
$\call$ passing through $z$. A cycle in $\call_z$  is, by
definition, a closed path (a loop) $\gamma :[0,1]\to \call_z$. A
cycle $\gamma\subset \call_z$ is called {\slsf a vanishing cycle} if
the following two conditions hold:

\begin{itemize}
\item $\gamma$ is not homotopic to zero in $\call_z$;

\smallskip
\item there exist a sequence of points $z_n\to z$ and a sequence of
loops $\gamma_n:[0,1]\to \call_{z_n}$ such that $\gamma_n$ uniformly
converge to $\gamma$ and each $\gamma_n$ is homotopic to zero in $\call_{z_n}$.
\end{itemize}
Classically vanishing cycles became the object of study in foliation
theory since the seminal paper of Novikov \cite{N}, where he used
them to produce a compact leaf in every smooth foliation by surfaces
on $\sss^3$, see also \cite{H}.

\smallskip Apart of the question of existence of compact leaves vanishing
cycles come into a play as obstructions to the simultaneous uniformization
of leaves. Following Il'yashenko, see \S2 in  \cite{Iy2}, take a
smooth complex hypersurface $D$ in $X$ transverse to the leaves of
$\call$. Such $D$ will be called simply {\slsf a transversal} in the
sequel. Set $\call_D=\bigcup_{z\in D}\call_z$ and call this open subset
of $X$ the {\slsf cylinder} of $\call$ over the transversal $D$.

\smallskip Let $\tilde\call_D=\bigcup_{z\in D}\tilde \call_z$ be the union
of the universal coverings of the leaves $\call_z$ equipped with the
natural topology, see Section 3. Let's call $\tilde\call_D$ the {\slsf
universal covering cylinder} (or, simply the {\slsf covering cylinder}
if no misunderstanding can occur) of $\call$ over $D$. It is clear
(see Section 3 for more details) that a leaf $\call_z\subset\call_D$
containing a vanishing cycle exists if and only if the natural
topology of $\tilde\call_D$ is not separable (\ie is not Hausdorff).
Separability of $\tilde\call_D$ means that the leaves of $\call$
which cut $D$ can be simultaneously uniformized. Therefore a
vanishing cycle in some leaf $\call_z\subset\call_D$ is an
obstruction to such simultaneous uniformization. $\call$ is called
{\slsf uniformizable} if for any transversal $D$ the cylinder
$\call_D$ can be uniformized. Therefore $\call$ is uniformizable if
and only if it doesn't contain a vanishing cycle in any of its
leaves. This explains one more reason for the interest in studying
of vanishing cycles.

\newprg[prgINT.sph-she]{Vanishing cycles and foliated shells}

One of the main goals of this paper is to show that  a
vanishing cycle generates a very rich complex geometric object - a
{\slsf foliated shell}.

\smallskip Let $P=\{ z=(z_1,z_2)\in \cc^2: \max \{|z_1|,|z_2|\}\leq 1\}$ be
the unit bicylinder in $\cc^2$ and $B = \{z=(z_1,z_2)\in \cc^2:
\max\{|z_1|,|z_2|\} = 1\}$ its boundary. For some $0<\eps<1$ let
$B^{\eps} = \{ z\in \cc^2: 1-\eps < \max\{|z_1|,|z_2|\}  < 1+\eps
\}$ be a shell around $B$. Denote by $\pi :\cc^2\to\cc$ the
canonical projection $\pi (z)=z_1$ onto the first coordinate of
$\cc^2$. Note that $B^{\eps}$ is foliated by $\pi$ over the disc
$\Delta_{1+\eps}$ of radius $1+\eps$ ($\Delta_r$ denotes the disc of
radius $r>0$ in $\cc$). Denote this foliation by $\call^{\v}$ and call it
{\slsf a vertical foliation}. Its
leaves $\call^{\v}_{z_1}:=\pi^{-1}(z_1)$ are discs $\Delta_{1+\eps}$
if $1-\eps <|z_1|<1+\eps  $ and are annuli
$A_{1-\eps,1+\eps}:=\Delta_{1+\eps}\setminus\bar\Delta_{1-\eps}$ if
$|z_1|\le 1-\eps$.
\begin{nndefi}
The pair $(B^{\eps}, \call^{\v})$ will be called the {\slsf
standard foliated shell}. \label{st-shell}
\end{nndefi}
\begin{figure}[h]
\centering
\includegraphics[width=2.0in]{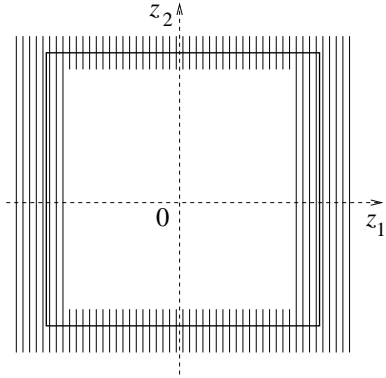}
\caption{The {\slsf  standard foliated shell} is foliated by discs
and annuli over the disc $\Delta_{1+\eps}$. In particular,
$(B^{\eps},\call^{\v})$ is a foliated manifold.}
\label{shell_bid-fig}
\end{figure}
By a foliated manifold in this paper we shall understand  a pair
$(X,\call)$, where $X$ is a complex manifold (separable and
countable at infinity) and $\call$ is a holomorphic foliation by
curves on $X$. Let $(X,\call)$ be a foliated manifold and let
$h:(B^{\eps},\call^{\v})\to (X,\call)$ be a foliated holomorphic
immersion of the standard foliated shell into $(X,\call)$ (an
immersion between two foliated manifolds is called {\slsf  foliated}
if it sends leaves to leaves). Denote by $\Sigma$ the image of the
boundary $B$ under $h$.

\newpage

\begin{nndefi}
The image $h(B^{\eps})$ is called {\slsf  a foliated shell} in
$(X,\call )$ if:

\smallskip
1) immersion $h$ is a generic injection, \ie is such that for all
$z_1\in\Delta_{1+\eps}$ except of a finite set  the restriction
$h|_{\call^{\v}_{z_1}}:\{z_1\}\times A_{1-\eps,1+\eps}\to X$ is an
imbedding;

2) $\Sigma$ is not homologous to zero in $X$.
\label{fol-shell1}
\end{nndefi}

\noindent Roughly speaking the condition (1) means that $h$ is
(much) better then simply an immersion. The main point is of course
the condition (2). It is very strong and our corollaries will
demonstrate this.

\begin{nnexmp}\rm
\label{hopf-st}
The reader should think about the Hopf surface $H^2=\cc^2\setminus
\{0\}/z\sim 2z$. The same vertical foliation $\call^{\v}$ is
invariant under the action $z\sim 2z$ and therefore projects to a
foliation $\call$ on $H^2$. Let $h:\cc^2\setminus \{0\}\to H^2$ be
the canonical projection. It obviously induces a ``foliated
inclusion'' $h:(B^{\eps},\call^{\v})\to (H^2,\call)$. $\Sigma =
h(B)$ is of course not homologous to zero in $H^2$.
\end{nnexmp}

Let $\omega$ be a $(1,1)$-form on $X$. $\omega$ is called {\slsf
pluriclosed} if $dd^c\omega = 0$. Sometimes one calls such $\omega$
also $dd^c$-closed. Recall that $d^c\deff \frac{i}{4\pi}(\bar\d -\d)$
and therefore $dd^c=\frac{i}{2\pi}\d\bar\d$ (in particular
$dd^c\ln |z|^2 = \delta_{0}$). We call a form $\omega$ {\slsf a
taming form for} $\call$ if $\omega|_{\call}>0$. Foliations admitting
a pluriclosed taming form we shall call {\slsf pluritamed}. Our first
result is the following:

\begin{nnthm}
\label{imm-shel}
Let $\call$ be a holomorphic foliation by curves on a compact
complex manifold $X$ which admits a pluriclosed taming form and let
$D$ be a transversal to $\call$ in  $X$. Then the following statements  are
equivalent:

\smallskip\sli Some leaf $\call_{z}\subset\call_D$ contains a
vanishing cycle.

\smallskip\slii The cylinder $\call_D$ contains a foliated shell.
\end{nnthm}

\begin{nnrema}\rm (a)  Statement (ii) means that the mapping $h:B^{\eps}
\to X$, which ``supports'' the foliated shell in $X$, actually takes
values in the cylinder $\call_D$ (but $\Sigma = h(B)$ is not
homologous to zero in the whole of $X$!).

\smallskip\noindent (b) A transversal $D$ is irrelevant in this
theorem: if $\call_{z}$ contains a vanishing cycle then (\slii is
true for {\slsf  every} transversal $D\ni z$.

\smallskip\noindent (c) \rm Recall that a two-dimensional shell in a
complex manifold $X$ is a holomorphic image $\Sigma$ of  $B$ such
that $\Sigma$ is not homologous to zero in $X$. Such shells can
exist only in non-K\"ahler $X$ by the Hartogs-type extension theorem
for K\"ahler manifolds, see \cite{Iv3} (and therefore foliations on
K\"ahler manifolds don't have vanishing cycles). We want to stress
here that $X$ may contain a two-dimensional shell, but it may not be
a {\slsf foliated} shell for the given foliation $\call$. A simple
example is the elliptic fibration on the same Hopf surface $H^2$.
This fibration doesn't admit a foliated shell, while $H^2$ itself
does contain a two-dimensional shell.

\smallskip\noindent (d) In fact in the process of the proof of
Theorem 1 we establish the following useful characterization of
shells:

\begin{nnprop}
\label{imm-shel-prop}
Let $w$ be a $dd^c$-closed
taming form for $\call$. A holomorphic foliated immersion
$h:B^{\eps}\to X$ represents a foliated shell if and only if it is a
generic injection and
\begin{equation}
\int_Bd^c(h^*\omega) \not= 0.
\end{equation}
\end{nnprop}

I.e. not only $h(B)$ is not homologous to zero in $X$ but, moreover,
the distinguished closed $3$-form $d^c\omega$ doesn't vanish on
$h(B)$.  From Proposition \ref{imm-shel-prop} we immediately obtain
the following:

\begin{nncorol}
\label{d-closed} If the taming form $\omega$ of the foliation
$\call$ is $d$-closed then $\call$ has no vanishing cycles.
\end{nncorol}

The strategy of the proof follows that developed for the K\"ahler
case in \cite{Br3} with the replacement of the Thullen type
extension theorem of Siu, \cite{Si1}, by Theorem 1.5 from
\cite{Iv6}, and since all the results of this paper are valid for
disc-convex manifolds it also includes the Stein case as found in
\cite{Iy2}.

\smallskip\noindent (e) The boundary $B$ is topologically the
three-dimensional sphere $\sss^3$. It is not difficult to produce
algebraic (and therefore K\"ahler) manifolds with nontrivial
$\pi_3$, but none of them contains a shell. The reason is that a
shell is a global {\slsf pseudoconvex} object in the complex
manifold $X$ and not simply an element of $\pi_3(X)$.

\smallskip\noindent (f) The meaning of the Theorem 1 is that a
topological property of $(X,\call)$ to contain a vanishing cycle is
equivalent to a complex geometric (even analytic) property to
contain a foliated shell.
\end{nnrema}

\newprg[prgINT.imbed]{Imbedded cycles and imbedded shells}

Note that our foliated shells are, after all, an {\slsf immersed}
objects in $X$ (even if they are ``generic injections''). It would
be definitely preferable to have really an {\slsf imbedded} ones.
However, let us stress at this point that not all foliations with
shells contain {\slsf  imbedded} foliated shells, \ie such that
$h: B^{\eps}\to X$ is an imbedding. The reason is that the
underlying manifold $X$ may not contain an imbedded two-dimensional
shell at all.

\begin{nnexmp}\rm
\label{hopf-lens}
Let, for example, $H^2/(z\sim -z)$ be the quotient of our Hopf
surface by the antipodal involution. The vertical foliation
$\call^{\v}$, described in the Example \ref{hopf-st}, is stable under this
involution and we obtain a foliated manifold $(H^2/\zz_2,\call
/\zz_2)$. The standard foliated shell immerses to $H^2/\zz_2$ and
$B$ maps onto the  quotient $B/\zz_2$ which is topologically a lens
space. I.e., we have here an immersed foliated shell. Due to a
result of Kato, see \cite{K1},  would $H^2/(z\sim -z)$ contain an
imbedded shell, it would be a deformation of a blown-up primary Hopf
surface, \ie its fundamental group would be $\zz$. And this is not
the case, because $\pi_1(H^2/\zz_2)=\zz\rtimes\zz_2$.
\end{nnexmp}
Nevertheless one can find an {\slsf  imbedded} foliated shell provided that:

\smallskip

\begin{itemize}
\item the vanishing cycle $\gamma$ is {\slsf  imbedded} into its leaf
$\call_{z}$;

\smallskip
\item the shell itself is allowed to have somewhat more complicated
topology.
\end{itemize}

Let us more carefully explain what does it mean that
$\gamma\subset \call_{z}$ is imbedded. Let $d$ be the order of
the holonomy of $\call$ along the imbedded loop $\gamma$. It
should be finite, otherwise $\gamma$ cannot be approximated by the
loops $\gamma_n$ in the nearby leaves which are homotopic to zero.
But then for a generic nearby leaf $\call_{z_n}$ the nearby loop
$\gamma_n\subset\call_{z_n}$ will approximate $d\cdot \gamma$ (not
just $\gamma$!) Therefore in the definition of an {\slsf  imbedded
vanishing cycle} one should specify that $\gamma_n\to d\cdot
\gamma$ where $d\ge 1$ is the order of the holonomy of $\call$
along $\gamma$.

\smallskip Now let us turn to the topology of shells. Recall that a
cyclic surface quotient is a normal complex space $\calx^{l,d}$
which is the quotient of $\cc^2$ by the finite group $\Gamma_{l,d}$
of transformations given by $(z_1,z_2)\to (e^{\frac{2\pi i
l}{d}}z_1, e^{\frac{2\pi i }{d}}z_2)$. Here $1\le l<d$ is relatively
prime with $d$. This action preserves the vertical foliation on
$\cc^2$ and therefore $\calx^{l,d}$ is equipped with the
``vertical'' foliation to, which we denote by $\call^{\v}$ again.
Note that the standard ``vertical'' projection $\pi : \calx^{l,d}\to
\cc/<e^{\frac{2\pi i l}{d}}> = \cc$ is well defined and its fibers
are still the leaves of our vertical foliation. Take some smoothly
bounded domain $G\comp \Delta$ such that $\d G\not\ni 0$ but $G\ni
0$ and consider the domain $P=\bigcup_{z\in G}\Delta_z\subset
\calx^{l,d}$ (here $\Delta_z\deff \{z\}\times \Delta$). Remark that
the boundary $B$ of $P$ lies in the smooth part of $\calx^{l,d}$.
For some $\eps >0$ denote by $B^{\eps}$ the $\eps$-neighborhood of
$B$.

\begin{nndefi}
\label{lens-shell} A foliated cyclic shell in $(X,\call)$ is a
foliated holomorphic immersion $h:(B^{\eps},\call^{\v}) \to
(X,\call)$ such that:

1) $h$ is a generic injection;

2) $\Sigma := h(B)$ is not homologous to zero in $X$.
\end{nndefi}

With this notion at hand we can state the following:

\begin{nnthm}
\label{imb-shel}
Let $\call $ be a holomorphic
foliation by curves  on a compact complex manifold $X$ which admits
a pluriclosed taming form and let $D$ be a transversal to $\call$ in
$X$. Then the following conditions are equivalent:

\smallskip\sli Some leaf $\call_{z}\subset\call_D$ contains an imbedded
vanishing cycle.

\smallskip\slii The cylinder $\call_D$  contains an imbedded foliated
cyclic shell.
\end{nnthm}

\smallskip We should point out that the topology of cyclic
shell as we define it can be quite complicated. It is not just a
lens space, \ie is not simply a quotient of $\sss^3$ by a free
action of a finite group.

\smallskip Now we must to explain when the existence of a vanishing cycle in
some leaf $\call_{z}$ of $(X,\call)$ implies the existence of an
{\slsf imbedded} one (in the same leaf). It occurs to depend on
certain ``almost Hartogs'' property of the foliated pair $(X,\call)$, see
Definition  \ref{alm-hart-def}. For the time being let us mention that if
$\omega$ is actually a metric form on $X$ then a foliated pair $(X,\call)$ is
almost Hartogs for every $\call$.

\begin{nnthm}
\label{imb-cycle}
Let $(X,\call,\omega)$ be a pluritamed compact foliated manifold. Suppose additionally
that $\omega$ is a metric form. Then $(X,\call)$ is almost Hartogs. In particular, if
some leaf of $(X, \call)$ contains a vanishing cycle then it contains also
an imbedded vanishing cycle.
\end{nnthm}

\smallskip Almost Hartogs are also all foliated pairs $(X,\call)$, where
the manifold $X$ admits a rational or elliptic fibration, see
Propositions \ref{turb-alm-hart} and \ref{ric-alm-hart}. It is our understanding
of the subject that if an immersed shell in a pluritamed foliated pair $(X,\call)$
is found then the almost Hartogs property of $(X,\call)$ is responsible for
the presence also of an imbedded (but may be cyclic) shell.

\smallskip Let's say a few more words about foliations with shells.
First we remark that shells do come in families. Intuitively
speaking we want to say that if our foliated manifold $(X,\call)$
contains a foliated shell then it breaks into a complex
($\dim_{\cc}X - 2$) - parameter family of ``foliated universes''
each containing a foliated shell.

\begin{figure}[h]
\centering
\includegraphics[width=1.8in]{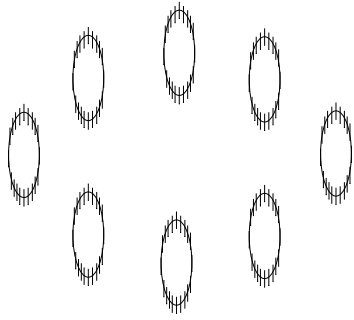}
\caption{Shells persist transversely to $\call$: if there exists a
foliated shell in $(X,\call)$ then it is not disappearing, one can
move it transversely to $\call$.} \label{persist-fig}
\end{figure}

More precisely, the following is true:

\begin{nnprop}
\label{pers-shel}
Let $\call$ be a
holomorphic foliation by curves on a compact manifold $X$ of complex
dimension $n\geq 3$ which admits a pluriclosed taming form. Suppose
that $(X,\call)$ contains a foliated shell
$h:(B^{\eps},\call^{\v})\to (X,\call)$ (imbedded or immersed). Then
there exists a smooth family $\{h_{\lambda}\}_{\lambda\in
\Delta^{n-2}}$ of foliated shells containing $h$ and  transversal to
$\call$ in the sense that:

\begin{itemize}
\item $h_0=h$;
\item $D_{\lambda}h_{0} \left( T_0\Delta^{n-2}\right)\cap D_zh_0(T_zB^{\eps}) =
\{0\}$ for every $z\in B^{\eps}$.
\end{itemize}
\end{nnprop}

Such families of shells clearly come out in our proofs of Theorems 1
and 2. Remark also that due to the equivalence between shells and
vanishing cycles the Proposition \ref{pers-shel} reads also as persistence of
vanishing cycles in an obvious sense. If a two-dimensional
``foliated universe'' is, moreover, compact then it can be listed
explicitly. Namely, the remarkable result of Kato in \cite{K2} (but
even more the ``pseudoconvex surgery'' invented there) allows us to
describe all possible pairs $(X,\call)$, where $X$  is a compact
complex surface, and $\call$ is a holomorphic foliation on $X$ which
contains a vanishing cycle:

\begin{nncorol}
\label{dim-2}
Let $X$ be a compact
complex surface and $\call$ a (singular) holomorphic foliation by
curves such that some leaf $\call_{z}$ of $\call$ contains a
vanishing cycle $\gamma$. Then:

\smallskip \sli  either  $X$ is a modification of a Hopf surface and
$\call_{z}$ is an elliptic curve;

\smallskip \slii or, $X$ is a modification of a Kato surface and
the closure of $\call_{z}$ is a rational curve.
\end{nncorol}

\begin{nnrema}\rm (a) In both cases (\sli and (\slii of this Corollary
the foliated shell in question is either spherical or a (holomorphic)
quotient of the standard $\sss^3$ by $\Gamma_{l,d}$ for some $l,d$.
For the definition of a foliated {\slsf spherical}
shell see Subsection 4.1 (in fact it means that as  the boundary $B$
one can take the standard sphere $\sss^3\subset \cc^2$).  This
Corollary  we state for singular foliations, the reason is that the case (\slii
occurs only for a singular $\call$.

\smallskip\noindent (b) Remark that in the case of surfaces we obtain
a Novikov-type result, \ie the compactivity of the closure of the
leaf supporting a vanishing cycle.

\smallskip\noindent (c) The result, stated in this Corollary, is
obtained also in \cite{Br5}.
\end{nnrema}

\newprg[prgINT.plu-exc]{Pluriexact foliations}

Now let us clarify our assumption on a taming form $\omega$ to be
pluriclosed. Let $T$ be a $(1,1)$-current on $X$ with measurable
coefficients, \ie locally $T=T_{\alpha ,\bar\beta}\frac{\d}{\d z^{\alpha}}
\wedge \frac{\d}{\dbar z^{\beta}}$, where $T_{\alpha ,\bar\beta}$ are
measures. Then there exists a $(1,1)$-vector field $\hat T$
and a complex Radon measure $\| T\|$ such that
\[
<\phi ,T> = \int\limits_X\phi (\hat T_z) d\| T\|(z)
\]
for every test $(1,1)$-form $\phi$, see \cite{HL}.

\begin{nndefi} Following \cite{Su} we shall use the following terminology throughout
this paper:

\begin{itemize}
 \item A $(1,1)$-current $T$ is said to be {\slsf directed by} $\call$ (or,
{\slsf tangent to} $\call$) if for $\| T\|$ - a.a. $z\in X$ one has
$\hat T_z = \frac{i}{2}\v\wedge\bar \v$, where $\v\in T_z\call$.

 \item A {\slsf foliated cycle} on $(X,\call)$ is a {\slsf closed} $(1,1)$-current
directed by $\call$.
\end{itemize}
\end{nndefi}
Remark that a  smooth $(1,1)$-form $\omega$ on $X$ is a {\slsf taming form } for
$\call$ if and only if for every non-trivial positive $(1,1)$-current $T$ directed by
$\call$ one has $<\omega ,T>>0$. In the spirit of \cite{HL} one can prove the following:

\begin{nnprop}
\label{harv-laws}
Let $\call$ be a holomorphic foliation by curves on a compact
complex manifold $X$. Then

\sli either $(X,\call)$ admits a pluriclosed taming form,

\slii or, here exists a non-trivial, positive, $dd^c$-exact
$(1,1)$-current on $X$ directed by $\call$.
\end{nnprop}

\smallskip
I.e., in the case (\slii $(X,\call)$ carries a non-trivial, positive,
$dd^c${\slsf-exact foliated cycle}. A foliated manifold $(X,\call,\omega)$ admitting
a non-trivial, positive, $dd^c$-exact, bidimension $(1,1)$ current
$T$ tangent to (or, directed by) $\call$ we shall call {\slsf
pluriexact}.

\smallskip
Via the aforementioned duality  the characterization result of the
Theorem 1 shows that the class of all holomorphic foliations by
curves on compact complex manifolds splits naturally into the
following three non-intersecting subclasses: the class $\cals$ of
{\slsf  shelled foliations}, the class $\calu$ of {\slsf
uniformizable} foliations, and the class $\cale$ of {\slsf
pluriexact} foliations. A {\slsf  shelled foliation} or a {\slsf
foliation with shells} is a foliation on a compact manifold which
contains foliated shells.

\begin{figure}[h]
\centering
\includegraphics[width=5.0in]{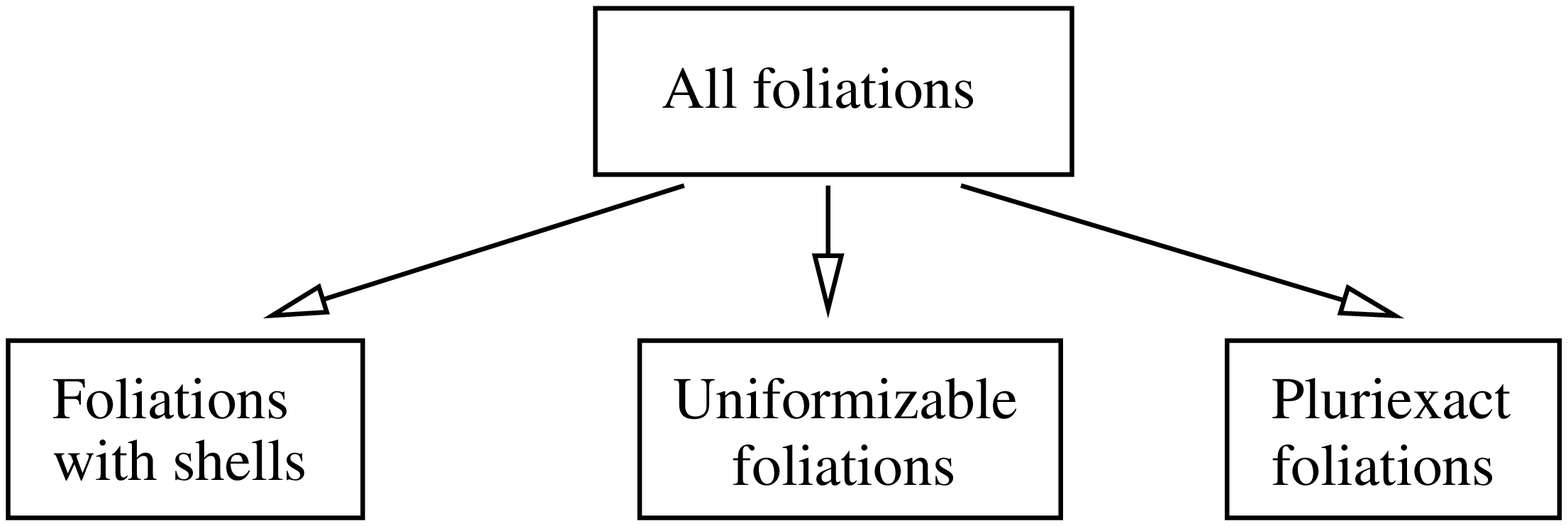}
\caption{} \label{boxes1-fig}
\end{figure}

Note that in the definition of classes $\cals$ and $\calu$ we
require both a pluriclosed taming form and a shell/or absence of
shells. The point is that a foliated shell is of real importance
only in the presence of such a taming form.  In that case it turns
out to be a dominating object in $(X,\call)$. As well as the
uniformizability condition on $\call$ implies more information about
this foliation provided $\call$ admits a pluritaming form. While in
the absence of such a form a $dd^c$-exact current tangent to $\call$
is (eventually) of much greater importance.

\medskip We see the future development of the subject as the study of
each of these classes separately, eventually with the very different
tools, and $\cale$ to be certainly further subdivided. Let us outline
one of the possible ways of doing that. For that remark that a foliated cycle $T$
in a standard way, see \cite{Go}, defines a transverse invariant measure.
Therefore:

\begin{nncorol}
\label{trans-mes}
A pluriexact holomorphic foliation by curves on a compact complex manifold
admits a transversal invariant measure.
\end{nncorol}

Let us try now to say more about this measure. In order to do so let us
define subclasses $\cale_-, \cale_+, \cale_0$ such that $\cale =
\cale_-\bigsqcup \cale_0\bigsqcup \cale_+$ and state our results for
each of them.

\medskip\noindent{\slsf Class $\cale_-$.} The first is the class $\cale_{-}$
 of pluriexact foliations carrying  a
plurinegative taming form, \ie $(X,\call)\in\cale_-$ if it is
pluriexact and if there exists a $(1,1)$-form $\omega$ on $X$ such
that $\omega|_{\call}>0$ and $dd^c\omega\le 0$. Such foliated
manifolds we shall call {\slsf plurinegative}. In Section 5,
Proposition \ref{exact-curr} we shall see that in the  presence of
plurinegative taming form on $(X,\call)$ every vanishing cycle
$\gamma\subset \call_z$ results to a non-trivial exact foliated
cycle, which will be denoted as $T_z$. Its support is contained in
$\bar\call_z$.

\smallskip Now let us describe the way a vanishing cycle appears in a foliated
manifold of class $\cale_-$. In the following theorem we suppose
that the ambient manifold $X$ carries a strictly positive
$dd^c$-closed $(2,2)$-form (this is always the case if $X$ is
compact and $\dim_{\cc}X=3$ for example). Remark finally that $\cale
$ doesn't contain compact complex manifolds of dimension two.
Therefore in the following theorem we deal with $\dim_{\cc}X\ge 3$.

\begin{nnthm}
\label{plurineg} Let $(X,\call , \omega_1)$ be a disc-convex
foliated manifold of class $\cale_-$ and suppose the manifold $X$
itself admits a strictly positive $dd^c$-closed $(2,2)$-form
$\omega_2$. Let $\gamma\subset\call_z$ be a vanishing cycle in the
leaf $\call_z$ of $\call$ and let $T_z$ be the corresponding exact
foliated cycle. Then:

\sli either $(X,\call)$ admits a complex $(\dim_{\cc} X - 2)$-parameter
family of distinct exact foliated cycles, which contains the cycle $T_z$;

\slii or $(X,\call)$ contains a $(\dim_{\cc} X - 3)$-parameter
family of three-dimensional foliated shells.
\end{nnthm}

The notion of three-dimensional foliation shell being clear let us
explain the item (\slii of this Theorem by an example.

\begin{nnexmp}
\label{hopf-three}
Take $H^3=\cc^3/z\sim 2z$ - the Hopf threefold. Let $\call^{\v}$ be again the
vertical foliation $\call_c^{\v}=\{ z_1=c_1, z_2=c_2\}$. $(H^3,\call^{\v})$ admits a
plurinegative taming $(1,1)$-form but doesn't admit a pluriclosed one. It
also contains a $3$-dimensional foliated shell but not a $2$-dimensional one.
\end{nnexmp}

\noindent (i) Indeed, set $z^{'}=(z_1,z_2)$ and consider the following
$(1,1)$-form on $H^3$:

\begin{equation}
\theta = \frac{i}{2}\frac{(dz^{'},dz^{'})}{\norm{z^{'}}^2}.
\end{equation}
$\theta$ is a well defined positive bidimension $(2,2)$-current on
$H^3$. One easily checks that $dd^c\theta = -c_4[\call_0^{\v}]$, where
$[\call_0^{\v}]$ is the current of integration over the central fiber
$\call_0^{\v}$ of $\call^{\v}$ and $c_4$ is the volume of the unit ball in
$\cc^2$. $\theta$ is a clear obstruction to the existence of
pluriclosed $\call^{\v}$-taming form. Indeed, would $\omega$ be such a form
than one would have: $0 = <dd^c\omega , \theta> = <\omega , dd^c\theta> < 0$ -
contradiction.

\smallskip\noindent (ii) At the same time the $(1,1)$-form
\begin{equation}
\omega = \frac{i}{2}\frac{(dz,dz)}{\norm{z}^2},
\end{equation}
where $z=(z_1,z_2,z_3)$, is strictly positive on $H^3$ (not only on
$\call$, \ie is a metric form) and one easily checks that
$dd^c\omega\leq 0$ (but $dd^c\omega \not=0$ contrary to the
two-dimensional case). I.e. $\omega$ serves as a plurinegative
taming form for any foliation by curves on $H^3$. Our foliated Hopf
manifold $(H^3,\call^{\v})$ contains an obvious foliated
$3$-dimensional shell, but doesn't contain {\slsf any}
$2$-dimensional shell because $H_3(H^3,\zz) = 0$. I.e., the case
(\slii of Theorem \ref{plurineg} realizes here.

\begin{nnrema} \rm
A behavior as in the case (\sli of Theorem \ref{plurineg} appears in
Example 3.6 of \cite{Iv6}, where also some other relevant examples
can be found.
\end{nnrema}

\medskip\noindent{\slsf Class $\cale_+$.} A foliated manifold $(X ,\call)$ we call
{\slsf pluripositive} if there exists a non-trivial $(1,1)$-current $T$ tangent to
$\call$ such that $T=dd^cR$ for some {\slsf positive} $(2,2)$-current $R$.  We denote
the class of pluripositive foliations as $\cale_+$. Note that obviously
$\cale_-\cap \cale_+=\emptyset$. Our main result on class $\cale_+$ is the following:

\begin{nnthm}
\label{pluripos}
Let $(X,\call)$ be a disc-convex foliated manifold which possesses a
non-trivial, positive $(1,1)$-current $T$ directed by $\call$ such that $T=dd^cS$
for some positive current $S$. Then:

\sli $\chi_{\call^{s}}T=0$,

\slii the transversal measure $\mu$ induced by $T$ on $X^0\deff X\setminus \call^{s}$
has finite logarithmic potential.
\end{nnthm}
Here $\call^s \deff \call^{\sing}$ is the singular locus of $\call$ (\ie it is admitted in
this formulation that $\call$ can be singular). Remark also that for the Hopf foliated pair
$(H^3,\call^{\v})$ of Example \ref{hopf-three} the transverse measure is the delta function.
This makes contrast to the statement of Theorem \ref{pluripos} for the pluripositive foliations.

\medskip\noindent{\slsf Class $\cale_0$.}  We define this class simply as
$\cale_0 \deff \cale\setminus \left(\cale_-\bigsqcup
\cale_+\right)$. I.e., $\cale_0$ is the class of foliated manifolds
which do not admit neither a plurinegative taming form no a positive
$dd^c$-exact foliated cycle $T$. Surprisingly this can happen.
Namely we shall see that:

\begin{nnexmp}
\label{mojsheson}
There exists a compact Moishezon $3$-fold $X$ and a holomorphic foliation by curves $\call$
on $X$ such that $(X,\call)\in \cale_0$.
\end{nnexmp}

We also give a characterization of class $\cale_0$ in terms of
currents, see Propositions \ref {harv-laws-neg} and
\ref{harv-laws-fol}.

\newprg[prgINT.unif]{Uniformizable foliations}

For uniformizable  foliations tamed by a pluriclosed form we expect
more or less the same results as for foliations on compact K\"ahler
(or algebraic) manifolds. Let's give some typical statements.

\begin{nncorol}
\label{raz-fibr}
Let $\call$ be a holomorphic
foliation by curves on a compact complex manifold $X$ which admits a
plurinegative taming form. Suppose that $\call$ contains a leaf
whose universal cover is $\cc\pp^1$. Then the universal cover of
every leaf is $\cc\pp^1$ and, moreover, $\call$ is a rational
quasi-fibration.
\end{nncorol}

\smallskip
For foliations on K\"ahler manifolds this result is proved in
\cite{Br3}. In Section \ref{sect.EXOQ} we prove also the following version of the
Reeb stability theorem:

\begin{nnprop}
\label{reeb}
Let $\call$ be a
holomorphic foliation by curves on a compact complex manifold $X$
admitting a $dd^c$-negative taming form.

\sli If $\call$ has a compact leaf with finite holonomy then all
leaves of $\call$ are compact with finite holonomy.

\slii If every leaf of $\call$ is compact then every leaf has finite
holonomy. In that case there is an upper bound on volumes of leaves
and the leaf space is Hausdorff.
\end{nnprop}

For the foliations on K\"ahler manifolds this result is well known,
see \cite{Ga,P}. Both statements (and some others) follow from the
compactness property of the corresponding cycle space in the
presence of a plurinegative taming form, see Remark \ref{ratgeo2}.

\smallskip Take a cycle $\gamma$ on some leaf $\call_z$. Following
\cite{LP}, see also \cite{Iy1}, we define in Subsection 6.3 the
domain of preservation of the homotopy class $[\gamma]$ and prove
the following:

\begin{nnprop}
\label{land-petr}
Let $\call$ be a holomorphic foliation by curves on a compact
complex manifold $X$ admitting a pluriclosed taming  form. Then:

\smallskip
1) either the domain of preservation $\Omega_{\gamma }$ is Hausdorff
(and therefore is a complex manifold) for every loop $\gamma $,

\smallskip
2) or $(X,\call )$  contains a foliated shell.
\end{nnprop}

\newprg[prg1.6]{The structure of the paper, notes} This paper is organized
in the following way.

\smallskip\noindent 1. First, in Section 2 we develop the main technical
tool - a meromorphic extension Theorem \ref{reparalldim2} from
generalized Hartogs domains. This theorem should be viewed as a
generalization of Theorem 2.2 from [Iv6] and Proposition 4.1 from
[Br5] at a time. The new points here are replacement of {\slsf
metric} forms by {\slsf taming} ones (this includes the replacement
of the Siu's  Thullen type extension theorem by Lemma \ref{area1} -
a nonparametric version of Theorem 1.5 from \cite{Iv6}) and
detecting the obstructions to the extension as {\slsf foliated
shells} and not simply as {\slsf shells} like in \cite{Iv6} - the
corresponding arguments are gathered in Lemmas \ref{ext-shel} and
\ref{descrete}.

\smallskip\noindent 2. We start the Section 3 with recalling the necessary
definitions and notions around uniformization of foliations,
vanishing ends and covering cylinders (tubes), which were developed
by M. Brunella in \cite{Br1}-\cite{Br4}. The right notions were
worked out in the sited papers in part after appearance of example
in [CI], see discussion before Theorem 3.1 in \cite{Br4}. This
example is recalled and enhanced in subsection 2.3 of the present
paper in the context of extension theorems. Its relevance to the
vanishing ends is discussed in Remark \ref{exempl2}. After that  we
prove a more precise version of Theorem 1, namely the Theorem 3.1,
which includes the case of singular foliations, non-compact ambient
manifolds and, more crucially, specifies the location of a foliated
shell. A more precise version of Theorem 2 - the Theorem 3.2 about
imbedded cycles and shells follows in Subsection 3.6.

\smallskip\noindent 3. Section 4 contains the proof of Theorem 3. More
generally it describes the techniques to obtain an {\slsf imbedded}
vanishing cycle from a non-imbedded one. The main achievement is
Theorem 4.1, which relates an "almost Hartogs" property of a
foliated pair $(X,\call)$ with imbedded vanishing cycles.
It contains also the description of examples of complex
dimension two, \ie of complex surfaces. In particular, that means
the Corollary \ref{dim-2}.

\smallskip The logic how to obtain from immersed cycles  the imbedded
ones is explained in \cite{Br5}, where the case of nonuniformizable
foliations on complex surfaces is studied. Our paper follows this
logic, but applies it in all dimensions for disc convex manifolds,
and more crucially for taming (and not only metric forms). The
advantage of this becomes clear via Proposition 3, \ie the
complementary class $\cale$ also is quite understandable.

\smallskip\noindent 4. Section 5 is devoted to pluriexact foliations
and contains the proof of Theorems \ref{plurineg} and \ref{pluripos}
as well as Example \ref{mojsheson}.

\smallskip\noindent 5. The last Section 6 contains multidimensional
examples relevant to the subject of this paper, proofs of Corollary
\ref{raz-fibr} and Propositions \ref{reeb} and \ref{land-petr}. Here
we also formulate several open questions.

\medskip\noindent{\slsf  Acknowledgement.} 1. This research in its final
part was done during the Authors stay in the Max-Planck-Institute
f\"ur Mathematik, Bonn. I would like to thank this Institution for
the hospitality.

\noindent 2. I would like also to give my thanks to Stefan
Nemirovski and Vsevolod Shevchishin for useful discussions on the
subject of this paper.

\noindent 3. The results of this paper had been reported  at the
conference ``Recent  Developments in non-K\"ahler Geometry'' held on
the honor of Masahide Kato at the University of Hokkaido, Sapporo on
5-7 March 2008, see \cite{NK}. I would like to thank the organizers
for the possibility to give this talk and for the useful post-talk
discussions.

\noindent 4. I'm also grateful to the Referee of this paper for the
remarks and useful suggestions.

\break

\newsect[sect.EXT]{Pluriclosed Taming Forms and Foliated Immersions.}

\newprg[prgEXT.hart-f]{Generalized Hartogs figures}

Let us start with some definitions. In the definition of a {\slsf
foliated manifold} $(X,\call)$ from now we aloud $\call$ to be a
{\slsf singular} holomorphic foliation by curves on $X$. One of the
ways to define such $\call$ is the following. Take a sufficiently
fine open covering $\{\Omega_{\alpha}\}$ of $X$. Then $\call$ will
be defined by the nonvanishing identically holomorphic vector fields
$v_{\alpha}\in \calo (\Omega_{\alpha}, TX)$ which are related on a
non-empty intersections
$\Omega_{\alpha,\beta}\deff\Omega_{\alpha}\cap \Omega_{\beta}$ as
$v_{\alpha}=h_{\alpha,\beta} v_{\beta}$. Here $h_{\alpha,\beta}\in
\calo^*(\Omega_{\alpha,\beta})$. After contracting the common
factors one immediately sees that the singular set $\call^{\sing}$
of $\call$, which is defined as $\call^{\sing}=
\{z:v_{\alpha}(z)=0\}$, is an analytic subset of $X$ of codimension
at least two. Set $X^0\deff X\setminus\call^{\sing}$. The leaves of
$\call$ are, in the first approximation, defined as the leaves of
the smooth foliation  $\call^0\deff \call|_{X^0}$, \ie they are
entirely off the singular set of $\call$. Then, depending on the
someone goals, one adds to them some ``ends''. We shall do that in
the following Section.

\smallskip A particular class of foliated manifolds are fibrations by
curves, \ie triples $(W,\pi ,V)$ where $W$ is a complex manifold
of dimension $\dim V +1$ and $\pi :W\to V$ is a surjective
holomorphic submersion with connected fibers. A holomorphic mapping
$f:(X,\call)\to (X',\call')$ is said to be a foliated immersion if
it is an immersion and sends leaves to leaves. In the case of fibrations,
\ie if $f:(W,\pi,V)\to (W',\pi',V')$, one can be more precise: there
exists a holomorphic map $f_{\v}:V\to V'$ such that for all $z\in V$
one has $f(W_z)\subset W^{'}_{f_{\v}(z)}$. Dimension of $W'$
might be bigger then that of $W$. If $V'=V$ one often supposes also
that $W_z$ goes to $W^{'}_z$ for all $z\in V$. This will be clear from
the context.

\begin{defi}
\label{genhart}
A generalized Hartogs figure is a quadruple $(W,\pi ,U,V)$, where
$W$ and $V$ are complex manifolds, $U$ an open subset of $V$ and
$\pi :W\to V$ is a holomorphic submersion such that:

\smallskip \sli for all $z\in V\setminus U$ the fiber $W_z=\pi^{-1}(z)$
is diffeomorphic to an annulus;

\slii for $z\in U$ the fiber $W_z$ is diffeomorphic to a disc.
\end{defi}
Generalized Hartogs figures are foliated manifolds (even fibrations)
of a special type: they are concave in the most na\"ive and clear
sense.  Manifold $W$ has a distinguished part of the boundary formed
by the outer boundaries $\d_0 W_z$ of annuli $W_z$. We shall suppose
that $W$ is smooth up to this part of its boundary and denote it by
$\d_0 W$, \ie $\d_0 W=\cup_{z\in V}\d_0 W_z$. Projection $\pi$ is
also supposed to be smooth up to $\d_0 W$ and therefore $\pi :\d
W_0\to V$ is a circle fibration. For $z\in U$ the outer boundary
$\d_0 W_z$ is actually the boundary of the disc $W_z$.

\smallskip Recall that the {\slsf standard Hartogs figure} is
the open subset of $\cc^{n+1}$ of the form

\begin{equation}
H_{\eps} = \left( \Delta^n_{1+\eps}\times A_{1-\eps,1+\eps}\right)
\cup \left( \Delta^n_{\eps}\times\Delta_{1+\eps}\right)
\eqqno(hart)
\end{equation}
for some $\eps >0$. $H_{\eps}$ likewise carries our ``vertical
foliation'' $\call^{\v}$. This time the leaves $\call^{\v}_{z'}$ are
discs $\Delta_{1+\eps}$ if $||z'||< \eps$ and annuli for $\eps \le
||z'||< 1+\eps$. Here $z'=(z_1,...,z_n)$ and $||\cdot||$ is the
polydisc-norm in $\cc^n$. Remark now that $(H_{\eps},\call^{\v})$
fits, of course, into the Definition \ref{genhart} with
$V=\Delta^n_{1+\eps}$, $U=\Delta^n_{\eps}$ and $\pi$ being the
restriction of the canonical ``vertical''  projection $\cc^{n+1}\to
\cc^n$ onto $H_{\eps}$. Remark furthermore that the standard
foliated shell is also a generalized Hartogs figure. Namely it can
be viewed as $(B^{\eps},\pi , A_{1-\eps,1+\eps},\Delta_{1+\eps})$.

\begin{defi}
If $U=\emptyset $ we call $(W,\pi ,\emptyset ,V)$ trivial, if
$U=V$ we call $(W,\pi ,V,V)$ complete and in the latter case
often denote it as $(W,\pi ,V)$.
\end{defi}
The standard Hartogs figure is newer trivial by definition, \ie it is
commonly accepted that always $\eps >0$. Let $D$ be a non-empty open
subset of $V$. Set $W|_D=\pi^{-1}(D)$ and consider it also as a
generalized Hartogs figure $(W|_D,\pi|_{D},D\cap U,D)$ - a subfigure
of $(W,\pi,U,V)$. Moreover, if $S\subset V$ is a submanifold of $V$
one can consider the restriction $(W|_S,\pi|_{S},S\cap U,S)$ and it
is again a generalized Hartogs figure. In the sequel we shall often
avoid the word ``generalized'' and call our figures simply {\slsf Hartogs
figures}, specifying over what $V$ they are considered.

\begin{rema}\rm
The necessity of considering generalized Hartogs figures in this
paper comes from the simple observation  that: {\it every vanishing
cycle produces a natural generalized (or topological) Hartogs figure
around it}. This will become clear in Section 3. The fact that the
standard Hartogs figure is not sufficient for our considerations
will be explained by an example in Subsection \ref{prg2.3}.
\end{rema}

\newprg[prgEXT.rep]{Extension after a reparametrization}

The following notion comes back to \cite{Ti}, see also \cite{Bl}.
Let $f:A_{1-\eps,1}\to X$ be a holomorphic immersion.
\begin{defi}
We say that $f$ extends to $\Delta$ after a reparametrization if for
some $\delta >0$ there exists an imbedding $h:A_{1-\delta,1}\to
A_{1-\eps,1}$ sending $\d\Delta$ to $\d\Delta$ and preserving the canonical
orientation of $\d\Delta$, such that $f\circ h$
holomorphically extends to $\Delta$.
\end{defi}

It is clear that such $h$, if exists, should be holomorphic. We
shall use also the following form of this notion. Let $\gamma$ be a
simple oriented loop on a bordered Riemann surface $W$. The latter
should be viewed simply as a collar adjacent to $\gamma$. Let
$f:W\to X$ be a holomorphic immersion. Suppose that there exist a
Riemann surface $\widetilde{W}$ which is a bordered disc with
boundary $\tilde \gamma$ (canonically oriented) and a biholomorphic
mapping $h$ from a collar adjacent to $\tilde\gamma$ to $W$ (smooth
up to the boundaries) and sending $\tilde\gamma$ onto $\gamma$,
preserving orientations, such that the composition $f\circ h$
holomorphically extends onto the disc $\widetilde{W}$. Then we shall
say that $f$ extends onto the disc $\widetilde{W}$ after a
reparametrization. If such $\widetilde{W}, \tilde\gamma$ and $h$ do
exist but are not specified we shall say simply that $f$
holomorphically extends onto a disc after a reparametrization.
\begin{rema}\rm
(a) There is one case when the extension of an immersion after a
reparametri\-zation  may be not unique in the sense that there may not
exist an automorphism $\phi$ of $\Delta$ such that one extension is equal
to the another one composed with $\phi$. For example, take a function
$f(z)=4z+\sqrt{z^2-1}$ and consider it as a holomorphic mapping from
a thin annulus around $\d \Delta (2)$ - the circle of radius $2$,
into $\cc\pp^1$. Then $f$ has two extensions after a
reparametrization:

\smallskip 1) An injective one. Indeed, $f$ is an
imbedding of $\d \Delta (2)$ into $\cc\pp^1$ and therefore bounds a
disc, say $D$. Let $r:\Delta (2)\to D$ be a Riemann mapping (it is
biholomorphic in a neighborhoods of the closed discs). Set
$h=f^{-1}\circ r$ - a reparametrization  of $\d \Delta (2)$. Then
$f\circ h = r$ is the extension of $f$ onto $\Delta (2)$ after a
reparametrization.

\smallskip 2) A non-injective one. This is given by the formula defining
$f$. It has two ramification points $\pm 1$ and extends onto the union
$\Delta (2)\cup \cc\pp^1$ appropriately glued along the slit
$[-1,1]$. The Riemann surface obtained is again a disc. This second
extension is non-injective.

\smallskip\noindent (b) At the same tame, if $f$ was a {\sl generic
injection} (i.e. injective outside of a finite set) then its extension
after reparametrization, which we also require to be a generic injection,
is unique (if exists). Uniqueness means here up to a biholomorphic
automorphism of the disc.
\end{rema}

Now let's turn ourselves to the families of immersions.

\begin{defi}
\label{geninject}
A holomorphic mapping $f:(W,\pi ,V)\to X$  of a fibration $(W,\pi,V)$ into
a complex manifold $X$ is called {\it generically injective} if for all
$z\in V$ outside of a proper analytic subset $A\subset V$ the restriction
$f_z:=f|_{W_z}$ is a generic injection.
\end{defi}
Note that we do not ask $f$ to be "generically injective" itself but
only its restrictions onto "generic" fibers. Actually $f$ may not be
even an immersion. However in most cases mappings appearing in this
paper will be both immersions and generic injections. We shall also
need a corresponding notion for the meromorphic case.
\begin{defi}
\label{merinj}
A meromorphic mapping $f:W\to X$ between complex manifolds is a meromorphic
immersion if it is an immersion outside of its indeterminacy set $I_f$. If,
moreover, $(W,\pi,V)$ is a holomorphic fibration  then a meromorphic mapping
$f$ is called generically injective if $f|_{W_z}$ is a generic injection for
$z$ outside of a proper analytic subset of $V$.
\end{defi}
Here and always in this paper writing $f(C)$ for some meromorphic map
and some complex curve $C$ we mean that the restriction $f|_C$
of $f$ onto $C$ is well defined (this means that $C$ in not contained
in the indeterminacy set of $f$) and $f(C)$ is actually $f|_C(C)$.
Again we will mostly work with meromorphic maps which are both meromorphic
immersions and generic injections.

\smallskip Let a holomorphic generic injection $f:(W,\pi ,U ,V)\to X$ of a
generalized Hartogs figure into a complex manifold $X$ be given and let
$\hat U$ be some open subset of $V$ containing $U$.

\begin{defi}
\label{reparexten}
We say that $f$ meromorphically extends onto the
Hartogs figure $(\widetilde{W}, \pi , \hat U,$ $ V)$ over (the
same!) $V$   {\sl after a reparametrization } if there exists a
foliated biholomorphism of trivial figures $h:(\d_0\widetilde{W},
\pi ,\emptyset , V)\to (\d_0W,\pi ,\emptyset , V)$ (\ie $h$ is
defined in an one-sided neighborhood of $\d_0\widetilde{W}$ and
$h(z)$ tends to $\d_0W$ when $z$ tends to $\d_0\widetilde{W}$) such
that $f\circ h$ extends to a generically injective meromorphic map
$\tilde f : (\widetilde{W} , \pi , \hat U, V)\to X$. A
reparametrization map $h$ is supposed to be constant over $V$, \ie
$h(\d_0\widetilde{W}_{z_1})\subset W_{z_1}$ for all $z_1\in V$.
\end{defi}

\begin{rema}\rm
(a) We say simply that $f$ extends after a reparametrization
if the data as in Definition 2.6 do exist (but, may be, are not
specified). In that case we often omit tildes over the extended
objects, such as $W$ and $f$ ($\pi$ will newer come with a tilde in
this context). I.e., we often say that $f$ extends onto $(W,\pi ,
\hat U, V)$ after a reparametrization.

\smallskip\noindent (b) Remark that if $f$ extends as a meromorphic map being
a generic injection on $(W,\pi,U,V)$ with $U\not=\emptyset$ then its
extension will be automatically a generic injection. However in the
definition above we do not exclude the case when $U=\emptyset$.
\end{rema}

\begin{defi}
If for any point $z\in V$ there exists a neighborhood $V(z)$
such that the restriction $f|_{W|_{V(z)}}$ extends after a
reparametrization onto a complete Hartogs figure
$\big(\widetilde{W}|_{V(z)}, \pi , V(z)\big)$ over $V(z)$ then we
say that $f$ locally extends after a reparametrization.
\end{defi}

Let us be very precise at this point: by saying that the ``restriction
$f|_{W|_{V(z)}}$ extends'' we mean here that one is taking the restriction
of $f$ onto the Hartogs subfigure $(W|_{V(z)},\pi,V(z)\cap U, V(z))$ and this
restriction extends onto the complete figure $(W|_{V(z)}, \pi , V(z))$. I.e., a
reparametrization is supposed to be made near $\d_0W|_{V(z)}$ only.

If a generically injective (!) mapping $f$ extends locally then it extends
globally. Namely the following is true:

\begin{lem}
\label{localization}
Let $V \supset U_1\cup U_2$ with  (may be empty) intersection
$U_{12}\deff U_1\cap U_2$. Let $(W,\pi ,\emptyset , V)$ be a trivial
generalized Hartogs figure over $V$ and
$f:W\to X$ be a generically injective  holomorphic map into a
complex manifold $X$ such that $f$ meromorphically extends onto a
complete Hartogs figures $(W|_{U_k},\pi ,U_k)$ for $k=1,2$ after a
reparametrization. Then $f$ extends after a reparametrization onto a
figure $(W,\pi , U_1\cup U_2, V)$.
\end{lem}
\proof {\sl Step 1. Extending onto a complete figure $(W|_{U_1\cup
U_2},\pi,U_1\cup U_2)$.} Denote by
$h_k:(\d_0\widetilde{W}_k,\pi_k,\emptyset ,U_k)\to (\d_0W|_{U_k},\pi
,\emptyset ,U_k)$ the corresponding foliated biholomorphisms. The
fact that $f$ is a generic injection imply that $h_2^{-1}\circ
h_1=(f\circ h_2)^{-1}\circ (f\circ h_1) :\d_0\widetilde{W}_{1,z}\to
\d_0\widetilde{W}_{2,z}$ extends for every $z\in U_{12}$ onto a
corresponding disc and therefore extends to a foliated
biholomorphism between complete figures
$(\widetilde{W}_1|_{U_{12}},\pi_1 , U_{12})$ and
$(\widetilde{W}_2|_{U_{12}},\pi_2 ,U_{12})$. Therefore complete
figures $(\widetilde{W}_1,\pi_1 , U_1)$ and $(\tilde{W_2},\pi_2 ,
U_2)$ glue together to a complete figure $(\widetilde{W}_3,\pi ,
U_1\cup U_2)$ and reparametrization maps $h_k$ glue to a
reparametrization map $h:(\d_0\widetilde{W}_3,\pi , \emptyset ,
U_1\cup U_2)\to (\d_0W|_{U_1\cup U_2},\pi , \emptyset ,U_1\cup
U_2)$. Mapping $f$ extends, after being reparameterized by $h$ onto
the complete figure $(W|_{U_1\cup U_2},\pi,U_1\cup U_2)$.

\smallskip\noindent{\sl Step 2. Extending to $(W,\pi,U_1\cup U_2,V)$.}
Reparametrization $h:\d_0\widetilde{W}_3\to \d_0W|_{U_1\cup U_2}$
constructed in the first Step allows us to glue figures
$(\widetilde{W}_3,\pi , U_1\cup U_2)$ and $(W,\pi,\emptyset ,V)$
together to a figure $(W,\pi,U_1\cup U_2,V)$ and $f$ extends onto
it.

\smallskip\qed

From this lemma we obtain the following

\begin{corol}
\label{globalization} Let  $f:W\to X$ be a generically injective
holomorphic immersion of a generalized Hartogs figure $(W,\pi ,U,V)$
into a complex manifold $X$. Then there exists a maximal open
$U\subset \hat U\subset V$ such that $f$ meromorphically extends
onto  $(W,\pi , \hat U, V)$ after a reparametrization.
\end{corol}

\begin{figure}[h]
\centering
\includegraphics[width=2.5in]{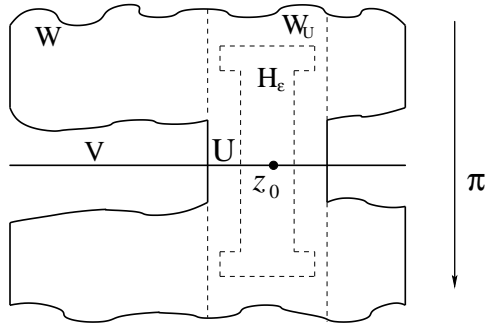}
\caption{This is the  standard Hartogs figure imbedded into  a
generalized Hartogs figure $W$ constructed in  \cite{CI}. Every
attempt to imbed $H_{\eps}$ into this $W$ will look like on the
picture: if the fiber over the origin in $H_{\eps}$ is mapped to a
fiber over some point  $z_0\in U$ then the image of $H_{\eps}$ will
newer leave $W|_U$.} \label{hartogs-fig}
\end{figure}

\newprg[prg2.3]{Hartogs figures and reparametrizations.}
The necessity of considering generalized Hartogs figures in
foliation theory will be very clearly seen along this paper. The
fact that the work with them cannot be reduced to the standard
Hartogs figures in $\cc^2$ (or $\cc^n$) comes from the example
constructed in \cite{CI}. This example is explained on the Figure
\ref{hartogs-fig}. Namely the following can happen. There exists a
generalized Hartogs figure $W$ over a disc (\ie both $\emptyset\not
= U \subset V$ are discs in $\cc$) with the following property:
whenever a holomorphic foliated imbedding $h:(z_1,z_2)\to
(h_1(z_1),h_2(z_1,z_2)$) of $H_{\eps}$ into $W$ is given such that
$h_1(0)=z^0\in U$ then necessarily $h_1(\Delta_{1+\eps})\subset U$
(whatever $\eps > 0$ is).

This feature of the example in question is not explicitly stated in
\cite{CI} (a somewhat weaker property of it was sufficient there)
and therefore we shall enhance (and simplify considerably) our example
in this subsection. Let us briefly recall the construction
from \cite{CI}.

\smallskip\noindent{\slsf Step 1. Construction.} Consider the
following domain
\[
W=\Delta\times\Delta\setminus \{ (z_1,z_2): 1/3\le |z_1|\le 2/3
\text{ and } z_2^2=z_1 \lambda (|z_1|^2) \}.
\]
Here $\lambda\in \calc^{\infty} (\rr), 0\le \lambda\le 1, \lambda
(t) \equiv 0 \text{ for } t<1/9 \text{ and } \lambda (t) \equiv 1
\text{ for } t>4/9\}$.

\smallskip In $W$ consider the (almost) complex structure $J$
defined by the basis of $(1,0)$-forms $dz_1$ and $dz_2 +
a(z_1,z_2)d\bar z_1$ where
\begin{equation}
\eqqno(a) a(z_1,z_2) =
\begin{cases}
\frac{z_2z_1^{2}\lambda'(|z_1|^2)}{ z_2^2 - z_1\lambda (|z_1|^2)}
\text{ if } 1/3\le |z_1|\le 2/3, \cr 0 \text{ otherwise} .
\end{cases}
\end{equation}
In \cite{CI}, see Lemma 1, it was proved that:

\sli $J$ is integrable;

\slii $J=J\st$ on $W\setminus \big(\bar
A_{\frac{1}{3},\frac{2}{3}}\times\Delta \big)$ - this is obvious
from the definition of $J$;

\sliii functions $g(z)=z_1$ and $f(z)=z_2 + \frac{z_1}{z_2}\lambda
(|z_1|^2)$ are $J$-holomorphic;

\sliv $\ind_{|z_2|=1-\delta}f(z_1,z_2)=-1$ for every $|z_1|=1$ and
every $\delta >0$ small enough ($0<\delta <1/6$ is sufficient).

\begin{rema}\rm (a) The construction of this example could be easily
understood looking on the function $f$. It is designed to have
properties incompatible with being "holomorphic". Then one computes
its differential and gets
\[
df(z_1,z_2) = \big(\frac{\lambda}{z_2} +\lambda'\frac{z_1}{z_2}\bar
z_1\big)dz_1 + \big(1-\frac{z_1}{z_2^2}\lambda \big)\big(dz_2 +
\frac{z_2z_1^2\lambda'}{z_2^2 - z_1\lambda}d\bar z_1\big).
\]
Differential $dz_1$ should be of the type $(1,0)$ in order to have a
holomorphic fibration $(z_1,z_2)\to z_1$. All what is left is to
take the structure such that $df$ is a $(1,0)$-form. Now the choice
for $a$ made in \eqqref(a) becomes clear.

\smallskip\noindent (b) Item (\sli follows immediately from the fact
that $J$ admits two holomorphic functions. Item (\sliv will be not
used here.
\end{rema}

\smallskip Remark that $(W, \pi , \Delta_{1/3}, \Delta_{2/3})$, equipped
with the structure $J$, is a generalized Hartogs figure ($\pi$ is
the natural projection $(z_1,z_2)\to z_1$).

\smallskip\noindent{\slsf Step 2. Properties of $(W,J)$.}  Let us
list some important properties of this example.

\medskip
\sli Remark that $W_{z_1}$ has the conformal structure of a
(pluri-punctured) disc. This follows from the observation that
$f(z_1,z_2)$ ($z_1$ is fixed) is holomorphic (meromorphic at zero)
on $z_2$ with respect to both structures: $J$ and $J\st$.

\smallskip
\slii The most important feature is the following

\begin{prop}
For any $\eps>0$ and any foliated holomorphic imbedding $h(z_1,z_2)
= (h_1(z_1),h_2(z_1,z_2))$ of the standard Hartogs figure $H_{\eps}$
into $W$ such that $h_1(0)\in \Delta_{1/3}$ one has
$h_1(\Delta_{1+\eps})\subset \Delta_{1/3}$.
\end{prop}
\proof Set $\beta (z_1) = \lambda (|h_1(z_1)|^2)$. A disc
$\Delta_{z_1}\deff \{z_1\}\times \Delta_{1+\eps}$ (if $|z_1|<\eps$)
or an annulus $A_{z_1}\deff \{z_1\}\times A_{1-\eps,1+\eps}$ (if
$|z_1|\ge \eps$) is mapped by $h$ to the fiber $W_{h_1(z_1)}\deff
\pi^{-1}(h(z_1))$.

\smallskip The holomorphic mapping $\Phi\deff (g,f):W\to
\cc\times\cc\pp^1$ realizes $W$ as a Riemann domain over
$\cc\times\cc\pp^1$. The composition $\Phi\circ h$ extends to a
bidisc $\Delta^2_{1+\eps}$ by the classical Hartogs-Levi theorem and
this extension is a foliated meromorphic immersion. The last is
because the zero divisor of the Jacobian (if non-empty) should have
intersect $H_{\eps}$, but there our map is locally biholomorphic.

\smallskip This is a contradiction, because on each fiber $\Delta_{z_1}$,
$1/3 <|z_1|\le 2/3$, the meromorphic function $(f\circ h)(z_1,z_2) =
h_2(z_1,z_2) + \frac{h_1(z_1)}{h_2(z_1,z_2)}\lambda (|h(z_1)|^2)$ is
a composition of a holomorphic function $h_2$ and Jukovskiy function
$z+\frac{\text{ const }}{z}$ and the last is not an immersion.

\smallskip\qed

\smallskip
\sliii Along the same lines as in the proof of this Proposition one
can show that: for any $z_1^0\in \d\Delta_{1/3}$ and any $\eps
>0$ there exists {\slsf no imbedding} of $W|_{\Delta (z_1^0,\eps)}$ into
$\cc^2$, \ie there no holomorphic coordinates on $W|_{\Delta
(z_1^0,\eps)}$.

\newprg[prgEXT.key]{Key lemma}

\medskip The following  statement is the key lemma for this section.
It replaces the Lemma 2.3 from \cite{Iv6} for the ``unparametric''
case of the present paper. In fact it states even more general result
that we need in this paper, see Remarks \ref{general1} and \ref{general2}.
But we are including
it for the future references. We suppose that our complex manifold
$X$ is equipped with some Riemannian metric. The condition (\slii in
the following lemma, where this metric is used, doesn't depend, in
fact, on a particular choice of a metric.

\begin{lem}
\label{area1} Let $f:W\to X$ be a generically injective holomorphic
map of a trivial Hartogs figure $(W,\pi,\emptyset,V)$ into a complex
manifold $X$. Suppose that $\dim V=1$ and that for some sequence
$z_n\to z_0\in V$ restrictions $f|_{W_{z_n}}:W_{z_n}\to X$
holomorphically extend as generic injections onto a discs
$\widetilde{W}_{z_n}$ after a reparametrization. Suppose
additionally that:

\smallskip
\sli $\tilde f|_{\widetilde{W}_{z_n}}(\widetilde{W}_{z_n})$  stay in
some compact of $X$;

\slii $\area\left( \tilde
f|_{\widetilde{W}_{z_n}}(\widetilde{W}_{z_n})\right)$ are uniformly
bounded.

\smallskip Then there exists a neighborhood $D\ni z_0$ such that $f$
extends meromorphically onto a figure $(\widetilde{W},\pi,D,V)$
after a reparametrization. Moreover, the extension $\tilde f$ is a
generically injective meromorphic map.
\end{lem}
\proof Writing $\tilde f|_{\widetilde{W}_{z_n}}$ in the statement of
this lemma we mean that for every $n$ a reparametrization map
$h_{z_n}:\d \widetilde{W}_{z_n}\to\d W_{z_n}$ is given such that
$\tilde f|_{\d\widetilde{W}_{z_n}}\deff  f|_{\d W_{z_n}}\circ
h_{z_n}$ extends generically injectively and  holomorphically to the
disc $\widetilde{W}_{z_n}$. The proof will use in a crucial way the
description of convergence of analytic discs obtained in \cite{IS3}
and the structures of Banach neighborhoods of stable curves obtained
in \cite{IS1}.

\smallskip Set $f_n=f|_{\widetilde{W}_{z_n}}$ and consider them as
complex discs over $X$, parameterized by a fixed disc $\Sigma$ (see
\S3 from \cite{IS1} or \S2 from \cite{IS3} for exact definitions).
Applying Theorem 1 from \cite{IS3} we can find a subsequence from
$\{f_n\}$ that converge in the sense of Definition 2.5 from
\cite{IS3} to a stable curve $f_0$ over $X$, parameterized again by
a disc. Be careful, this $f_0(\Sigma)$ may have compact components.

\smallskip By Theorem 3.4 from \cite{IS1} the space of discs over $X$
which are close to $f_0$ is a Banach analytic set of finite
codimension. Denote it by $\calc$. By the Theorem of Ramis, see
\cite{Ra}, $\calc$ is the union of finitely many irreducible
components $\calc_j$ and each $\calc_j$ is a finite ramified
covering over a Banach ball. Take a component which contains
infinitely many of $f_n$-s. In order not to complicate our notations
we suppose that $\calc$ is irreducible itself and contains all
$f_n$.

\smallskip For the sequel it is important to understand how $\calc$ was
constructed in \cite{IS1}:

\smallskip 1) The parametrizing disc $\Sigma$ is covered by finite number
of discs, annuli and pants $\Sigma_{\alpha}$ (the boundary annulus is one
of them, denote it as $\Sigma_{\alpha_0}$). This covering has that property
that each intersecting pair $\Sigma_{\alpha}, \Sigma_{\beta}$ intersect by
an annulus denoted as $\Sigma_{\alpha ,\beta}$.

\smallskip 2) For each $\alpha$ a Banach manifold $H_{\alpha}$ of holomorphic
maps from $\Sigma_{\alpha}$ to $X$ is considered. The same type manifolds
$H_{\alpha ,\beta}$ of holomorphic maps $\Sigma_{\alpha , \beta}\to X$ for
intersecting $\Sigma_{\alpha}$ and $\Sigma_{\beta}$ are considered.

\smallskip 3) For every pair of intersecting $\Sigma_{\alpha}$ and $\Sigma_{\beta}$
a transition map $\psi_{\alpha ,\beta}: H_{\alpha }\times H_{\beta}\to H_{\alpha ,
\beta}$ is defined.

\medskip Now $\calc$ comes out as the zero set of some ``gluing'' holomorphic map
$\Psi = \{ \psi_{\alpha ,\beta}\} : \bigcup_{\alpha}H_{\alpha }
\to \bigcup_{\alpha}H_{\alpha , \beta} $. By construction $\calc$ restricts as a
Banach analytic subset to each of $H_{\alpha}$.

\medskip All what is left to do is to replace $H_{\alpha_0}$ (the manifold of maps
from the annulus adjacent to the boundary) by a $1$-dimensional
manifold $\calw :=\{f|_{\d W_z}: z \text{ in a neighborhood}$
$\text{ of } z_0\}$ ($\d W_z$ stays here for an annulus adjacent to
$\d_0W_z$). The obtained Banach analytic set, we still denote it as
$\calc$, is of finite dimension (the proof goes along the same lines
as the proof of Lemma 1.1 from \cite{Iv6}). In fact it is clearly of
dimension not more than one. But since it contains the sequence
$\{f|_{\widetilde{W}_{z_n}}\}$ its dimension is actually one.
Therefore $\calc$ is a usual analytic set by Barlet-Mazet theorem,
\cite{M}, \ie is a complex curve in our case. Restriction
$\calc\to\calw$ is an analytic map and it is {\sl proper} (!),
because a nondegenerate analytic maps between complex curves are
always proper. Therefore its image is the whole $\calw$. We get an
extension $\tilde f_z$ for all $z$ close to $z_0$ as a family by a
tautological map $\tilde f:\tilde\calw \to X$. Here $\tilde\calw$ is
a tautological family of discs over $\calw$.

\smallskip

\qed

\begin{rema}\rm
\label{general1} (a) An analogous statement can be proved also in the
case $\dim V\geq 2$, but then one should require the boundedness of
rational cycle geometry of $X$ as in \cite{Iv6} (only cycles tangent
to $\call$ are relevant). We shall do this later, see Lemma
\ref{area2}.

\smallskip\noindent (b) One can seriously simplify the proof of this
lemma if one imposes {\sl ad hoc} the condition that (some
subsequence of) the sequence $\{ \tilde
f|_{\widetilde{W}_{z_n}}(\widetilde{W}_{z_n})\}$ converges to $f_0$
without bubbles. The proof is then almost immediate, since then only
one Banach manifold $H_0$ appears (no Banach analytic sets), that of
deformations of $f_0$ and it has dimension at least one because it
contains a converging sequence.
\end{rema}

\newprg[prgEXT.dim-two]{Two-dimensional case}
Recall that a complex manifold $X$ is called disc-convex if for any
compact $K\subset X$ the exists a compact $\hat K$ in $X$ such that
for any holomorphic map $\phi :\overline{\Delta}\to X$  such that
$\phi (\d \Delta)\subset K$ one has $\phi (\overline{\Delta})
\subset \hat K$. Let's adapt this notion to the foliation theory:

\begin{defi}
\label{disc-conv} A complex foliated manifold $(X,\call)$ is called
disc-convex if for any compact $K\subset X$ the exists a compact
$\hat K_{\call}$ in $X$ such that for any holomorphic map $\phi
:\overline{\Delta}\to X$ tangent to $\call$ and  such that $\phi (\d
\Delta)\subset K$ one has $\phi (\overline{\Delta}) \subset \hat
K_{\call}$.
\end{defi}

A holomorphic mapping $\phi :\overline{\Delta}\to X$ is called
tangent to $\call$ if it takes (almost all) its values in some leaf
of $\call$. Note that for disc-convex $(X,\call)$ and foliated
mappings $f:(W,\pi,\emptyset,V)\to (X,\call)$ the condition (\sli
in Lemma \ref{area1} is satisfied automatically.

\smallskip Let  $\omega$ be a $(1,1)$-form on $X$.
\begin{defi}
We call $\omega$ plurinegative ($dd^c$-negative) if $dd^c\omega\leq 0$.
We call $\omega$ pluriclosed ($dd^c$-closed) if $dd^c\omega = 0$.
\end{defi}

\smallskip Denote by
$\cale_{\rr}^{p,p}$ the sheaf of smooth real $(p,p)$-forms and by
$\cale^{p,g}$ the sheaf of smooth complex valued $(p,q)$-forms on
$X$. Likewise by $\cale_{p,q}$ we denote the dual to $\cale^{p,g}$
space of currents of bidimension $(p,q)$  and by $\cale^{\rr}_{p,p}$
the space of real currents of bidimension $(p,p)$.

Fix some strictly positive $(1,1)$-form $\Omega$ on $X$. Given a
holomorphic foliation by curves $\call$ on $X$ define the following
convex compact $K_{1,1}(\call) \subset \cale_{1,1}^{\rr}(X)$. For
every point $z\in X^0$ take a  $(1,1)$-vector $\frac{i}{2}\v\land
\bar\v$ tangent to $\call_z^0$ such that $<\frac{i}{2}\v\land \bar\v
,\Omega>=1$.

For $z\in\call^{\sing}$ take any sequence $z_n\to z$, $z_n\in X^0$
and any $\frac{i}{2}\v_n\land \bar\v_n$ tangent to $\call_{z_n}^0$
such that $<\frac{i}{2}\v_n\land \bar\v_n ,\Omega>=1$ for all $n$.
Subtract a converging subsequence from $\frac{i}{2}\v_n\land
\bar\v_n$ and denote by $\frac{i}{2}\v_0\land \bar\v_0$ its limit.
In this way we obtain all positive bidimension $(1,1)$
$\delta$-currents tangent to $\call$. $K_{1,1}(\call)$ is the
closure of the convex hull of these $\delta$-currents.

\begin{defi}
A $(1,1)$-form  $\omega$ we call a taming form for $\call$ if
$<T,\omega>>0$  for every $T\in K_{1,1}(\call)$.
\end{defi}

The Theorem below plays the role of  Lemmas 3.1 from \cite{Iv4} and 2.4
from \cite{Iv6} for the generalized Hartogs figures. The proof closely
follows \cite{Iv6}.

\begin{thm}
\label{repar2dim}
Let $(X,\call)$ be a disc-convex foliated complex
manifold which admits a $dd^c$-negative taming form $\omega $ and let
$f:(W,\pi ,\emptyset,V)\to (X,\call)$ be a generically injective
foliated holomorphic mapping from a trivial two-dimensional Hartogs
figure $W$ over a disc $V$ into $X$. Suppose that:

\sli  for some sequence  $z_n\to z_0\in V$ restrictions $f_n\deff
f|_{W_{z_n}}$ extend  onto discs $\widetilde{W}_{z_n}$ after a
reparametrization as generic injections;

\slii $\area\left( \tilde
f|_{\widetilde{W}_{z_n}}(\widetilde{W}_{z_n})\right)$ are uniformly
bounded.

\smallskip\noindent Then mapping $f$ extends after a reparametrization
as a generically injective meromorphic map onto a complete Hartogs
figure $(\widetilde{W}, \pi ,V)$ over $V$ minus a closed "essential
singularity" set $S$ of the form $S=\bigcup_{z\in S_1}S_z$, where
$S_1$ is a closed complete polar subset of $V$ and $S_z$ for every
$z\in S_1$ is a compact in the disc $\widetilde{W}_z$.
\end{thm}
\proof   Denote by $\hat U$ the maximal open subset of $V$ such that
$f$ meromorphically extends onto a Hartogs figure $(W, \pi
,\hat U,V)$ after a reparametrization (tildes are skipped everywhere).
By Corollary \ref{globalization} we know
that $\hat U$ exists and by Lemma \ref{area1} that $\hat U\not=\emptyset$.
All we need to prove is that $\d \hat U\cap V$ is a polar set. The question
is local with respect to the base $V$ and therefore we fix
$z_1^0\in \d\hat U\cap \Delta$ and suppose that $V$ is a disc around $z_1^0$.
We denote by $T$ the current $f^*\omega $. It is a smooth $(1,1)$-form
outside of a discrete set $A$ of points of indeterminacy of $f$
which is $dd^c$-closed everywhere on $W$ as a current.

\smallskip\noindent{\slsf Step 1. Laplacian of the area function.}
For points $z_1\in \hat U$ the following area function is well
defined:

\begin{equation}
\label{areaf1}
a(z_1) = \area\left(f|_{W_{z_1}}(W_{z_1})\right) = \int\limits_{W_{z_1}}T,
\end{equation}
here $W_{z_1}=\pi^{-1}(z_1)$. Let $\rho$ be a bump function, equal
to $\adyn$ ion a neighborhood of $\d_0W_V$. Fix a disc $D\subset
\hat U$ and denote by $\pi : \bar W_D\to D$ the natural projection.
Here witting $\bar W_D$ we mean $W_D\cup \d_0W_0$. The restriction
$T|_{\bar W_D}$ we denote still as $T$. Write

\begin{equation}
\eqqno(areaT)
dd^c a =  dd^c\left[\pi_*\rho T + \pi_*(1-\rho)T\right] = dd^ca_{\rho} +
dd^c[\pi_*(1-\rho)T],
\end{equation}
where $\pi_*$ stands for the push forward operator of the
restriction $\pi|_{\bar W_D}:\bar W_D\to D$.  $a_{\rho}$ denotes here the
area function, which corresponds to the form $\rho T$. Remark that
$\pi_*(\rho T)$ is well defined and smooth on the whole of $V$ and
therefore such is also $dd^ca_{\rho}$. What concerns the second term in
\eqqref(areaT) it is equal to $
-\pi_*\left[dd^c\rho\wedge T + 2d\rho \wedge d^cT + (\rho -1)dd^cT\right]$.
As a result we see that

\begin{equation}
\eqqno(laplac-a)
dd^ca = dd^ca_{\rho} -\pi_*\left[dd^c\rho\wedge T + 2d\rho \wedge d^cT +
(\rho -1)dd^cT\right].
\end{equation}
Since $T$ is plurinegative we get that
\begin{equation}
\eqqno(laplac-a-neg)
dd^ca \le dd^ca_{\rho} -\pi_*\left[dd^c\rho\wedge T + 2d\rho \wedge d^cT\right]
\end{equation}
Remark that the right hand side of this expression is also well
defined and smooth on the whole of $V$. \eqqref(laplac-a-neg) shows,
in particular, that the Laplacian of $a$ is bounded from above by a
function, denote it as $c$, which smoothly extends from $\hat U$ to
$V$.

\begin{rema}\rm
\label{smooth} {\bf (a)} All computations on this stage were done on
$\d_0 W|_V$ minus a discrete set of points of indeterminacy of $f$.
At points $z_1$ such that $W_{z_1}$ cuts an indeterminacy point of
$f$ the area function $a$ jumps. Usually this small "non-smoothness"
of $a$ plays no role in the forthcoming considerations. But, while
working with pluriclosed forms it will be convenient to remark that
$a$ can be considered to be smooth everywhere due to the following
observation. First: $a$ is clearly smooth outside of a discrete set.
Second: subtracting from it the Poisson integral of the right hand
side of \eqqref(laplac-a-neg) (which becomes to be an equality if
$\omega$ is pluriclosed), we get a bounded harmonic function which
smoothly extends through these discrete set. Therefore we can
"correct" $a$ on a discrete set of points to make it genuinely
smooth.

\smallskip\noindent{\bf (b)} Alternatively, following \cite{Br5},
one can exploit approach from \cite{Iv6} and compute the
Laplacian of $a$ in such a way that the resulting expression
involves only the boundary integrals.
\end{rema}

\smallskip\noindent{\slsf Step 2. Polarity of $\d\hat U\cap \Delta$.}
Area function $a$ writes as $a(z_1) = b(z_1) + h(z_1)$, where $b$ is
the Poisson integral of $c$ (and therefore is smooth on $V$), and
$h\deff a-b$ has nonpositive Laplacian, \ie is superharmonic in
$\hat U$. In addition $h$ is bounded from below on $\hat U$. All
what is left to do is to remark that $a(z_1)$ and therefore $h(z_1)$
tend to $+\infty$ when $z_1\to \d\hat U\cap V$ by Lemma \ref{area1}.
Therefore we are in the conditions of the proof of Lemma 3.1 from
\cite{Iv4} or Lemma 2.4 from \cite{Iv6}. Indeed $\d \hat U\cap V$
occurs to be a $+\infty$ set for a superharmonic in $\hat U$
function $h$. But then setting

\[
h_n(z_1) = \min\{n,h(z_1)\}
\]
we get an increasing sequence of superharmonic function {\slsf in
the whole of} V (!), which therefore converges to a superharmonic
function on $V$, equal to $h$ on $\hat U$ and to $+\infty $ on $\d
\hat U\cap V$. One concludes that $\d\hat U \cap V$ is complete
polar in $V$ (see Lemma 2.4 from \cite{Iv6} for more details about
this) and one sets $S_1\deff \d\hat U \cap V$.

\smallskip\qed

\begin{rema}\rm
\label{general2}
One can modify the proof of Theorem \ref{repar2dim} along the lines
of \cite{Iv3} and then the Remark \ref{general1} (b) would be sufficient.
But anyway, this doesn't make the proof shorter and we shall not do that.
\end{rema}

\newprg[prgEXT.shells]{Condition $\int d^cT=0$ and foliated shells}

In this subsection we suppose the timing form $\omega$ is pluriclosed
and that $f$ is already extended onto
$W\setminus S$, see Theorem \ref{repar2dim} (tildes are skipped everywhere).
The polar set $\d\hat U\cap V$ we had denote as $S_1$. Therefore  the
"essential singularity" set $S$ of the extended map is actually
$S=\bigcup_{z\in S_1}S_z$. Now we shall see how this leads to a
foliated shell. We suppose therefore that $S$ is not empty. Take a
point $s_0\in S_1$ and let $V$ be a small disc around $s_0$. Shrinking
$W_V$ if necessary we can suppose that $W_V$ is a bidisc $\Delta^2$
in $\cc^2$ with coordinates $z_1,z_2$. Write $T=it^{\alpha \bar\beta}
dz_{\alpha}dz_{\bar\beta}$ for $T\deff f^*\omega$.

\begin{lem}
\label{ext-shel}
Suppose that the timing form $\omega$ is pluriclosed and that for a relatively
compact disc $D\comp V$ such that $\d D\cap S_1=\emptyset$ one has
\begin{equation}
\eqqno(decete)
\int\limits_{\d W|_D}d^cT = 0.
\end{equation}
Then $f$ meromorphically extends onto $W_D$.
\end{lem}
\proof We know already that  $dd^c a$ smoothly extends onto $D$ (and
$a$ is positive!). Denote by $\widetilde{dd^ca}$ this extension.
Again supposing that $W_D$ is a bidisc we can compute the integral
of $d^cT$ over components $\d_0W_D$ and $W_{\d D}$ of the boundary
$\d W_D$ of $W_D$, and using the expression (2.2.5) from \cite{Iv6}
for the Laplacian $dd^ca$, get:

\begin{equation}
\eqqno(deceT1) 2\int\limits_{\d_0 W_D} d^c T = - \int\limits_{D}
\widetilde{dd^ca} + \frac{1}{4\pi} \int\limits_{\d_0W_{
D}}\left(\frac{\d t^{1,\bar 1}}{\d z_2}dz_2 - \frac{\d t^{1,\bar
1}}{\d\bar z_2}d\bar z_2\right)\wedge dz_1\wedge d\bar z_1
\end{equation}
and

\begin{equation}
\eqqno(deceT2)
\int\limits_{W_{\d D}} d^c T = \int_{\d D}d^ca + \frac{1}{4\pi}
\int\limits_{\d_0 W_{\d D}}t^{1\bar 2}dz_1\wedge d\bar z_2 - t^{2\bar 1}
dz_2\wedge d\bar z_1.
\end{equation}
The second term in the right hand side of \eqqref(deceT2) is the
integration over $\d_0 W_{\d D}$ as the boundary of $W_{\d D}$.
Considering it as the boundary of $\d_0W_D$ (and thus changing its
orientation) and applying Stokes we get after summing up
\eqqref(deceT1) with \eqqref(deceT2) that

\begin{equation}
2\int\limits_{\d_0 W_D} d^c T + \int\limits_{W_{\d D}} d^c T=\int_{\d D}d^ca
- \int\limits_{D} \widetilde{dd^ca} + \int\limits_{\d_0 W_D} d^c T,
\end{equation}
\ie that
\begin{equation}
\eqqno(dca-ddca)
\int\limits_{\d W_D} d^c T = \int_{\d D}d^ca - \int\limits_{D} \widetilde{dd^ca}.
\end{equation}

Therefore for the negative measure  $\mu_h \deff dd^c h$ supported
on $S_1$ we obtain (using smoothing by convolutions) the identity:
\begin{equation}
\eqqno(laplac-h)
\mu_h(S_1\cap D) = \int\limits_Ddd^c h = \int\limits_{\d D}d^ca -
\int\limits_D\widetilde{dd^c a} = \int\limits_{\d W_D} d^c T.
\end{equation}
The assumption of our Lemma mean now that $\mu_h(S_1\cap D) = 0$, \ie that $h$
is harmonic. Therefore $a$ is smooth and Lemma \ref{area1} implies now the
extendibility of $f$ onto $W_D$.

\smallskip\qed

We conclude with the following

\begin{corol}
If in the conditions of Theorem \ref{repar2dim} the mapping $f$ is additionally
supposed to be an immersion, taming form $\omega$ to be pluriclosed, and the
singularity set $S$ of the extended mapping is non-empty then $(X,\call)$ contains
a foliated shell.
\end{corol}

Really, the extended mapping $\tilde f:\widetilde{W}\setminus S \to
X$ might fail to be an immersion in this case only on a discrete
subset of $\widetilde{W}\setminus (S\cup I_{\tilde f})$. We add this
subset to $S$ together with indeterminacy set $I_{\tilde f}$ of
$\tilde f$ to get $\tilde S$. The projection $\tilde S_1=\pi (\tilde
S)$ will stay complete polar. Let $0\in S$. Take a small disc
$\Delta_r$ around $0$ in such a way that $\d\Delta_r\subset
V\setminus \tilde S_1$ and such that $\widetilde{W}|_{\Delta_r}$ is
a bidisc (after a slight shrinking of its outer boundary
$\d_0\widetilde{W}|_{\Delta_r}$),\ie
$\widetilde{W}|_{\Delta_r}=\Delta_r\times\Delta$ as foliated
manifolds. In $\widetilde{W}|_{\Delta_r}$ take a bidisc
$P=\Delta_r\times\Delta_{1-\eps}$ for $\eps>0$ small enough to
insure the immersivity of $\tilde f$ near the boundary $B$ of $P$.
Lemma \ref{ext-shel} says now that $\int_{\tilde f(B)}\omega\not=0$,
\ie we got a foliated shell.

\smallskip Remark that we also proved  the Proposition \ref{imm-shel-prop}
from the Introduction:

\begin{corol}
A generic holomorphic injection $h:(B^{\eps},\call^{\v})\to
(X^0,\call^0,\omega)$ into a disc-convex pluritamed foliated
manifold defines a foliated shell if and only if it is an immersion
and
\begin{equation}
\int_{h(B)}d^c\omega \not= 0.
\end{equation}
\end{corol}

Finally let us proof that

\begin{lem}
\label{descrete}
If the timing form $\omega$ is pluriclosed then the set $S_1$ is
at most countable.
\end{lem}
\proof We use the notations of the proof of Lemma  \ref{ext-shel}, \ie $h$
denotes the
superharmonic extension of $a-b$ to $\Delta = V$. Consider the following
representation of $\pi_1(\Delta\setminus S_1)$ in $\rr$:
\begin{equation}
\eqqno(pres-1) R_1:\pi_1(\Delta\setminus S_1)\ni \gamma \to
\int\limits_{\gamma}d^ch.
\end{equation}
We can always suppose that $\gamma$ has only normal crossings and
denote by $\Gamma$ its interior, \ie a finite number of discs. By
\eqqref(laplac-h) we have that
\[
\int\limits_{\gamma}d^ch = \int\limits_{\d W_{\Gamma}}d^cT =
\int\limits_{f(\d W_{\Gamma})}d^c\omega .
\]
Therefore the image of the representation $R_1$ is contained in the
image of the representation
\begin{equation}
\eqqno(pres-2) R_2: H_3(X,\zz) \ni M \to \int\limits_Md^c\omega ,
\end{equation}
and the last is countable. Therefore $\im R_1$ is at most countable.
Would $S_1$ be uncountable one would be able to find a nontrivial $\gamma \in
\pi_1(\Delta\setminus S_1)$ such that
\[
\int_{\gamma} d^ch = 0 \text{ and therefore } \int_{\Gamma} dd^ch = 0.
\]
But then for any component $D$ of $\Gamma$ we would have
\[
\int_{\d W_D}d^cT= \int_{D} dd^ch = 0,
\]
and this by Lemma \ref{ext-shel} implies that $f$ extends to $W_D$,
\ie that $S_1\cap D=\emptyset$. Since this is true for every component
of $\Gamma$ we get that $\gamma$ is trivial in $\pi_1(\Delta\setminus S_1$).
Contradiction.

\smallskip\qed

\newprg[prgEXT.non-pam]{Nonparametric extension in all dimensions}

Let $n\geq 1$ be the dimension of the base $V$. In the following
formulations tildes are skipped everywhere.

\begin{thm}
\label{reparalldim1} Let $(X,\call)$ be a disc-convex foliated
manifold admitting a plurinegative taming form $\omega $ and let
$f:(W,\pi ,U,V)\to (X,\call)$ be a generically injective foliated
holomorphic map from a non-trivial Hartogs figure (\ie
$U\not=\emptyset $) into $X$. Then $f$ extends after a
reparametrization  to a foliated meromorphic map $\tilde
f:(\widetilde{W},\pi,V)\setminus S\to X$ of a complete Hartogs
figure minus a closed subset $S$ of the form $S=\cup_{z\in S_1}S_z$,
where:

\smallskip $(a_1)$ $S_1$ is a complete $(n-1)$-polar subset of $V$ of
Hausdorff dimension $2n-2$.

\smallskip $(a_2)$ $S_z$ is a compact in the disc $\widetilde{W}_z$ for
every $z\in S_1$.

\smallskip $(a_3)$ If $\dim W = \dim X$ and $f$ was an immersion then
the extended map $\tilde f$ is a meromorphic immersion.
\end{thm}

\noindent Complete $(n-1)$-polarity of $S_1$ means that every point
$0\in S_1$ admits a neighborhood $\Delta^n=
\Delta^{n-1}\times\Delta$ with coordinates $(\lambda ,z_1)$ such
that for every $\lambda$ the disc $\Delta_{\lambda} \deff
\{\lambda\}\times \Delta$ intersects $S_1$ by a complete polar
compact set. Hausdorff zero-dimensionality of $\Delta_{\lambda}\cap
S_1$ follows. For the purposes of this paper we will need to know
more about the behavior of the extended map $\tilde f$ near the
essential singularity set $S$. Supposing that $S$ is nonempty take a
point $s_0\in S$ and find a coordinate $n$-disc $\Delta^n\ni 0 = \pi
(s_0)\in S_1$ and a neighborhood $P$ of $s_0$ in $W$ biholomorphic
to $\Delta^n\times\Delta$ such that $\pi|_{P}:P\to \Delta^n$ is the
natural vertical projection $\Delta^n\times \Delta\to\Delta^n$. In
what follows $z_2$ will denote the coordinate along the fiber of
$\pi$.

\begin{thm}
\label{reparalldim2}
Under the conditions of Theorem 2.2 suppose
additionally that the taming form $\omega$ is pluriclosed and that
the singular set $S$ is nonempty. Then $S$  has the following structure:

\smallskip $(b_1)$ The coordinate polydisc $\Delta^n$ as above can be
presented as $\Delta^n=\Delta^{n-1}\times\Delta$ with coordinates
$(\lambda,z_1)$ in such a way that the restriction to $S\cap
\Delta^{n+1}$ of an another vertical projection
$\pi_1:\Delta^{n+1}\to \Delta^{n-1}$, \ie of $(\lambda,z_1,z_2)\to
\lambda$ is proper and surjective. In another words for every
$\lambda\in \Delta^{n-1}$ the intersection $S_{\lambda}\deff
\Delta^2_{\lambda} \cap S$ is nonempty.

\smallskip $(b_2)$ For every $\lambda\in \Delta^{n-1}$ let
$B_{\lambda}=\d\Delta^2_{\lambda}$. Then $\tilde f(\d B_{\lambda}))$
is not homologous to zero in $X$, \ie it is a foliated shell in
$(X,\call)$, provided that $f$ was in addition an immersion.

\smallskip $(b_3)$ Denote by $\pi_2:\Delta^n\to \Delta^{n-1}$ the
natural projection $(\lambda , z_1)\to \lambda $. Then for every
$\lambda $ the set $S_{1,\lambda }\deff \pi_2^{-1}(\lambda ) \cap
S_1$  is a most countable.
\end{thm}

\noindent In all applications/formulations $S$ will be supposed to
be minimal closed such that $f$ extends onto $\widetilde{W}\setminus
S$. The proof is not a direct generalization of the two-dimensional
case. First of all we need to introduce one object relevant to a
complex foliated manifold $(X,\call)$. Denote by $\calr_{\call}$ the
analytic space of rational cycles on $X$ tangent to $\call$. Recall
that a rational cycle is a finite linear combination of rational
curves with integer coefficients: $C=\sum_jn_jC_j$. Here each $C_j$
is a rational curve in $X$. We fix a Hermitian metric on $X$ and
denote by $\omega$ its associated $(1,1)$-form. The area of $C$ is
defined as
\begin{equation}
\label{areafunction1}
\v_{\omega}(\calc) = \sum_jn_j\int_{C_j}\omega
.
\end{equation}
\begin{defi}
Let us say that $(X,\call)$ has unbounded rational cycle geometry if
there exists a path $\gamma :[0,1[\to \calr_{\call}$ such that

\smallskip 1) $C_{\gamma (t)}$ stays in some compact $K$ of $X$ for
all $t\in [0,1[$;

\smallskip 2) $\v_{\omega} (C_{\gamma (t)})\to +\infty $ when
$t\nearrow 1$.
\end{defi}
\noindent Here $C_{\gamma (t)}$ is the rational cycle in $X$
corresponding to the point $\gamma (t)\in \calr_{\call}$. This
notion doesn't depend on the particular choice of $\omega$ and
represents from itself a pure complex-geometric property of
$(X,\call)$.

\smallskip In Lemma \ref{ratgeo} we shall prove that if $(X,\call)$
admits a plurinegative taming form then the rational cycle geometry
of $(X,\call)$ is bounded. Recall finally, that a subset $A$ of a
complex manifold $V$ is said to be {\sl thick at} the point $z^0\in
V$ if for any neighborhood $U$ of $z^0$ $A\cap U$ is not contained
in a proper analytic subset of $U$. Now we can state the needed
lemma:

\begin{lem}
\label{area2} Let $f:W\to X$ be a generically injective foliated
holomorphic mapping  of a trivial Hartogs figure
$(W,\pi,\emptyset,V)$ into a complex foliated manifold $(X,\call)$.
Suppose that $\dim V\ge 2$ and that for all $z$ in some subset
$A\subset V$ thick at $z^0$ all restrictions $f|_{W_{z}}:W_{z}\to X$
holomorphically extend onto discs $\widetilde{W}_{z}$  after a
reparametrization. Suppose additionally that:

\smallskip
\sli $\area\left( \tilde
f|_{\widetilde{W}_{z}}(\widetilde{W}_{z})\right)$ is uniformly
bounded for $z\in A$;

\slii  $\tilde f|_{\widetilde{W}_{z}}(\widetilde{W}_{z})$  stay in
some compact $K$ of $X$;

\smallskip \sliii $(X,\call)$ has bounded rational cycle geometry.

\smallskip Then there exists a neighborhood $D\ni z^0$ such that $f$
extends meromorphically onto a complete Hartogs figure
$(\widetilde{W},\pi,D)$ after a reparametrization.
\end{lem}
\proof We keep the notations used in the proof of Lemma \ref{area1}.
Only for the annulus $\Sigma_{\alpha_0}$ adjacent to the boundary of
the disc $\Sigma$ the manifold $\calw :=\{f|_{\widetilde{W}_z}: z
\text{ in a neighborhood}$ $\text{ of } z^0\}$ now has dimension
$n=\dim V\ge 2$.

Let $\nu >0$ be the minimum of areas of rational curves tangent to
$\call$ and contained in the compact $K$. We divide $A$ into a
finite union of increasing closed subsets: $A_1\subset ....\subset
A_k\subset ... \subset A_K=A$ where $A_k= \{z\in A:
\area_{\omega}\big(\tilde f|_{\widetilde{W}_z}\big)\le
k\frac{\nu}{2}\}$. For some $k$ the set $A_k\setminus A_{k-1}$ is
thick at origin. In the sequel we take it as $A$. As a result every
converging sequence $\{\tilde f|_{\widetilde{W}_{z_n}}:z_n\in A,
z_n\to z^0\}$ has the same limit. Really two different limits should
differ by a rational cycle. Therefore their areas should differ at
least by $\nu$. Contradiction.

\smallskip  The Banach analytic set $\calc$ obtained literally as
in the proof of Lemma \ref{area1} is again finite dimensional. But
the problem is that the restriction map $r:\calc\to\calw$ might be
not proper. That mean that for some $z$ close to $z^0$ the preimage
$r^{-1}(\tilde f|_{\widetilde{W}_{z}})$ is not compact. But a cycle
in this preimage is different from $\tilde f|_{\widetilde{W}_{z}}$
itself by a rational cycle tangent to $\call$. Therefore we got a
contradiction with the boundedness of rational cycle geometry
condition (iii) of this Lemma. So $r$ is proper and by Remmert
proper mapping theorem  $r(\calc)$ is an  analytic set in $\calw$.
Since it contains a thick subset it is the whole $\calw$. We again
get extension of all $f_z$ for $z$ close to $z^0$ as a family by a
tautological map $\tilde f:\tilde\calw \to X$.

\smallskip

\qed
\begin{rema}\rm
Note that in this Lemma no taming conditions on $(X,\call)$ need to be
imposed (it is weakened, in fact, to a boundedness of the rational cycle
geometry).
\end{rema}

\medskip\noindent{\sl Proof of Theorems \ref{reparalldim1} and
\ref{reparalldim2}.} Let again $\hat U$
be the maximal open subset of $V$ such that $f$ extends onto $(W,\pi
,\hat U)$ after a reparametrization.

\smallskip Now we can proceed exactly as in the proof of Steps 1 and
2 on the page s 817-818 of \cite{Iv6}. Our present situation is even
somewhat simpler because $f$ is holomorphic on the Hartogs figure.
Reparametrizations do not cause any additional difficulties. In this
way we get we get that $S_1\deff\d\hat U\cap V$ is of Hausdorff
dimension $2n-2$. By the ``rationalization trick'' we extend our map
$f$ onto a complete Hartogs figure $(W,\pi ,V)$ over $V$ minus the
closed set of the form $S=\bigcup_{(\lambda,z_1)\in
S_1}S_{\lambda,z_1}$ with $S_{\lambda,z_1}$ being compact subsets of
the discs $\Delta_{\lambda,z_1}$. Due to our localization Lemma
\ref{localization} we need to work in a neighborhood of a point
$(\lambda,z_1)\in S_1$ only. Lemma \ref{area2} shows that for every
natural $N$ the set $\{ z\in V: f|_{W_z} \text{ extends onto a disc
} \widetilde{W}_z \text{ and }\area (\tilde f|_{\widetilde{W}_z})
\le N\}$ is thin in a neighborhood of $S_1$. The rest is obvious. In
particular, one gets that for any two-dimensional submanifold
$U\subset V$ the domain $U\setminus (S_1\cap U)$ is the maximal
domain over which the restricted map $f|_{W|_U}$ extends after a
reparametrization.

\smallskip Starting from this point no further reparametrizations are
needed. Therefore the proof of $(a_1), (a_2)$ and $(b_1) $ is done.

\smallskip $(a_3)$ is clear, because $\tilde f$ could fail to be
immersion only along a divisor, which should then intersect
$\widetilde{W}_U$. But this is not the case.

\smallskip $(b_2)$ is exactly the Corollary 2.2 from the preceding
subsection. Remark that this item  easily implies the following:

\begin{corol}
\label{minimal} In the conditions (and notations) of Theorem
\ref{reparalldim2} denote by $S^0_{\lambda}$ the minimal closed
subset of $\Delta^2_{\lambda}$ such that the restriction
$f_{\lambda}
\deff f|_{\Delta^2_{\lambda}}$ extends onto $\Delta^2_{\lambda}
\setminus S^0_{\lambda}$. Then $S^0_{\lambda}=S_{\lambda}$.
\end{corol}

Item $(b_3)$ follows from Lemma \ref{descrete}.  Theorems \ref{reparalldim1}
and \ref{reparalldim2} are proved.

\smallskip\qed

\smallskip
Let's repeat once more that $S$ in Theorem \ref{reparalldim2} is
always understood as being the minimal closed subset that $f$
meromorphically extends to its complement.

Suppose now that a polydisc $P=\Delta^n\times\Delta_{1+\eps}$ is
fixed, a closed subset $S\subset P$ of the form $S=\bigcup_{z\in
S_1}S_z$ in $P$ is given, where $S_z$ is a compact subset of the
leaf $P_z\deff \{z\}\times \Delta_{1+\eps}$ for every $z\in S_1$.
Suppose that $S_1\ni 0$ and that $0$ is an accumulation point for
$\Delta^n\setminus S_1$. Finally, let a meromorphic foliated generic
injection $f:(P\setminus S,\pi)\to (X,\call)$ into a disc-convex
foliated manifold be given. We shall make use of the following:

\begin{corol}
For a fixed constant $C>0$ let $A_C$ denote a set of $z$ in a
neighborhood of $0\in\Delta^n$ such that $(\lambda,z_1)\not\in S_1$
and that

\begin{equation}
\label{areainfty} \int\limits_{\Delta_{P_z}}f^*\omega \le C.
\end{equation}
If $(X,\call)$ has bounded rational cycle geometry (ex. admits  a
plurinegative taming form) then $A_C$ is contained in a germ of a
proper analytic subset of $\Delta^n$ at $0$.
\end{corol}

The proof follows immediately from Lemmas \ref{area2} and \ref{ratgeo}.

\begin{rema}\rm
\label{exempl1} (a) The Hartogs type extension theorems for meromorphic
maps with values in  manifolds carrying
pluriclosed metric forms were proved in \cite{Iv4} and \cite{Iv6}.
The K\"ahler case was previously done in \cite{Iv3}.

\smallskip\noindent (b) ``Unparametric'' versions of these results
are due to Brunella, see  \cite{Br3} and \cite{Br5}.

\smallskip\noindent (c) Here we prove such type of results for the
forms positive only along a given foliation, which represents
additional difficulties. One of them is that the singularity set $S$
can be ``massive'' along the $z_2$-direction (we are using notations
of Theorem \ref{reparalldim2}). Another is the less obvious
appearance of shells in Lemma \ref{ext-shel}.

\smallskip\noindent (d) Example of Subsection \ref{prg2.3} has all necessary
features to appear in the proofs of this paper. It is a generalized Hartogs
figure and, moreover, it admits a holomorphic foliated immersion, namely
$(g,f):W\to \cc\times\cc\pp^1$, into an algebraic (!) surface.
Therefore the Hartogs figure $(W,\pi, U,V)$ of this example could be
well an open subset of some holonomy covering cylinder
$\hat\call^0_D$ for some foliation $\call$ and, say,
$\gamma_0=\{\frac{2}{3}\}\times \d\Delta$ could be a vanishing
cycle.

\smallskip\noindent (e) It is instructive to see how the mapping
$(g,f):W\to \cc\times\cc\pp^1$ in this example extends to a complete
Hartogs figure after a reparametrization. It is almost tautological:
remark that $f(z_1, \cdot)$ is an imbedding into $\{z_1\}\times \cc
\subset \cc\times\cc\pp^1$ near $\{z_1\}\times \d\Delta$ and
therefore $\Phi = (g,f):\d_0W\to \d_0\widetilde{W}$ is a
biholomorphism for an obvious $\widetilde{W}\subset \cc\times\cc$.
Therefore $\Phi^{-1}$ will be a reparametrization.
\end{rema}

\newsect[sect.VCFS]{Vanishing Cycles, Covering Cylinders and Foliated
Shells.}

\newprg[prgNVFS.van-end]{Vanishing ends and holonomy covering cylinders}
To the classical definitions of holonomy coverings, eg. \cite{Iy1,Iy2,Sz},
we will further employ a subtle extension of \cite{Br3,Br4} to take
account of ``removable singularities''. Let $(X,\call)$ be a foliated
manifold. We use the notations introduced at the beginning of the
Section 2. Take a point $z^0\in X^0$ and denote by $\call^0_{z^0}$
the leaf of $\call^0$ passing through $z^0$. Recall that a parabolic
end of $\call^0_{z^0}$ is a closed subset $E\subset \call^0_{z^0}$
which is biholomorphic to the closed punctured disc
$\bar\Delta^{*}=\{ \zeta\in \cc : 0<|\zeta|\leq 1\}$. By $\d E$ we
shall denote the biholomorphic image of the circle $\{ |\zeta|=1\}$
- the outer boundary of the end $E$. Foliation $\call$ may have a
nontrivial holonomy along $\d E$, which can be finite or infinite.

\smallskip
Consider the case when holonomy is finite. Recall what does that
mean. Take a transversal $D$ through $z^0$, \ie a complex
hypersurface in $X^0$ which is everywhere transverse to the leaves
of $\call$. Transversalis will be always taken small enough, in
particular, we shall always suppose that $D\subset X^0$ and that $D$
is transversal to $\call^0$ ``up to its boundary $\d D$''. Take a
path $\gamma_{z^0}$ on $\call^0_{z^0}$ which goes from $z^0$ to some
point $q\in \d E$, then goes one time around $\d E$ and goes back to
$z^0$. If one takes a point $z\in D$ close to $z$ and travels on
$\call_{z}^0$ along the path $\gamma_{z}$ close to $\gamma_{z^0}$
then one certainly hits $D$ in a neighborhood of $z^0$ by a point
$g(z)$. This defines a local biholomorphism $g: (D,z^0)\to (D,z^0)$
which generates a subgroup $<g>$ of the group $Bihol(D,z^0)$ of
local biholomorphisms of $D$ fixing $z^0$. We suppose that $<g>$ is
finite, \ie $g^d = \id$ for some $d\geq 1$ and this $d$ is always
taken to be the minimal satisfying this property. This $d$ is called
the order of the holonomy of $\call$ along $\d E$.

\begin{lem}
\label{vanishend}
 Let $E\subset \call_{z^0}^0$ be a parabolic end.
Then for a sufficiently small $\eps>0$ there exists a foliated
holomorphic immersion $f:\Delta^{n}\times A_{1-\eps ,1+\eps}\to
\call^0_D$  such that:

\smallskip
\sli $f(\{ 0\}\times A_{1-\eps ,1+\eps})\subset \call^0_{z^0}$  and
the restriction $f|_{\{ 0\}\times A_{1-\eps ,1+\eps}}:\{ 0\}\times
A_{1-\eps ,1+\eps}\to  \call^0_{z^0}$ is a regular covering of order
$d$ (\ie covers $d$-times some imbedded annulus in $\call^0_{z^0}$
and $f(\{0\}\times \d\Delta)=d\cdot \d E$).

\slii For all $z\in \Delta^{n}$ outside of a proper analytic subset
$A\subset \Delta^{n}$ the restriction $f|_{\{ z\}\times A_{1-\eps
,1+\eps}}:\{ z\}\times A_{1-\eps ,1+\eps}\to  \call_{z}$ is an
imbedding.
\end{lem}
\proof It will be convenient to move the point $z^0$ (and the
transversal $D$) to $\d E$. Take an annulus $A_0$ on $\call^0_{z^0}$
around $\d E$. Let $g\in Bihol(z^0,D)$ generates the holonomy of
$\call$ along $\d E$ as above. Denote by $A$ the germ of a proper
analytic subset of $D$ at $z^0$ which consists from those $z\in D$
that the orbit of the corresponding holonomy has cardinality $l>1$
($l$ necessarily divides $d$). When one travels from $z\in D$ to $z$
in the leaf $\call_z^0$ along a curve close to $\d E$ one cuts an
imbedded annulus $A_z$ on $\call_z^0$. For $z$ in the exceptional
set $A$ one sweeps $A_z$ $d/l$ times (if the orbit of $z$ has
cardinality $l$), for $z$ outside from $A$ only once. The union
$W=\bigcup_{z\in D}A_z$ has a natural structure of a complex
manifold and possesses a natural foliated holomorphic immersion $f:
D\times A_{1-\eps,1+\eps}\to W$ coming from the construction, which
is a generic injection (all this provided $D$ is shrinked to a small
neighborhood of $z^0$). $f$ sends each annulus $\{z\}\times
A_{1-\eps,1+\eps}$ onto the corresponding $A_z$ with corresponding
multiplicity. For that one might need to shrink $D$ and annuli $A_z$
for $z\in D$ once more. Now we can suppose that $D$ is biholomorphic
to $\Delta^{n}$. The rest is obvious.

\smallskip\qed

As we remarked in the proof our  $f$ is a {\sl generic injection} of
the trivial Hartogs figure $\Delta^{n}\times A_{1-\eps, 1+\eps}$
over a polydisc in the sense of Definition \ref{geninject} and
results of the previous section are applicable to such $f$.

\begin{defi}
\label{van-end-def1} A parabolic end $E$ is called a {\slsf
vanishing end} of order $d$ if:

\sli the holonomy of $\call$ along $\d E$ is finite of order $d\geq
1$;

\slii the generic injection $f:\Delta^{n}\times A_{1-\eps
,1+\eps}\to \call^0_D$, constructed above,  extends as a foliated
meromorphic immersion $\tilde f:\widetilde{W} \to X$ from a complete
Hartogs figure $(\widetilde{W},\pi,\Delta^n)$ over $\Delta^{n}$ to
$X$ after a reparametrization;

\sliii the intersection of $\widetilde{W}_0\deff \pi^{-1}(0)$ with
the set of points of indeterminacy $I_{\tilde f}$ of $\tilde f$
consists of a single point $a\in \widetilde{W}_0$.
\end{defi}
The point $q=\tilde f|_{\widetilde{W}_0}(a)$ will be called the
{\slsf endpoint} of the vanishing end $E$ (or of the leaf
$\call^0_z$). Following Brunella, see \cite{Br3}, we add all
vanishing endpoints to the leaf  $\call^0_{z^0}$ and call the curve
obtained a {\slsf completed leaf} through $z^0$. Completed leaf will
be denoted as $\call_{z^0}$.

\begin{rema}\rm
\label{van-end-rem1} Let us give two very simple examples explaining
this notion.

\smallskip\noindent{1.} Consider the radial foliation in $\cc^2$,
\ie $\call_c=\{z_2/z_1=c\}$ for $c\in\cc\pp^1$. The origin of
$\cc^2$ is a parabolic end for every leaf $\call_c^0$. But it is
never a vanishing end! Really, one cannot construct a foliated
meromorphic immersion as in Definition \ref{van-end-def1} in this
case. Any $\tilde f$ will contract some complex curve to a point.

\smallskip\noindent{2.} Let $\call^v$ be the vertical foliation in
$\cc^2$, \ie $\call_c=\{z_1=c\}$ for $c\in \cc$. Blow-up the origin
$\pi :\hat\cc^2\to\cc^2$ and lift our foliation to $\hat\cc^2$. The leaf
$\call_0^0$ has now a parabolic end at its point of intersection with
the exceptional divisor and this end is a vanishing end. The role
of $f=\tilde f$ plays $\pi^{-1}$.

\smallskip\noindent{3.} Let $D$ be a transversal through $z$ and let
$E$ be a vanishing end of $\call_z^0$. Remark that for points $z'$
close to $z$ on $D$ only those ones which belong to  some proper
analytic subset could be such that $\call_{z'}^0$ has a vanishing
end $E'$ with $\d E'$ close to $E$. Really, such $z'$ should lie in
the projection $\tilde \pi :\widetilde{W}\to D$ of a point of
indeterminacy of $\tilde f$. Therefore a generic leaf of $\call$ has
no vanishing ends at all.

\smallskip\noindent{4.} If the holonomy along $\d E$ is infinite then
$E$ is never a vanishing end by definition.
\end{rema}

\begin{rema} \rm
\label{exempl2} Example, which we discussed in the previous
Section, showed the necessity of modifying the notion of a vanishing
end, see discussion before the Theorem 3.1  in \cite{Br4}. The necessary
modification  was undertaken in \cite{Br3} and it is this notion which we
use along this paper.
\end{rema}

\smallskip
For each $z\in D$  take a holonomy cover $\hat\call_z^0$ of the leaf
$\call_z^0$. Recall that a holonomy cover of $\call_z^0$ is a cover
with respect to the holonomy subgroup $\hol(z,\call_z^0)$ of the
fundamental group $\pi (z,\call_z^0)$. That means that in the
construction of $\hat\call_z^0$ two pathes $\gamma_1, \gamma_2$ from
$z$ to some $w\in\call_z^0$ define the same point of $\hat\call_z^0$
if and only if $\gamma_1\circ\gamma_2^{-1}\in \hol(z,\call_z^0)$,
\ie if the holonomy along $\gamma_1\circ\gamma_2^{-1}$ is trivial.

\smallskip Set

\begin{equation}
\eqqno(holonomcylinder) \hat\call^0_D = \bigcup_{z\in
D}\hat\call_z^0.
\end{equation}
This set (introduced by Suzuki in \cite{Sz} under the name of ``tube
normaux'') has the natural structure of a complex manifold together
with the natural projection $\pi :\hat\call^0_D\to D$ which sends
$\hat\call_z^0$ to $z$. It admits also the natural locally
biholomorphic foliated map $p:\hat\call^0_D \to \call^0_D\subset
X^0$ which sends $\hat\call_z^0$ to $\call_z^0$ with
$p|_{\hat\call_z^0}:\hat\call^0_z\to \call^0_z$ being the canonical
holonomy covering map. Call $\hat\call^0_D $ the {\slsf holonomy
covering cylinder} of $\call$ over $D$.

\smallskip
Vanishing ends of $\hat\call_z^0$ are defined similarly to that of
$\call_z^0$. Let $E$ be a parabolic end of $\hat\call_{z^0}^0$ Take
$f:\Delta^{n}\times A_{1-\eps,1+\eps}\to \hat\call^0_D$ such that:

\medskip
\sli $f:\Delta^{n}\times A_{1-\eps,1+\eps}\to \hat\call^0_D$ is an
imbedding;

\slii  $f(\{0\}\times \d\Delta)=\d E$ (note that $d=1$ in this case).

\medskip The only difference that now $f$ takes values in $\hat\call_D^0$
and $f$ is an imbedding. The last is because the holonomy of the
foliation $\hat\call^0$ on $\hat\call_D^0$ is trivial.

\begin{defi}
\label{van-end-def2} $E$ is called a vanishing end of
$\hat\call_{z^0}^0$ if $h=p\circ f$ extends to a meromorphic
foliated immersion $\tilde h:\widetilde{W}\to X$ after a
reparametrization (not $f$ itself as in Definition
\ref{van-end-def1}) and $\widetilde{W}_0$ intersects the
indeterminacy set $I_{\tilde h}$ of $\tilde h$ by exactly one point.
\end{defi}
The union of $\hat\call^0_z$ with all its vanishing endpoints
equipped with an obvious complex structure  will be denoted as
$\hat\call_z$. We shall call it also {\slsf a completed holonomy
covering leaf} of the leaf $\call^0_z$. Set
$\hat\call_D:=\bigcup_{z\in D} \hat\call_z$ and call it the {\slsf
completed holonomy covering cylinder} over $D$. Now let us bring
together the principal properties of $\hat\call_D$, which will be
repeatedly used along this paper.

\begin{lem}
\label{imbed}
\sli The completed holonomy covering cylinder possesses the natural
structure of a foliated complex manifold with foliation given by the
natural projection $\pi :\hat\call_D\to D$ defined as above by $\pi
(\hat\call_z) = z$.

\slii The natural foliated holomorphic immersion $p:\hat\call^0_D\to
\call^0_D$ extends to a meromorphic foliated immersion $p
:\hat\call_D\to X$ and its restrictions
$p|_{\hat\call_z}:\hat\call_z\to \call_z$ are ramified at vanishing
ends.
\end{lem}
\proof (\sli Cylinder $\hat\call^0_D$ has a natural complex
structure. Therefore we need to add vanishing ends to some leaves
and extend this structure to a neighborhood of each added end. Take
a vanishing endpoint $a\in\hat\call_{z^0}$. Let $f:\Delta^{n}\times
A_{1-\eps,1+\eps}\to \hat\call_D^0$ be an imbedding from the
Definition \ref{van-end-def2} with $h=p\circ f$ already extended to
a meromorphic foliated immersion of $(\Delta^{n}\times
\Delta_{1+\eps},\call^{\v})$ into $(X,\call)$. Let $I_h$ be the
indeterminacy set of $h$. For $z\not\in A:=\pi (I_h)$ the
restriction
$h|_{\{z\}\times\Delta_{1+\eps}}:\{z\}\times\Delta_{1+\eps}\to X$ is
an imbedding and therefore so is also the
$f|_{\{z\}\times\Delta_{1+\eps}}:\{z\}\times\Delta_{1+\eps} \to
\hat\call_D^0$. This implies that $f$ is an imbedding on
$\Delta^n\setminus (A\times\Delta )$. This immediately implies that
$f$ is an imbedding on $\Delta^n\setminus I_h$. Therefore we can
complete $\hat\call^0_D$ by $I_h$ over the image $(\pi\circ
f)(\Delta^n) \subset D$. This defines the structure of a complex
manifold on $\hat\call_D$. The rest is obvious.

\smallskip\noindent (\slii This item follows readily from the
construction above.

\smallskip\qed

\begin{rema}\rm
\label{orbifold1}
Let us make a remark which will be important for
the future. The covering $p_{z^0}:\hat\call_{z^0}\to\call_{z^0}$ is
an orbifold covering. That means that its ramification index at
point $a$ depends only on $b\deff p_{z^0}(a)$. This is also an
unbounded covering in the sense that for every $a$ there exists a
disc-neighborhood $V\ni b$ such that $p_{z^0}^{-1}(V)$ is a disjoint
union of discs $W_j$ with centers $a_j$ - preimages of $b$, such
that every restriction $p_{z^0}|_{W_j}:W_j\to V$ is a proper
covering ramified over $b$.
\end{rema}

\newprg[prgVCFS.exemp]{Vanishing ends in dimension 3 vs dimension
2: an example.} We want to exploit the example constructed in
\cite{Iv5}. This example of a compact (!) complex threefold $X$ has
the following very strange features:

\smallskip\sli{\it First: for every domain $D\subset\cc^2$ every
meromorphic map $f:D\to X$ extends to a meromorphic map $\hat f:\hat
D\to X$ of the envelope of holomorphy $\hat D$ of $D$ into $X$.} For
example, if $D=H_{\eps}$ is the standard Hartogs figure then every
$f$ extends from $H_{\eps}$ to the bidisc $\Delta^2_{1+\eps}$. One
can easily prove (following the lines in \cite{Iv5} on the pp.
99-105) also a non-parametric version of this statement in the
spirit of Theorem \ref{repar2dim} of the present paper.

\smallskip\slii{\it Second: but there exists a meromorphic map
$f:\Delta^3_*\to X$ of the punctured $3$-disc into $X$ which doesn't
extend to the origin.}

\smallskip The construction goes (very roughly) as follows.

\smallskip\noindent a) Take the standard three-ball
$\bb^3\subset\cc^3$ and blow-up the line $l_0=\{z_2=z_3=0\}$ in it.
Denote by $E_0$ the exceptional divisor and by $Y_0$ the threefold
obtained.

\smallskip\noindent b) Denote by $l_1$ the intersection of the
exceptional divisor with the proper transform of the plane $L_1\deff
\{ z_1=0\}$. Proper transforms will be denoted with the same
letters. Therefore $l_1=L_1\cap E_0$. This is a copy of a projective
line $\pp^1$. Now blow-up $l_1$. Denote by $E_1$ the exceptional
divisor of this second blow-up.

\smallskip\noindent c) Take some $\pp^1$ on $E_1$ (for example the
intersection of $E_1$ with the proper transform of the plane
$L_2\deff \{z_2=0\}$), denote is as $l_2$ and blow it up again.
$E_2$ is again the exceptional divisor and by $l_3$ denote the
projective line, which is the intersection of $E_2$ with (the proper
transform of) $E_0$.

\smallskip\noindent d) Fix a point $0_3$ on $l_3$ and consider the
biholomorphism $g$ of $\bb^3\subset\cc^3$ onto a neighborhood of
$0_3$ which sends the origin $0_0$ of the space $\cc^3$ to $0_3$ and
which, moreover, in the natural coordinates of $\cc^3$ and of the
resulting blowing-up is the identity. If one writes $g$ in the
natural coordinates $(z_1,z_2,z_3)$ of $_cc^3$ only then it has the
form
\[
g : (z_1,z_2,z_3) \to (z_1,z_2,z_2z_2^2z_3,z_1,z_2,z_3),
\]
see \cite{Iv5} for much more details.

\smallskip\noindent e) Blow-up $l_3$ and denote the resulting
threefold as $Y_3$. $g$ lifts to an imbedding of $Y_0$ into $Y_3$,
denote this lift as $\hat g$. Remove $\hat g(Y_0)$ from $Y_3$, \ie
consider the threefold $X_0\deff Y_3\setminus \hat g(Y_0)$. Identify
the boundary components of $X_0$ by $\hat g$ and obtain a compact
complex threefold $X$. This $X$ is our example.

\smallskip Universal covering $\tilde X$ of $X$ can be obtained by
gluing infinitely many copies $\{X_0^j\}_{j=-\infty}^{+\infty}$ of
$X_0$ one to another by $\hat g$ (more precisely by corresponding
powers of $\hat g$, \ie $X_0^{-1}$ is attached to $X_0^0\deff X_0$
by $\hat g$ and $X_0^1$ is attached to $X_0^0$ by  $\hat g^{-1}$ and
so on). Remark that $\tilde X$ can be naturally viewed as a
blowed-up $\cc^3\setminus \{0\}$. For us the following feature of
$X$ will be sufficient:

\smallskip\sliii {\it if one takes a (singular in general) complex
surface $Z$ in $\bb^3$ and lifts $Z\setminus \{0\}$ naturally to
$\tilde X$ then the closure $\tilde Z$ of this lift is contained in
$\bigcup_{j=-N}^NX_0^j$ for some $N$, \ie in some "finite" part of
the universal covering $\tilde X$.}

\smallskip For example, it is not difficult to see that the closure
of the lift of $(L_1\cap \bb^3)\setminus \{0\}$ is a three times
blowed-up two-ball imbedded into $X_0^{-1}\cup X_0^0\cup X_0^1$. In
fact (\sliii explains why (\slii is true.

\smallskip Take a vertical foliation $\call\deff \{z_1=\const,
z_2=\const\}$ on $\bb^3$. It lifts to a foliation (denoted with the
same letter) on $\tilde X$. Take the transversal $D\deff \{ z_3=1/2,
|z_1|^2 + |z_2|^2 < 1/2\}$. The leaf $\call_0$ has an obvious end -
the image of the punctured disc $\{(0,0,z_3) : 0<|z_3|<\frac{1}{2}\}$
in $X$. This end is {\slsf not} a vanishing end, because
$\tilde\call_D$ has a puncture at this end.  At the same time if one
takes any one-disc $S\subset D$ then $\tilde\call_S$ lifts to some
finite part
$\bigcup_{j=-N}^NX_0^j$ of $\tilde X$ ($N$ crucially depends on the
choice of $S$ and there is no any bound on it). Therefore our end
will be a vanishing end for any $\tilde\call_S$ as above.

\begin{rema}\rm
This example explains how subtle is the definitions of a vanishing end
when working in the manifolds of higher dimension.
\end{rema}

\newprg[prgVCFS.van-cyc]{Vanishing cycles}

\smallskip Let now $\hat\gamma :[0,1]\to \hat\call_{z}^0$ be a loop
in $\hat\call_{z}^0$ which is not homotopic to zero in
$\hat\call_{z}^0$.

\begin{defi}
\label{vanishcycle} We say that $\hat\gamma $ is a vanishing cycle
if for some sequence $z_n\to z$ there exist loops $\hat\gamma_n$ in
$\hat\call_{z_n}$ uniformly converging to $\hat\gamma$ which are
homotopic to zero in the corresponding leaves $\hat\call_{z_n}$.

(a) We say that $\hat\gamma$ is an algebraic vanishing cycle if
$\gamma$ is not homotopic to zero in $\hat\call_{z}^0$ but is
homotopic to zero in the completed leaf $\hat\call_{z}$.

\smallskip (b) If $\hat\gamma$ is not homotopic to zero also in the
completed leaf $\hat\call_{z}$ we call it an {\slsf essential
vanishing cycle}.
\end{defi}
There is an analogy (rather deep in fact) between {\slsf
algebraic}/{\slsf essential} vanishing cycles and poles/essential
singularities of meromorphic functions. Really, pole of a
meromorphic function $f$ becomes a regular point if one completes
$\cc$ to $\cc\pp^1$ and considers $f$ as a holomorphic mapping into
the latter manifold. However, an essential singular point stays to
be a singularity of $f$ also after this operation. The same with
cycles. For the moment let us say that:

\begin{itemize}
\item If $\call^{\sing}=\emptyset$, \ie if $\call$ has no singularities,
then every vanishing cycle is an essential vanishing cycle, more
precisely projects to a vanishing cycle under the holonomy covering
map $\hat\call_z\to\call_z$, see Remark \ref{equiv}.

\item Algebraic vanishing cycles in the leaf $\hat\call_z^0$ can be
removed (\ie one can make these cycles homotopic to zero) by adding
to $\hat\call_z^0$   {\sl vanishing ends}.

\item It is known also (it follows from \cite{Br3}) that if $X$ is
K\"ahler, then all vanishing cycles (of any $\call$) are algebraic.
\end{itemize}

\smallskip In this paper we shall concentrate our attention on
essential vanishing cycles only. In this subsection, following
\cite{Br1}, we show that if $\hat\call_z$ contains an essential
vanishing cycle then it contains an {\slsf imbedded} essential
vanishing cycle. Take an immersed loop $\gamma$ in a Riemann surface
$R$ which has only transversal self-intersections. Denote by $N$ the
closure of a sufficiently small tubular neighborhood of $\gamma$.
Add to $N$ all discs bounded by circles - components of $\d N$, and
denote the obtained compact as $\bar N$.

\begin{lem}
\label{pi1-inject}
Imbedding $\overline{N}\subset R$ induces the
natural injection $\pi_1(\bar N)\to \pi_1(R)$.
\end{lem}
\proof Suppose that there exists a loop $\beta$ in $\bar N$ not
homotopic to zero in $\bar N$ which is homotopic to zero in $R$.
Then the homotopy of $\beta$ to zero is supported in a compact part
of $R$ and therefore we can suppose that $R$ has finite topology,
\ie finite number of handles and boundary circles. In the sequel the
trivial case when $\bar N$ or $R$  is a disc or an annulus will be
omitted. Now we perform the following manipulations which obviously
do not change the homotopy type of $\bar N$.

\smallskip (a) Every connected component of  $R'\deff R\setminus \bar N$
which is an annulus adjacent to $\d R$ we add to $\bar N$.

\smallskip (b) If some component $C$ of $R'$ is an annulus with both
boundary circles belonging to $\d\bar N$ then we cut $C$ on two
annuli $C_1$ and $C_2$. Each of them we add to $\bar N$ and think
about $\bar N$ as having $\d C_1$ and $\d C_2=-\d C_1$ as two
boundary components.

\smallskip
Denote by $g$ the Riemannian metric on $R$  of curvature $-1$ having
boundary circles as geodesics. Every loop $\gamma$ in $R$ is now
homotopic to a unique geodesic $\tilde \gamma$ in metric $g$ which
is either not intersecting $\d R$ or is a component of $\d R$, see
for example \cite{Bu} Theorem 1.6.6. We deform all boundary circles
of $\bar N$ one by one to geodesics. If in the process of
deformation a curve is touching $\beta$ we move $\beta$
appropriately enlarging (or contracting) $\bar N$ in a way to keep
$\beta$ inside.

We end up with having all boundary circles of $\bar N$ geodesics in
$g$. Now we do the same with $\beta$ getting from it a geodesic
$\tilde \beta$ in $\bar N$. Note that it stays in $\bar N$ and do
not intersect also $\d C_2=-\d C_1$ from (b) (or coinciding with one
of them). But this But $\tilde\beta$ stays to be geodesic in $g$ on
the whole of $R$ and therefore is not homotopic to zero.
Contradiction.

\smallskip\qed

\medskip Now we are going to reduce the question of existence of
essential vanishing cycles in  $\hat\call_D$ to the existence of
{\slsf imbedded} essential vanishing cycles in $\hat\call_D$.
Namely, we shall prove that the following is true:

\begin{lem}
\label{imbvancyc1} If there exists an essential vanishing cycle in
$\hat\call_{z^0}$ then there exists an imbedded essential vanishing
cycle in $\hat\call_{z^0}$.
\end{lem}
\proof Let $\hat\gamma_0:[0,1]\to \hat\call_{z^0}$ be our essential
vanishing cycle. After perturbing it, if necessary, we can suppose
that $\hat\gamma_0$ is an immersion with only transversal
self-intersections. For every point $\hat\gamma_0(t)$ take an
$(n-1)$-disc $Q_{\hat\gamma_0(t)}$ in $\hat\call_D$ transversal to
the leaf $\hat\call_{z^0}$ and cutting it by the point
$\hat\gamma_0(t)$. Make these discs depend smoothly on
$\hat\gamma_0(t)$ in such a way that for
$\hat\gamma_0(t_1)\not=\gamma_0(t_2)$ the corresponding
$(n-1)$-discs do not intersect. Let's stress explicitly that
$Q_{\hat\gamma_0(t)}$ depends only on the image point
$\hat\gamma_0(t)$ on the curve and not on $t$. We have therefore a
natural projection $\Pi :
\bigcup_{\hat\gamma_0(t)}Q_{\hat\gamma_0(t)}\to \hat\gamma_0(t)$.
Extend these data over a closure of a small tubular neighborhood
$N_0$ of $\hat\gamma_0$. I.e., set $Q:=\bigcup_{\tau\in
N_0}Q_{\tau}$ and now $\Pi$ maps this $Q$ onto $N_0$.

\smallskip For every $z$ in our transversal $D$, which is close to
$z^0$ each $Q_{\tau}$ cuts the leaf $\call_z^0$ exactly by one point
and when $\tau$ runs over $N_0$ our discs $Q_{\tau}$ cuts a closure
of a tubular neighborhood $N_z$ of some closed curve $\gamma_z$
which covers $\hat\gamma_0$ under the projection
$\Pi|_{\gamma_z}:\gamma_z\to \hat\gamma_0$. Remark also that
$\Pi|_{N_z}:N_z\to N_0$ is bijective. Denote by $\bar N_0$ the union
of $N_0$ with all discs bounded by circles components of $\d N_0$.
Denote likewise by $\bar N_z$ the union of $N_z$ with all discs
bounded by circles components of $\d N_z$.

\smallskip Take some component $\hat\gamma_0^{'}$ of $\d N_0$ bounding
a disc in $\hat\call_{z^0}$. Then the corresponding component
$\gamma_{z_n}^{'}$ of $\d N_{z_n}$ bounds a disc in
$\hat\call_{z_n}$, say $D_0^{'}$ and then $\Pi|_{N_{z_n}}:N_{z_n}\to
N_0$ extend to a homeomorphism $\Pi|_{N_{z_n}\cup
D_{z_n}^{'}}:N_{z_n}\cup D_{z_n}^{'}\to N_0\cup D_0^{'}$.

\smallskip If $\gamma_{0}^{'}$ doesn't bound a disc in $\hat\call_{z^0}$
but $\gamma_{z_n}^{'}$ do  bounds a disc in $\hat\call_z$  we get an
imbedded essential vanishing cycle in $\hat\call_{z^0}$.

\smallskip
So, unless an imbedded vanishing cycle was found in
$\hat\call_{z^0}$ we end up with extending $\Pi$ to a homeomorphism
$\tilde\Pi :\bar N_{z_n}\to \bar N_0$.

\smallskip Since $\gamma_{z_n}$ is homotopic to zero in $\hat\call_{z_n}$
it will be homotopic to zero in $\bar N_{z_n}$ by Lemma 3.1.
Therefore $\hat\gamma_0$ should be homotopic to zero in $\bar N_0$
and therefore in $\hat\call_{z^0}$. Contradiction. Therefore the
only possibility left is that some component $\gamma_0^{'}$ of $\d
N_0$ doesn't bound a disc while $\gamma_{z_n}^{'}$ do bound disc,
\ie $\gamma_0^{'}$ is an imbedded essential vanishing cycle in
$\hat\call_{z^0}$.

\smallskip\qed

\begin{rema}\rm
Remark that if $\hat\gamma_0$ is an imbedded essential vanishing
cycle in $\hat\call_{z^0}$ then a sequence $z_n\to z$ such that
there exists $\gamma_{z_n}$ bounding a disc in $\hat\call_{z_n}$ and
$\gamma_{z_n}$ uniformly converging to $\gamma_0$ when $z_n\to z^0$
can be taken generic.
\end{rema}

\newprg[prgVCFS.cov-cyl]{Universal covering cylinder}

\smallskip
Further, for $z\in D$ denote by $\tilde\call_z$ the universal cover
of the completed holonomy leaf $\hat\call_z$. I.e., we take the
orbifold universal covering of $\call_z$, see Remark
\ref{orbifold1}. On the union
\begin{equation}
\eqqno(covcylinder)
\tilde\call_D=\bigcup_{z\in D}\tilde\call_z
\end{equation}
one defines a natural topology in the following way. An element of
$\tilde\call_D$ is a path $\gamma $ in some leaf $\hat\call_z$
starting from $z$ and ending at some point $w\in \hat\call_z$.
$\gamma $ and $\gamma^{'}$ define the same point if their ends
coincide and they are homotopic (inside $\hat\call_z$) with ends
fixed. A neighborhood of $\gamma\subset\hat\call_z$ in
$\tilde\call_D$ is the set of pathes $\gamma'$-s in the leaves
$\hat\call_{z'}$ with $z'$ close to $z$ which are themselves close
to $\gamma$. $\gamma'$ ``close'' to $\gamma$ is understood here as
closed in the topology of uniform convergence in the space $\calc
([0,1],X)$ of continuous mappings from $[0,1]$ to $X$.

\begin{defi}
$\tilde\call_D$ with the topology just described is called the
universal covering
cylinder of $\call$ over $D$.
\end{defi}

The natural projection $\pi :\hat\call_D\to D$ lifts to $\pi
:\tilde\call_D\to D$ (and will be denoted with the same letter).
There is a distinguished section $\sigma : D\to \tilde\call_D$
sending $z$ to $z$. The mapping $p:\hat\call_D\to X$ lifts to
$\tilde\call_D$ and stays to be a {\sl meromorphic foliated
immersion} $\tilde p:\tilde\call_D\to X$ in the sense that it is a
foliated immersion outside of its indeterminacy set.

\medskip Due to the eventual presence of essential vanishing cycles
the natural topology on the covering cylinder might be not
Hausdorff. Let us explain this in more details. Non-separability of
the natural topology on $\tilde\call_D$ means that:
\begin{itemize}
\item  there exist $z\in D$ and $w\in \hat\call_z$ and two pathes
$\gamma_1,\gamma_2$ from $z$ to $w$ such that
$\gamma_1\circ\gamma_2^{-1}$ is not homotopic to zero in
$\hat\call_z$;

\item there exist some sequence $z_n\to z$ in $D$, some sequence
$w_n\in \hat\call_{z_n}$ converging to $w$, some sequences of pathes
$\gamma_1^n$ and $\gamma_2^n$ from $z_n$ to $w_n$ each converging
uniformly to $\gamma_1$ and $\gamma_2$ such that  $\gamma_1^n\circ
(\gamma_2^n)^{-1}$ are homotopic to zero in $\call_{z_n}$.

\end{itemize}
And that exactly  means that $\gamma_1\circ\gamma_2^{-1}$ is an
essential vanishing cycle. Vice verse, if $\gamma :[0,1]\to
\hat\call_z$ is an essential vanishing cycle starting and ending at
$z$, then $\gamma$ and the trivial path $\beta \equiv z$ represent
two non-separable points in $\tilde\call_D$.

\newprg[prgVCFS.proof1]{Proof of Theorem 1}

In the presence of a pluriclosed taming form the problem of the
separability of the topology of $\tilde\call_D$ can be resolved by
Theorem 1 from the Introduction. Now we shall state and prove
somewhat more general and precise statement which contains the
aforementioned result. To make the statement more precise let's turn
to the definition of a {\slsf foliated shell}, \ie to the Definition
\ref{fol-shell1} from the Introduction.
\begin{defi}
\label{fol-shell2} In general, when foliation $\call$ is singular we
require that the mapping $h:(B^{\eps},\call^{\v})\to (X,\call)$
which defines a foliated shell takes its values in $X^0$.
\end{defi}
By the Theorem \ref{reparalldim2} we know that mapping
$h:(B^{\eps},\call^{\v})\to (X^0,\call)$, which defines a foliated
shell in a pluritamed foliated manifold extends onto
$P^{\eps}\setminus \bigcup_{z_1\in S_1}S_{z_1}$, where $S_1$ is at
most countable compact in $\Delta$. One more remark: for a
transversal $D\subset X^0$ and an imbedded disc $\Delta\subset D$
the restriction $\call_{\Delta}^0\deff\bigcup_{z\in\Delta}\call_z^0$
is well defined (we don't need to give this set more structure then
this which it already has).

\begin{thm}
\label{immshell} Let $(X,\call)$ be a disc-convex foliated manifold
which admits a $dd^c$-closed taming form and let $z^0\in X^0$ be a
point. Then the following statements are equivalent:

\sli The leaf $\hat\call_{z^0}$ contains an essential vanishing
cycle.

\slii For every transversal $D\ni z^0$ there exists an imbedded disc
$z^0\in \Delta\subset D$ such that $\call^0_{\Delta}$ contains a
foliated shell.
\end{thm}

\begin{rema}\rm
\label{equiv} (a) Let us explain that the item (i) of this Theorem
is equivalent to the item (i) of Theorem 1 from the Introduction in
the case when $\call$ is smooth (\ie without singularities). In that
case vanishing ends do not exist and, in particular,
$p_{z^0}:\hat\call_{z^0}\to\call_{z^0}$ is an unramified covering.
Let $\gamma_0\subset \call_{z^0}$ be a vanishing cycle and
$\gamma_n\subset \call_{z_n}$ be cycles homotopic to zero and
converging to $\gamma_0$. All $\gamma_n$ lift to cycles
$\hat\gamma_n\subset \hat\call_{z_n}$ converging to the lift
$\hat\gamma_0\subset\hat\call_{z^0}$ of $\gamma_0$. All
$\hat\gamma_n$ are homotopic to zero. But $\hat\gamma_0$ cannot be
homotopic to zero. Therefore we get a vanishing cycle $\hat\gamma_0$
in $\hat\call_{z^0}$.  Vice verse, let $\hat\gamma_0$ and
$\hat\gamma_n$ be as above in the holonomy covering leaves. Then
$\hat\gamma_n$ project to cycles homotopic to zero in corresponding
leaves. But $\hat\gamma_0$ project to some $\gamma_0$ which cannot
be homotopic to zero because in the latter case its lift
$\hat\gamma_0$ (as lift of any curve homotopic to zero) should be
homotopic to zero itself. Therefore $\gamma_0$ is a vanishing cycle
in $\call_{z^0}$.

\smallskip\noindent (b) The item (ii) specifies that the ``support''
$\Sigma = h(B)$ of the foliated shell is in $\call_{\Delta}^0$ (but
it is not homologous to zero in the whole of $X$!) Remark also that
the existence of an essential vanishing cycle in $\hat\call_{z^0}$
is unrelated to the choice of a transversal $D\ni z^0$ (and also on
the imbedded disc $z^0\in\Delta\subset D$). Therefore if for some
transversal and disc in it $D\supset \Delta\ni z^0$ there is a shell
in $\call_{\Delta}^0$ then it persists in all others.

\end{rema}
\proof
\smallskip\noindent (\sli $\Rightarrow$ (\slii
By Lemma \ref{imbvancyc1} we can suppose that our vanishing cycle
$\hat\gamma_0$ is imbedded into $\hat\call_{z^0}\subset\hat\call_D$.
Deforming it, if necessary, we suppose that $\hat\gamma_0$ is
contained in $\hat\call^0_{z^0}$. Therefore we can suppose that for
an imbedded loop $\hat\gamma_0\subset \hat\call_{z^0}^0$ started at
$z^0$ the following holds:
\begin{itemize}
\item $\hat\gamma_0$ does not bound a disc in $\hat\call_{z^0}$ ;

\item but for a generic sequence $z_n\to z^0$ and a sequence of imbedded
loops $\hat\gamma_n\subset\hat\call_{z_n}$ uniformly converging to
$\hat\gamma_0$ every  $\hat\gamma_{n}$ bounds a disc $\Delta_{n}$ in
$\hat\call_{z_n}$.
\end{itemize}

\smallskip Take a neighborhood $U$ of some $z_{N}$ such that for every
$z\in U$ there is an imbedded loop $\hat\gamma_z$ close to
$\hat\gamma_{N}$ bounding a disc $\Delta_z$ in $\hat\call_z$. We can
suppose that $\hat\gamma_z$ smoothly depend on $z\in U$. Take some
open cell $V\subset D$ containing $U$ and $z^0$ and extend our
family $\Gamma :=\{\hat\gamma_z\}$ smoothly over $z\in V$ (after
shrinking it over $U$, if necessary) in such a way that
$\hat\gamma_{z^0}$ coincides with $\hat\gamma_0 $. Perturbing the
family $\Gamma $, if necessary, we can suppose that some
neighborhood $W$ of $\Gamma\cup \Delta_{z_{N}}$ in $\hat\call_D$
forms a generalized Hartogs figure $(W,\pi ,U,V)$. Projection $\pi
:W\to V$ here is the restriction to $W$ of the natural projection
$\pi :\hat\call_D\to D$.

\smallskip Mapping $p:\hat\call_D\to X$ restricted to $W$ will be
likewise denoted as $p :W\to X^0 \subset X$ and it is a holomorphic
foliated immersion, because  the  construction can be obviously
fulfilled in such a way that $W\subset\hat\call_D^0$. Note also that
$p$ is a generic injection because for generic $z_{N}$ our
$p|_{\gamma_N}$ is an imbedding. But $p|_{\hat\gamma_0 }$ might be
only an immersion in general. By Theorems \ref{reparalldim1} and
\ref{reparalldim2} $p$ extends after a reparametrization onto
$\widetilde{W}\setminus S$, where $\widetilde{W}$ is a complete
Hartogs figure over $V$ and $S$ is of the form $S=\bigcup_{z\in
S_1}S_z$ with $S_1$ being a complete $(n-1)$-polar subset of $V$ and
every $S_z$ is a compact subdisc of the corresponding disc
$\widetilde{W}_z$. This extension $\tilde p$ is a foliated
meromorphic immersion, \ie it is an immersion outside of its
indeterminacy set $I_{\tilde p}$ and takes values in $X$ (not more
in $X^0$). The family, which corresponds in $\widetilde{W}$ to our
family $\Gamma$ will be denoted still by $\Gamma$ and no new
notation for the loops $\hat\gamma_z$ will be introduced.

\smallskip Observe that $z^0\in S_1$. Otherwise take an $n$-disc
$\Delta^{n}$ around $z^0$ in $D$ such that $\Delta^{n}\cap
S_1=\emptyset$ and such that:
\begin{itemize}
\item $\pi^{-1}(\Delta^{n})$ is biholomorphic to $\Delta^{n}\times
\Delta_{1+\eps}$ with $\pi$ being the vertical projection
$\Delta^{n}\times \Delta \to \Delta^{n}$ (one might need to shrink
$\Delta^{n}$ and $\widetilde{W}$ to achieve this).

\item For $z\in\Delta^{n}$ circles $\hat\gamma_z=\d\Delta_z$ belong to
our family $\Gamma$ (for this one might need to perturb $\Gamma$).
\end{itemize}
Our $\tilde p$ now is meromorphically extended to $\Delta^{n}\times
\Delta_{1+\eps}$. But that means (by the very definition of
vanishing ends) that $p^{-1}\circ\tilde p$ lifts to a holomorphic
map $\tilde f:\Delta^{n}\times \Delta_{1+\eps}\to\hat\call_D$.
Therefore $\tilde f|_{\{z^0\}\times \bar\Delta}$ realizes the
homotopy of $\hat\gamma_0=\tilde f(\{z^0\}\times \d\Delta)$ to zero.
Contradiction.

\smallskip Denote by $A$ the proper analytic subset in a neighborhood
of $z^0$ on $D$ which consists from points $z$ such that $\tilde
p|_{\d_0\widetilde{W}_z}$ is not a generic injection. Again we
locally represent $\widetilde{W}$ as a product
$\widetilde{W}=\Delta^{n}\times \Delta_{1+\eps}$ with $\pi$ being
the vertical projection $\pi : \Delta^{n}\times \Delta_{1+\eps}\to
\Delta^{n}$ and with $z^0$ being the origin in these coordinates.
Decompose $\Delta^{n}=\Delta^{n-1}\times\Delta$ in such a way that
$(\{0\}\times\Delta)\cap A=\{ \text{ a finite set } \}$. Then by
Theorem \ref{reparalldim2} for $\lambda\in\Delta^{n-1}$ close to $0$
(if $n\ge 2$), or equal to $0$ (if $n=1$) we have $\tilde
p(\d\Delta^2_{\lambda})\not\sim 0$ (not homologous to zero in $X$).
Moreover, for every $z=(\lambda,z_1)$ in $\Delta^{n}$ such that
$(\lambda,z_1)\not\in A$ we have that $\tilde p|_{\{z\}\times
\d\Delta}$ is an imbedding.

\smallskip  Therefore we get a foliated shell  in $\call_{\Delta}^0$,
where $\Delta = \{0\}\times \Delta$ if $n=1$, or a family of
foliated shells $\tilde p(\d\Delta^2_{\lambda})$ if $n\ge 2$.

\medskip\noindent (\slii $\Rightarrow$ (\sli Suppose now that
$\call_D^0$  contains a foliated shell $h:(B^{\eps},\call^{\v})\to
(X^0,\call^0)$. By the Theorem \ref{reparalldim2} $h$ can be
extended to a foliated meromorphic map (in particular this extension
stays to be a generic injection) $h:(P^{\eps}\setminus S,\call^{\v})
\to (X,\call)$, where $P^{\eps}$ is the $\eps$-neighborhood of the
polydisc $P$ and $S=S_1\times \bar\Delta_{1-\eps}$ for some
non-empty (!) at most countable compact $S_1\subset
\Delta_{1-\eps}$. Note that we don't need to make any
reparametrizations here. Without loss of generality we suppose that
$S_1\ni \{0\}$.  Denote by $z^0$ the image under $h$ of the point
$q=(0,1)\in P^{\eps}$ - the future reference point for the leaf
$\call_{z^0}$ which contains $h(\{0\}\times\d\Delta)$. Since
$h|_{\{z_1\}\times \d\Delta}$ is not an imbedding only for finite
set of $z_1$-s we can shrink $\Delta$ and suppose that for all
$z_1\not= 0$ the restriction $h|_{\{z_1\}\times \d\Delta}$ is an
imbedding. In fact it will be an imbedding on some annulus
$A_{1-\eps,1+\eps}$ for some $\eps >0$ - the same for all $z_1\not=
0$ -  and therefore  it will be also an imbedding on the disc
$\{z_1\}\times\Delta_{1+\eps}$ provided $z_1\not\in S_1$.

\smallskip Now we can lift $h|_{\Delta_{\eps}\times
A_{1-\eps,1+\eps}}$ for some $\eps >0$, small enough, to an
imbedding $f\deff p^{-1}\circ h: \Delta_{\eps}\times
A_{1-\eps,1+\eps}\to \hat\call_D^0$. This should be explained in
more details. Consider $p^{-1}\circ h$ along $\d\Delta_0$. It cannot
be multivalued because for $z_1\sim 0, z_1\not\in S_1$ this map is
defined and singlevalued on the disc $\Delta_{z_1}$. Moreover
$(p^{-1}\circ h)|_{\d\Delta_0}$ is also univalent. This follows from
the same property of $(p^{-1}\circ h)_{\Delta_{z_1}}$, Rouche's
theorem and absence of the holonomy in $\hat\call^0_D$. The rest is
clear.

\smallskip If $n\ge 2$ we can extend this lifting to a holomorphic
foliated imbedding $f: \Delta^{n}_{\eps}\times A_{1-\eps,1+\eps}\to
\hat\call_D^0$ (taking a smaller $\eps >0$ if necessary). This
follows from the fact that $\Delta_{\eps}\times A_{1-\eps,1+\eps}$
is Stein, so $f(\Delta_{\eps}\times A_{1-\eps,1+\eps})$ has a Stein
neighborhood (after shrinking $\eps$, see \cite{Si2}) and from the
absence of holonomy on $\hat\call_D$. At this moment we fix an {\sl
imbedded} transversal $h(\Delta^n_{\eps}\times \{q\})$ and name it
as $D$. From that moment our mappings $f$ (respectively $h$) and
their future reparametrizations are mappings over $\Delta^n_{\eps}$
(in the regions where this makes sense), \ie fibers $\Delta_z$ are
mapped into fibers $\hat\call_z$ (or $\call_z$ respectively). That
means that by $z$ we denote further both a point in
$\Delta^n_{\eps}$ and its image $h(z,q)$ in $D\subset X^0$.

\smallskip We know already that for $z =(0,...,0,z_1)\in
\Delta^{n-1}_{\eps}\times \Delta_{\eps}, z_1\not\in S_1$
$h|_{\{z\}\times\d\Delta}$ extends to an imbedding of
$\{z\}\times\Delta_{1+\eps}$ to $\call_D^0$. Therefore it extends
after a reparametrization onto $\{z\}\times\Delta_{1+\eps}$ for
$z$-s in an open non-empty subset of $\Delta^{n}_{\eps}$ (a
neighborhood of any such $(0,...,0,z_1)$). The same is true
therefore for $f=p^{-1}\circ h$. Theorem \ref{reparalldim2} gives us
an extension $\tilde h$ of $h$ after a reparametrization onto
$\widetilde{W}\setminus\tilde S$ and this extension is a foliated
meromorphic immersion which is generically injective. Therefore the
same is true for $f$, \ie $f$ extends after a reparametrization to
$\tilde f:\widetilde{W}\setminus \tilde S\to \hat\call_D$. Remark
that $\tilde S$ is not empty and up to introducing new coordinates
(locally near the fiber $\widetilde{W}_{0}$) we can suppose that
$\widetilde{W} = \Delta^{n}\times \Delta_{1+\eps}$, $\tilde S=\tilde
S_1\times \Delta_{1-\eps}$ where $\tilde S_1$ as in Theorem
\ref{reparalldim2} and $\tilde S_1\ni 0$. The diagram on the Figure 5
could be useful here.

\begin{figure}[h]
\centering
\includegraphics[width=1.5in]{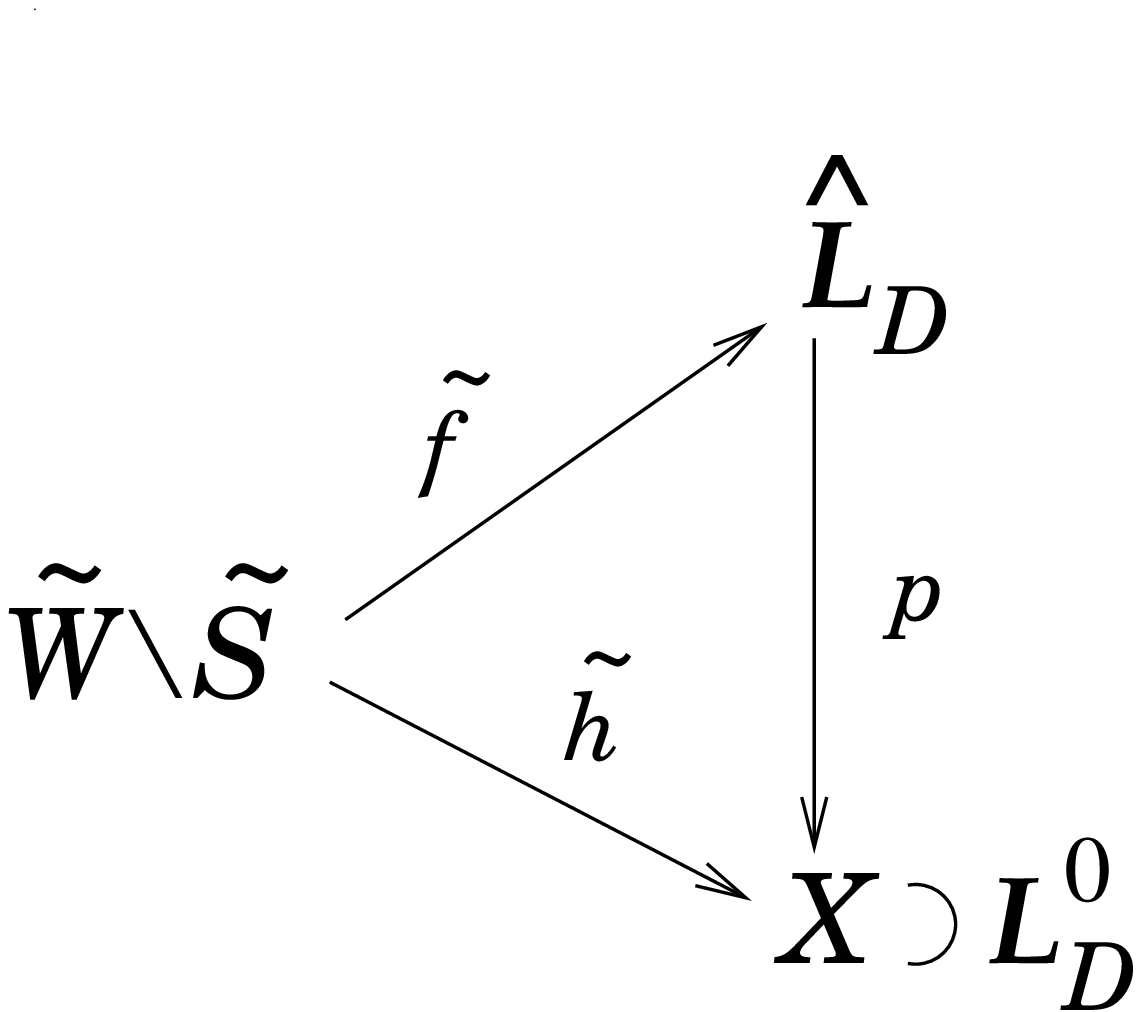}
\caption{Diagram relating $\tilde f$, $p$ and $\tilde h$: $\tilde f$
is a foliated imbedding (\ie it is {\sl holomorphic}), while both
$\tilde h$ and $p$ are, in general, meromorphic immersions.}
\label{diagram1}
\end{figure}

\smallskip We claim that $\hat\gamma_0\deff\tilde f|_{\Delta_0}(\d\Delta)$
is a vanishing cycle in $\hat\call_{z^0}$. Since for all $z_1\not\in
S_1$ the restriction $\tilde f|_{\{(0,z_1)\}\times
\d\Delta_{1+\eps}}$ is an imbedding, we get that $\tilde
f|_{\{(0,z_1)\}\times \Delta_{1+\eps}}$ is an imbedding to and
therefore $\hat\gamma_{z_1}\deff\tilde f|_{\{(0,z_1)\} \times
\d\Delta}$ is homotopic to zero in the corresponding leaf. All is
left to prove is that $\hat\gamma_0$ doesn't bound a disc in
$\hat\call_{z^0}$. But would $\hat\gamma_0$ bound a disc $\Delta^0$
in $\hat\call_{z^0}$ our foliation on $\hat\call_D$ in a
neighborhood of $\Delta^0$ would be biholomorphic to the product
$\Delta^{n}\times \Delta$ with $\Delta^0\deff\{0\}\times \Delta$ and
$\Delta^{z}\deff\{z\}\times\Delta$ being the leaves of
$\hat\call_D$. For all $(\lambda,z_1)\in \Delta^n_{\eps}$ (with
$\eps >0$ again to be taken small enough) $\tilde f$ sends
$\d\Delta_{\lambda,z_1}$ to some imbedded loop in some $\Delta_{\phi
(\lambda,z_1)}$, where $\phi :\Delta^n_{\eps}\to\Delta^n$ is some
holomorphic map sending $0$ to $0$.

\smallskip Now observe that $\area\left(p(\Delta_{\phi (\lambda,z_1})
\right)$ stays bounded when $\phi(\lambda,z_1)\to 0$.  This follows
from \cite{Ba}, see Corollary 2.4.2 in \cite{Iv5} for more details.
All is left is to remark that $p\left(\Delta_{\phi
(\lambda,z_1)}\right) = \tilde h(\Delta_{(\lambda,z_1)})$ for
$(\lambda,z_1)\not\in \tilde S_1$. Therefore this implies that
$\area (\tilde h(\Delta_{(\lambda,z_1)}))$ stays bounded as
$(\lambda,z_1)\to 0$ and $(\lambda,z_1)\not\in \tilde S_1$. But this
contradicts to (\ref{areainfty}) and to the fact that $0\in \tilde
S_1$ is an essential singular point of $\tilde h$.

\medskip\qed

In the process of proof of Theorem \ref{immshell} we saw that
vanishing cycles appear exactly in the fibers $\hat\call_z$ for $z$
belonging to the closed $(n-1)$ - polar set $S_1$ of Hausdorff
dimension $2n-2$. Therefore we obtain the following:
\begin{corol}
Let $\call$ be a holomorphic foliation by curves on a disc-convex
$(n+1)$-dimensional complex manifold $X$ which admits a pluriclosed
taming form and let $D$ be a transversal. Then the subset
$S_1\subset D$ of points $s$ such that the completed holonomy leaf
$\hat\call_s$ contains an essential vanishing cycle is complete
$(n-1)$-polar of Hausdorff dimension $2n-2$.
\end{corol}

\newprg[prgVCFS.proof2]{Imbedded vanishing cycles and proof of Theorem 2}

First of all let us make precise what we mean by an imbedded
essential vanishing cycle in the case of a {\slsf singular}
foliation. Let $\gamma_0\subset\call_{z^0}^0$ be an imbedded loop
and let $d\ge 1$ be the order of the holonomy of $\call$ along
$\gamma_0$. Denote by $\hat\gamma_0\subset \hat\call_{z^0}^0$ the
lift of $\gamma_0$. Then $p|_{\hat\gamma_0}:\hat\gamma_0\to\gamma_0$
is a regular covering of order $d$.

\begin{defi}
An imbedded essential vanishing cycle in $\call_{z^0}$ is a loop
$\gamma_0\subset\call_{z^0}$ for which the following items are
satisfied:
\begin{itemize}

\item  $\gamma_0$ is imbedded in $\call_{z^0}^0$, it admits a lift
$\hat\gamma_0$ which is imbedded in $\hat\call_{z^0}^0$ and
regularly covers $\gamma_0$ with degree $d$.

\smallskip
\item $\hat\gamma_0$ doesn't bound a disc on $\hat\call_{z^0}$.

\smallskip
\item For some (and therefore for a generic) sequence $\{z_n\}\subset D$
converging to $z^0$ there are imbedded loops $\hat\gamma_n$ in
$\hat\call_{z_n}$ uniformly converging to $\hat\gamma_0$, each
bounding a disc $D_{z_n}$ in $\hat\call_{z_n}$.

\end{itemize}
\end{defi}

\begin{rema}\rm
\label{innocent} The condition on $\gamma_0$ to be in
$\call_{z^0}^0$ and not just in $\call_{z^0}$ is not innocent at
all. One may not be able to perturb an imbedded
$\gamma_0\subset\call_{z^0}$ (which admits a lift) in the way that
this perturbation still admits a lift to $\hat\call_{z^0}$. And this
will be needed in the proof (and it is actually an important issue).
\end{rema}

Now we state the precise version of the Theorem 2 from the
Introduction.

\begin{thm}
\label{imbshell} Let $(X,\call)$ be a disc-convex foliated manifold
which admits a $dd^c$-closed taming form and let $D\subset X^0$ be a
transversal. Then the following statements are equivalent:

\smallskip
\sli Some leaf $\call_{z^0}\subset\call_D$ contains an imbedded essential
vanishing cycle.

\slii $\call_D$ contains an imbedded foliated cyclic shell.
\end{thm}
\proof (\sli $\Rightarrow$ (\slii  For a given
transversal $D\ni z^0$ we need to produce from an {\it
imbedded essential vanishing cycle} in $\call_{z^0}\subset\call_D$ an {\sl
imbedded foliated cyclic shell} in $\call_D^0$.

Take open cells $U\ni z_n$, $V\ni z^0$, $U\subset V\subset D$ such
that for an appropriate Hartogs figure $(W,\pi , U,V)\subset
\hat\call_V$ mapping $p:\hat\call_V\to X$ restricted to $W$ is a
foliated holomorphic immersion, which extends (after a
reparametrization ) to a foliated meromorphic immersion
$p:W\setminus S\to X$ as in Theorem \ref{reparalldim2} (we drop
tildes for the simplicity of our notations).

Note that $d$ is the maximal cardinality of the holonomy along loops
$\gamma_z\deff p(\d\Delta_z)\subset \call_z$ close to
$\gamma_{z^0}=\gamma_0$ for $z$ in a neighborhood of $z^0$. Find a
coordinate system $\Delta^{n-1}\times \Delta^2$ in a neighborhood of
$W_{z^0}$ in $W$ as in Theorem \ref{reparalldim2}, actually we shrink
$W$ to have $W=\Delta^{n-1}\times\Delta^2$ in the sequel. We keep
noting coordinates in $\Delta^{n-1}\times \Delta^2$ as $(\lambda ,
z_1, z_2)$.  Note that $(\lambda , z_1)$ are coordinates in a
neighborhood of $z^0$ on $D$.  Coordinates are chosen in such a way
that $z^0$ correspond to $(\lambda =0, z_1=0)$.

\smallskip Due to Theorem \ref{reparalldim2} the restriction to $S$ of
the natural projection $\pi_2:\Delta^{n-1}\times
\Delta^2\to\Delta^{n-1}$  is proper and surjective. Of course, for
that to be true one should remark here that $S\not=\emptyset$ and,
moreover, $\pi (S)=S_1\ni z^0$, because otherwise $\gamma_{z^0}$
would not be an essential vanishing cycle. Here, as usual, $S_1=\pi
(S)$ is the image of the singularity set $S$ under the natural
projection $\pi : \Delta^{n+1}\to \Delta^{n}$. By our assumption the
restriction $p|_{W_{z^0}} : W_{z^0}\setminus S_{z^0}\to \call_{z^0}
\subset X^0$ is a regular covering of order $d\geq 1$ between an
appropriate annuli in the source and target complex curves. For
$1\leq l\leq d$ denote by $A_l$ the analytic set in $\Delta^{n}$
which consists from points $q$ such that the cardinality of the
holonomy along $\gamma_q$ is at least $l$. Remark that $z^0=0\in
A_d$, $A_1=D$ and we set by definition $A_{d+1}=\emptyset$. Take a
minimal $l$ such that $S_1\cap \left(A_l\setminus
A_{l+1}\right)\not=\emptyset$. Call it $l_0$.

\medskip\noindent{\sl Case 1. $l_0=1$.}

\smallskip Take a point $s_1\in S_1\cap(A_1\setminus A_2)$ and
shrink our transversal $D$ once more to a polydisc $D=\Delta^{n}$ -
a neighborhood of $s_1$. We can suppose that $s_1=0$ in these
coordinates. If this neighborhood was taken small enough our
foliation has no holonomy along $\gamma_z$ for $z\in D$. Therefore
$p:\d_0 W|_{D}\to X^0$ is an imbedding. Consider the disc $\Delta_0
:= \{0\}\times\Delta\subset \Delta^{n-1}\times\Delta $ and consider
the restriction $W|_{\Delta_0}=\Delta_0\times\Delta$ and the
restriction of $p$ to $W|_{\Delta_0}$. Recall that $\Delta_0\cap
S_1$ is at most countable compact subset of $\Delta_0$.

\begin{lem}
\label{sard} There exists a finite union of imbedded loops $\beta
\subset \Delta_0$ which bound a relatively compact domain $G\subset
\Delta_0$ such that:

\smallskip a) $G\cap S_1\not=\emptyset$ and $\d G\cap S_1=\emptyset$.

\smallskip b) $p|_{\cup_{z\in\beta}W_z}$ is injective.

\smallskip c) Moreover, $p\left(\bigcup_{z\in\beta}W_z\cup \d_0
W|_{G}\right)$ is an imbedding.
\end{lem}
\proof As in Section 2 consider the area function
$a(z_1)=\int_{W_{z_1}}p^*\omega$ for $z_1\in \Delta_0\setminus S_1$.
Function $a$ is positive, smooth (see Remark \ref{smooth}) and tends
to infinity when $z_1\to S_1$, see Theorem \ref{areainfty} (by $S_1$
here we understand now $S_1\cap \Delta_0$ - but we not introduce any
new notations). By Sard's lemma for a generic positive $c$ the level
set $\beta_c = \{ z_1: a(z_1)=c\}$ is a union of smooth curves in
$\Delta_0$. In the sequel $c$ will be always taken bigger then $\inf
\{a(z_1):z_1\in \d\Delta_0\}$, \ie our curves will be all closed and
situated away from $\d\Delta_0$.

\smallskip\noindent{\sl Claim 1. $\beta_c$ has finite number of
irreducible components.}
Suppose not and denote by $\beta_c^i$ the sequence of irreducible
components of $\beta_c$. Let $q$ be an accumulation point of
$\beta_c^i$. $q$ belongs to $S_1$, because $\bigcup_i \beta_c^i$ is
a smooth manifold. But this contradicts Lemma \ref{area1}. Really,
$\bigcup_i \beta_c^i$ is thick at $q$ and therefore $p$ should
extend to a neighborhood of $W_q$. This contradicts to the fact that
$q\in S_1$.

Remark that we are working here with $p|_{W_{\Delta_0}}$ and use the
fact that $W_q$ contains a singular point of this restriction. This
follows from the homological characterization $(b_3)$ of essential
singularities of $p$ in Theorem \ref{reparalldim2}.

\smallskip\noindent{\sl Claim 2. $p$ is injective on $W|_{\beta_c}$.}
First of all $p$ is injective on each $W_{z_1}$,
$z\in\Delta\setminus S_1$ because it is injective on $\d W_{z_1}$.
Suppose that for some $z_1,z_2\in \beta_c$, $z_1\not=z_2$ one has
$p(W_{z_1})\cap p(W_{z_2}) \not= \emptyset$. Since $p(\d
W_{z_1})\cap p(\d W_{z_2}) = \emptyset$ we have that $p(W_{z_1})
\subset p(W_{z_2})$ (or vice verse). But this contradicts to the
fact that $\text{area} \left(p(W_{z_1})\right) = \text{area}
\left(p(W_{z_2})\right) = c$.

\smallskip
For every $i$ denote by $D^i$ the compact component of
$\Delta_0\setminus\beta_c^i$. Fix some point $s_1\in S_1$. Take one
of $D^i$-s, namely such that
$D^i\ni s_1$. Denote it as $D^1$ and its boundary
curve as $\beta^1$. If $p$ is not injective on $\d_0W|_{D^1}\cup
W|_{\beta^1}$ then there exists $z_1\in D^1$ such that
$p(W_{z_1})\cap p(W_{z_2})\not= \emptyset$ for some $z_2\in\beta^1$.
Since $p(\d W_{z_1})\cap p(\d W_{z_2}) = \emptyset$ we have two
possibilities. First: $p(W_{z_1})\supset p(W_{z_2})$ but this simply
doesn't imply that $p$ is not injective on $\d_0W|_{D^1}\cup
W|_{\beta^1}$. Therefore we are left with
the second one: $p(W_{z_1})\subset p(W_{z_2})$.

\smallskip\noindent{\sl Claim 3. If $p(W_{z_1})\subset p(W_{z_2})$ then
there exists $\beta_c^j\comp D^1$.} This is obvious, take a path from
$z_1$ to $S_1$ inside $D^1$. Then it will contain a point $z$ with
$a(z)=c$.

\smallskip If this $\beta_c^j$ surrounds our point $s_1$ call it
$\beta^2$ and
the compact component of $\Delta_0\setminus \beta^2$ call $D^2$. If
this is not the case call $\beta^1\cup\beta_c^j$ as $\beta^2$ and
the region bounded by them as $D^2$. Note that in both cases $D^2$
contains $s_1$.

\smallskip The process $D^1\supset D^2\supset ...$ is finite because
the number of $\beta_c^i$-s is finite. Therefore after a finite number
of steps we will get $D^N=:G$ and $\beta^N=:\beta=\d G$ such that $p$
injective on $\d_0W|_{G}\cup W|_{\beta}$ and $G$ has the required
properties.

\smallskip\qed

Since (taking initially $\Delta_0$ small enough) we can suppose that
$W|_{\Delta}$ is biholomorphic to $\Delta\times\Delta$ we get a
pseudoconvex domain $G\times\Delta\subset W|_{\Delta}$ such that $p$
has an essential singularity inside of this domain. By  Theorem
\ref{reparalldim2} this means that $p\left(\d (G\times\Delta)\right)$
is not homologous to zero in $X$. Set $B=\d (G\times\Delta)$, then
$p(G)$ is an imbedded foliated shell in $(X,\call)$.

\begin{rema}\rm
(a) Note that cyclic quotients didn't appear at this case, but the
topology of the shell became complicated.

\smallskip\noindent
(b) Note also that $G$ is found such that it contains an ad hoc
taken point $s_1\in S_1\cap (A_1\setminus A_2)$, \ie the constructed
shell is centered at this $s_1$. This will be used in the sequel,
see Remark \ref{less-precise}.
\end{rema}

\medskip\noindent{\sl Case 2. $l_0>1$.}

\smallskip Set $A=\bigcup_{l\geq 2}A_l$. This is a proper analytic
subset of $D$. Changing the slope of
$z_1$-coordinate and shrinking a neighborhood of $z^0$, if
necessary, we can suppose that the projection $\pi_1|_A:A\to
\Delta^{n-1}$ is proper. Here $\pi_1:\Delta^{n}\to\Delta^{n-1}$ is
the natural projection $(\lambda,z_1)\to \lambda$. $A_{l_0}\setminus
A_{l_0+1}$ contains a point $s_1\in S_1$. Shrinking $D$, if
necessary, we can suppose that $A_{l_0+1}\cap D=\emptyset$ and $D\ni
s_1$. From now on  intersect $A_{l_0}=A_{l_0}\cap D$.

\smallskip\noindent{\sl Claim 4. There exists an irreducible component
$A'$ of $A_{l_0}$ of pure dimension $n-1$ which is entirely contained
in $S_1$.}

\smallskip Choose coordinates $(\lambda ,z_1)$ in a neighborhood of
$s_1=(0,0)$ in $D=\Delta^{n-1}
\times\Delta$ in such a way that $\pi_1|_{A_{l_0}}$ is proper. If
$\dim A_{l_0}<n-1$ then $\dim \pi_1(A_{l_0})<n-1$. But we know that
for every $\lambda\in \Delta^{n-1} \setminus \pi_1(A_{l_0})$ there
exists at least one $z_1$ such that $(\lambda ,z_1)\in S_1$. Remark
also that the holonomy along $\gamma_{\lambda,z_1}$ for such
$z=(\lambda,z_1)$ is less then $l_0$. Contradiction to the
definition of $l_0$.

Therefore $\dim A_{l_0}=n-1$. Note that $S_1\subset A_{l^0}$ by the
definition of $l_0$. If there exists a point $q\in A_{l^0}\setminus
S_1$ then from homological characterization $(b_3)$ in Theorem
\ref{reparalldim2} it follows that no point of $A_{l^0}$ in a
neighborhood of $q$ belongs to $S_1$. Therefore all irreducible
components of $A_{l_0}$ intersecting this neighborhood do not belong
to $S_1$. In this way we find an irreducible component $A'$ of
$A_{l^0}$ which is entirely contained in $S_1$.

\smallskip From now on we can suppose that $A'=S_1$ is smooth and is
given by the equation
$z_1=0$ in $D$. Let $g:D\to D$ be a local biholomorphism generating
the holonomy along $\gamma_0$. Remark that $g|_{A'}\equiv \id$ and
$g^{l_0}\equiv \id$. $\gamma_0$ here is the boundary $\d\Delta_0$
and $s_1=0$.

\smallskip\noindent{\sl Claim 5. In an appropriate coordinates with
center at $s_1$ the
automorphism $g$ has the form $g (\lambda ,z_1) = (\lambda
,e^{\frac{2\pi i l}{l_0}}z_1)$ for some $l\in \{1,...,l_0\}$
relatively prime with $l_0$.}

\smallskip This is a nearly standard fact which easily follows
from the famous Bochner's linearization theorem, see \cite{Bo}. Really,
literally repeating the proof of Theorem 1 from \cite{Bo} one can find
coordinates in which $g $ is linear and still preserving $A'=\{z_1=0\}$.
Therefore in an appropriate coordinates $g$ has the form $g
(\lambda,z_1)=(\lambda,e^{\frac{2\pi i l}{l_0}}z_1)$ for some $l\in
\{1,...,l_0\}$ relatively prime with $l_0$.

\medskip
Factorize
$\Delta^{n-1}\times\Delta\times\Delta$ by the
action $(\lambda, z_1,z_2)\to (\lambda ,e^{\frac{2\pi
il}{l_0}}z_1,e^{\frac{2\pi i}{l_0}}z_2)$ to get $\Delta^{n-1}\times
\calx^{l,l_0}$, where $\calx^{l,l_0}$ is a surface with cyclic
quotient singularity. We get a holomorphic foliated immersion
$p:\Delta^{n-1}\times \calx^{l,l_0}\setminus S^{l,l_0} \to
X$, where  $S^{l,l_0}$ - image of $S$ under the factorization.

\smallskip Remark that $p|_{\d (\{0\}\times \calx^{l,l_0})}$ is now an
imbedding.
Therefore we can repeat arguments of Lemma \ref{sard} and prove that
$p|_{\d (\{0\}\times \calx^{l,l_0})}$ is injective in an
neighborhood of the boundary $B$ of the domain
$W_{l,l_0}=\bigcup_{z\in G} W_{z}$ for some $G\comp \Delta$. Would
$p(B)$ be homologous to zero in $X$ then by $(b_3)$ of the Theorem
\ref{reparalldim2} would imply the extensibility of $p$ onto our
domain $W_{l,l_0}$ and this is not the case.  I.e. we got an {\sl
imbedded} foliated cyclic shell.

\begin{rema}\rm
We silently used here a version of Theorem \ref{reparalldim2} in the
spaces with cyclic singularities. One can either prove such version
directly, or ``lift'' the problem to the covering of $\calx$ (which
is a bicylinder), apply extension there and push the extended map
down. This is possible, because the extended map will be also
invariant under the action of the cyclic group by the uniqueness
theorem for holomorphic functions.
\end{rema}

\medskip\noindent (\slii $\Rightarrow $ (\sli  Let $h:(B^{\eps},
\call^{\v})\to (X^0,\call^0)$
be an imbedded foliated cyclic shell. We can proceed literally as
in the proof of (\slii $\Rightarrow $ (\sli of Theorem 3.1. All we
need to do is to see that the cycle
$h|_{\{0\}\times\d\Delta}:\{0\}\times \d\Delta\to \call_{z^0}^0$ -
proved to be a vanishing one - was imbedded from the every beginning.
Further details will be omitted. Theorem is proved.

\smallskip\qed

\begin{rema}\rm
\label{less-precise} Note, that in the Theorem \ref{imbshell} the
place for the shell is less precise then in Theorem \ref{immshell}.
But let us still make a precision here. Let
$D=\Delta^{n-1}\times\Delta$ in a neighborhood of $z^0=(0,0)$ as
above. Then we proved, in fact, that we can find $\lambda\in
\Delta^{n-1}$ arbitrarily close to $0$ such that
$\call^0_{\Delta^2_{\lambda}}$ will contain an imbedded foliated
cyclic shell centered at given $s_1\in \Delta_{\lambda}$. If $n=1$
then this $\lambda$ is $0$. Remember that we were able to center our
shell in a generic point $s_1$ on $S_1$ near $z^0$, see Remark 3.6.
\end{rema}

\begin{rema} \rm
With the definition of a foliated shell of this paper, analogues of
implications (\sli $\Rightarrow$ (\slii of Theorems 3.1 and 3.2 can
be found for surfaces in \cite{Br5}.
\end{rema}

\newsect[sect.PCFS]{Pluriclosed Metric Forms and Foliated Spherical
Shells.}

\newprg[prgPCFS.sph-she]{Pluriclosed metric forms and foliated
spherical shells}

Up to now our immersed shells were boundaries of the bicylinder (or
pseudoconvex hypersurfaces close to it). One might ask if the
$CR$-geometry is relevant here? The test question would be: can one
take as shells the images of the standard spheres (with the standard
vertical foliation) and not such a Levi-flat objects as boundaries
of bicylinders? In the context of this paper this issue  goes
together in one line with reducing of the size of the essential
singularity set $S$ that is ``virtually present'' in the heart of
all our proofs. And this task is of capital importance. It appears
to be  crucial for getting from vanishing cycles the {\sl imbedded}
ones.

\smallskip At present we are able to reduce the size of $S$
(equivalently to pass to spheres as shells) only in the case when
our $dd^c$-closed taming form is actually a metric form on $X$, \ie
$\omega$ should be not just a $dd^c$-closed form positive in the
directions tangent to $\call$ but in all directions in $TX$.

\smallskip Let $B=\sss^3=\{ z=(z_1,z_2)\in \cc^2:||z||=1\}$ denote
the unit sphere in $\cc^2$, $P = \{ z\in \cc^2:||z||<1\}$ - the unit
ball. For some $0<\eps<1$ set $B^{\eps}=\{ z\in \cc^2: 1-\eps <
||z|| < 1+\eps \}$ - a shell around $\sss^3$. As before, denote by
$\pi :\cc^2\to\cc$ the canonical projection $\pi (z)=z_1$ onto the
first coordinate of $\cc^2$. Note that $B^{\eps}$ is foliated by
$\pi$ over the disc $\Delta_{1+\eps}$ of radius $1+\eps$. Denote
this foliation again as $\call^{\v}$. Its leaves
$\call_{z_1}:=\pi^{-1}(z_1)$ are discs if $1-\eps <|z_1|<1+\eps  $
and are annuli if $|z_1|<1-\eps$.

\begin{defi}\rm
The pair $(B^{\eps}, \call^{\v})$ we shall call the {\sl standard
foliated spherical shell}.
\end{defi}

Let $h:(B^{\eps},\call^{\v})\to (X^0,\call^0 )$  be some generically
injective foliated holomorphic immersion of the standard foliated
spherical shell into $(X^0,\call^0)$. Denote by $\Sigma$ the image
of the unit sphere $\sss^3$ under $h$.

\begin{defi}\rm
$h(B^{\eps})$ is called {\sl a foliated spherical shell} in
$(X,\call )$ if $\Sigma$ is not homologous to zero in $X$.
\end{defi}

\begin{figure}[h]
\centering
\includegraphics[width=2.0in]{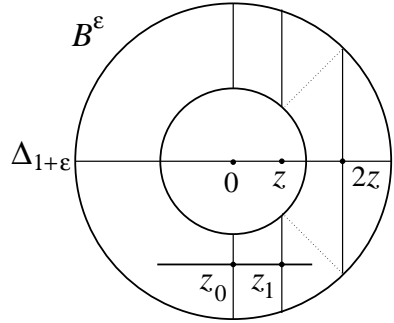}
\caption{The "vertical foliation" on the Hopf surface is again the
simplest example. The leaf $\call_{0}$ is a torus,
$\call_{z}=\call_{2z}$ is $\cc$ for $z\not= 0$. The cycle $\gamma =
\{(0,z): |z|=1\}$ is a vanishing cycle. Image of the $\eps$-
neighborhood of $\sss^3$ together with the "vertical foliation"
under the natural projection $\cc^2\setminus \{0\}\to H$ is a
foliated spherical shell in $(H^2,\call)$.} \label{shell_sp-fig}
\end{figure}

\begin{rema}\rm
\label{small} {\bf (a)} Let us recall the Main Theorem from
\cite{Iv6}, where we worked with pluriclosed {\sl metric} forms. We
proved there that the singularity set $S$ from the Theorem
\ref{reparalldim2} is "small" in the sense that for every
$\lambda\in \Delta^{n-1}$ (see notations in Theorem
\ref{reparalldim2}) the set $S_{\lambda}\deff S\cap
\Delta^2_{\lambda} $ is a complete pluripolar compact of
$\Delta^2_{\lambda}$ of Hausdorff dimension zero.

\smallskip\noindent{\bf (b)} The arguments of Lemma \ref{descrete}
can be repeated here and give that $S_{\lambda}$ are, in fact, at
most countable. The crucial issue, however, here is not a
countability (null-polarity is perfectly sufficient) but the size of
the sets $S_{\lambda, z_1}$ for $(\lambda , z_1)\in S_1$. If the
taming form is actually a metric form then these last sets are also
small, \ie at most countable. It is this fact which leads to the
foliated {\slsf spherical} shells in the Proposition below and
finally allows to produce imbedded vanishing cycles and shells.
\end{rema}

\begin{prop}
Let $(X,\call,\omega)$ be a disc-convex pluritamed holomorphic
foliation by curves. Suppose that $\omega$ is actually a metric form
and that a foliated manifold $(X,\call)$ contains a foliated shell
$h:(B^{\eps},\call^{\v})\to (X^0,\call^0)$. Then:

\sli $h$ extends to a foliated meromorphic immersion of $P^{\eps}
\setminus S$ where $S$ is at most countable compact subset of $P$.

\slii  $(X,\call)$ contains a foliated spherical shell.
\end{prop}
\proof The proof is  immediate because as a shell we can take a
standard $3$-sphere $\sss^3_r(s)$ around any point $s$ of $S$. A
radius $r$ should be chosen in such a way that $\sss^3_r(s)\cap
S=\emptyset$. And this is possible due to the null-dimensionality of
$S$.

\smallskip Countability of $S$ (absent in \cite{Iv6}) can be now achieved
due to the Part $(b_3)$ of Theorem \ref{reparalldim2}.

\smallskip\qed

\begin{rema}\rm
(a) In this proof we didn't use the condition on $h$ to be a
foliated map, because any holomorphic mapping $h:B^{\eps}\to X$,
where $X$ admits a pluriclosed metric form, extends to a meromorphic
map from $P^{\eps}\setminus S$ to $X$ with $S$ being at most
countable compact subset of $P$.

\smallskip\noindent (b) The conclusion of this Proposition remains valid
(with countability replaced by null-polarity) if $\omega$ is
supposed to be a pluriclosed taming form for $\call$ and there
exists some other plurinegative metric form $\omega_1$ on $X$
(irrelevant to $\call$). Really, all we need is to reduce the size
of the essential singularity set $S$ along ``$z_2$-direction'' and
this can be done with the help of $\omega_1$.
\end{rema}

\smallskip

\newprg[prgPCFS.alm-hart]{Almost Hartogs property of foliated pairs}
It occurs that the eduction of the size of $S$ already made is
exactly what one needs in order to produce imbedded vanishing
cycles. Let us formalize this by giving the following (we use
notations from Theorems \ref{reparalldim1} and \ref{reparalldim2}):

\begin{defi}
\label{alm-hart-def} A foliated manifold $(X,\call)$ of dimension
$n+1\ge 2$ is called {\sl almost Hartogs} if the following is
satisfied:

\sli  Every foliated holomorphic immersion $h:(W,\pi , U , D)\to
(X,\call)$ of a non-trivial generalized Hartogs figure of dimension
$\dim X$ extends to a  foliated meromorphic immersion of
$(W\setminus S,\pi , D)$ into $(X,\call)$ after a reparametrization,
where $S$ is a closed subset of $W$ of zero Hausdorff
$(2n-2)$-dimensional measure.

\slii Moreover, the essential singularity set $S$ (\ie the minimal
set with property as in (\sli) has the following structure:

\smallskip\quad a) for every point $s^0\in S$ there exists a neighborhood of
it biholomorphic to $\Delta^{n+1}$ such that the restriction of
$\pi$ to this neighborhood is the natural vertical projection
$\pi:\Delta^{n+1}\to \Delta^n$;

\smallskip\quad b) the restriction to $S\cap \Delta^{n+1}$ of the natural
projection $\pi_1:\Delta^{n+1}\to \Delta^{n-1}$ is proper and has at
most countable fibers.
\end{defi}

As usual ``meromorphic foliated immersion'' means here that the
extended $h$ is a foliated immersion outside of its indeterminacy
set. However, one should remark that the only point here is to
extend $h$: if a meromorphic extension of $h$ onto
$\Delta_{1+\eps}^{n+1}\setminus S$ is possible then it will be
automatically a foliated immersion outside of its indeterminacy set.
If $S$ happens to be empty for every such mapping into $(X,\call)$
then the latter is called simply ``Hartogs''. Needless to say that
the set $S$ appearing here is always closed. In the case of the
presence of a plurinegative taming form on $(X,\call)$ the item
(\sli is automatic by Theorem \ref{reparalldim1} and only the item
(\slii represents itself a condition.

\smallskip Our goal in this subsection is to reduce
the problem of finding imbedded vanishing cycles in a shelled
foliations to the proof of the almost Hartogs extension property  of
$(X,\call)$. And the latter can be proved in many interesting cases,
see \cite{Iv1,Iv2,Iv3,Iv4}. In particular, the Theorem 3.3 from
\cite{Iv6} (Proposition 4.1 of the present paper) can be restated in
the following form:

\begin{prop}
\label{metr-alm-hart}
Suppose that a foliated manifold $(X,\call)$ admits a pluriclosed
taming form $\omega$, such that $\omega $ is actually a metric form.
Then $(X,\call)$ is almost Hartogs.
\end{prop}

One more example is a result from \cite{Iv4} (it doesn't require
any special metric form on the total space $X$):

\begin{prop}
\label{turb-alm-hart}
Suppose that the manifold $X$ is an elliptic fibration (with
possibly singular fibers) over a disc-convex K\"ahler manifold $Y$.
Then every holomorphic foliation by curves on $X$ is almost Hartogs.
\end{prop}

Really, let $f:(W,\pi,U,V)\to X$ be a holomorphic map. If $p: X\to
Y$ is the holomorphic mapping defining the elliptic fibration then
the composition $p\circ f$ extends onto $W$ after a
reparametrization by \cite{Iv3} and \cite{Br3}. Following the
arguments in \cite{Iv4} one gets an extension of $f$ onto
$W\setminus S$ where $S$ is the indeterminacy set of $p\circ f$
(reparametrizations do not cause any problems here). One also
obviously has he following:

\begin{prop}
\label{ric-alm-hart}
Suppose that the manifold $X$ is a rational fibration (with
possibly singular fibers) over a disc-convex K\"ahler manifold $Y$.
Then every holomorphic foliation by curves on $X$ is almost Hartogs.
\end{prop}

\newprg[prgPCFS.imbed]{Imbedded vanishing cycles}

Recall that the classical result of Ohtsuka states the following: if
$h:\Delta^*\to P$ is a holomorphic map of a punctured disc to a
hyperbolic Riemann surface, then $h$ extends to zero as a
holomorphic mapping from $\Delta$ either to $P$, or to a {\slsf
completed} by one point surface $R\deff P\cup \{b_0\}$ with
$b_0\deff h(0)$. One can express this by saying that in the
situation described above $h$ cannot have an essential singularity
at zero. See \cite{Oh1,Oh2} and for a much simpler proofs see
\cite{Re, Ro}.

\smallskip We shall need an orbifold version of Ohtsuka's theorem
in this paper. The reason is the following: according to Remark
\ref{orbifold1} our holonomy covering map $h \deff
p_{z_0}:\hat\call_{z^0}\to\call_{z^0}$ is an {\slsf orbifold
covering} map (this follows from the very definition of the
vanishing end) and we like to prove that in the case when
$\hat\call_{z^0}$ is hyperbolic this map behaves near an eventual
puncture as in the Ohtsuka's theorem, \ie extends to a puncture
after a one point completion of $\call_{z^0}$. A proof using the
metric structure of orbifolds may be found in \cite{Br5}, \S 3.3,
while we present here a group theoretic alternative.

\smallskip For all notions and facts about orbifold Riemann surfaces
that we are going to use below we refer to the book of Milnor \cite{Mi1}
and references there. For the rudiments on Fuchsian groups see
\cite{Be}. Recall that for a Riemann surface orbifold $(R,\nu)$ the Euler
characteristic is defined as
\begin{equation}
\label{euler}
\chi (R,\nu) = \chi (R) + \sum_j\left( \frac{1}{\nu (z_j)} - 1\right),
\end{equation}
where $\chi (R)$ is the Euler characteristic of the underlying
Riemann surface $R$ and $\nu (z_j)$ is the value of ramification
function $\nu $ at ramification point $z_j$. Riemann surface
orbifold $(R,\nu)$ is called hyperbolic if its orbifold universal
covering $\tilde S_{\nu}$ is the unit disc and parabolic in the
opposite case. According to Lemma E.4 from \cite{Mi1} the Riemann
surface orbifold $(S, \nu)$ is hyperbolic if and only if $\chi (S,
\nu) <0$. A {\slsf regular} holomorphic map $h: S \to (R,\nu)$ from
a Riemann surface $S$ to a Riemann surface orbifold is by definition
a holomorphic map $h:S\to R$ such that for every $z\in S$ the branch
index of  $h$ at $z$ is equal to $\nu (h(z))$ (ex., no branching
whenever $\nu (h(z))=1$).

\begin{lem}
\label{ohtsuka} Let  $h:\Delta^*\to (P, \nu)$ be a regular
holomorphic map from the punctured disc to a hyperbolic Riemann
surface orbifold. Then $h$ cannot have an essential singularity at
the origin. More precisely:

\smallskip\sli either $h$ holomorphically extends to zero as a
mapping with values in $P$,

\smallskip\slii or, there exists a Riemann surface $R\supset P$ such
that $R\setminus P = \{b_0\}$ and $h$ holomorphically extends to
zero as a mapping to $R$ with $h(0) = b_0$.
\end{lem}
\proof We are following \cite{Ro}. Denote by $\Gamma$ the Fuchsian
group of deck transformations  of the universal covering $H \to
P=H/\Gamma$ ($H$ stands for the upper half plane). Wright $\Delta^*
= H/G$, where $G$ is generated by the translation $T_1(z) = z+1$.
Mapping $h$ induces a homomorphism $h_*:G \to \Gamma$. As in
\cite{Ro} we get that if for a circle $h_*(T_1) = 0$ then $h$ lifts
to $\tilde h:\Delta^*\to H$. In that case the Riemann extension
theorem applies and gives the extension of $\tilde h$ (and therefore
of $h$) to the origin.

\smallskip Suppose now that $h_*(T_1)\not= 0$ be not homotopic to zero in $P$.
Write $h_*T_1=T^n$, where $T$ is primitive. Then $h$ lifts to a
mapping $\tilde h : \Delta^*\to H/\Gamma_T$ where $\Gamma_T$ is the
cyclic subgroup of $\Gamma$ generated by $T$. From here and by Lemma
3 from \cite{Ro} we get that $T$ is parabolic and therefore
$H/\Gamma_T = \Delta^*$. If we prove that the natural mapping $\psi
: \Delta^* = H/\Gamma_T \to H/\Gamma = P$ extends to the puncture
(after, may be, completing $P$) our lemma will be proved, because
then $h=\psi \circ \tilde h$ will extend to.

\smallskip It will be convenient to break the proof into two cases.

\smallskip\noindent{\slsf Case 1. The orders of elliptic elements of
$\Gamma$ are uniformly bounded.} In that case we can apply the
result of Purzitsky, see \cite{Pr}: $\Gamma$ contains a torsion free
subgroup $\Gamma_1$ of finite index. Since $H \to H/\Gamma_1$ is an
unbranched covering the quotient $H/\Gamma_1$ is hyperbolic.
Therefore the theorem of Ohtsuka applies: for $r>0$ small enough
$\Delta^*_r$ projects properly   onto some puncture
$\Delta^*_{\rho}\subset H/\Gamma_1$. But $H/\Gamma_1$ is a {\slsf
finite} branched covering of $H/\Gamma = P$. Therefore under the
resulting covering $\psi$ our puncture $\Delta^*_r$ is mapped onto a
neighborhood of a puncture in $H/\Gamma $ and we are done.

\smallskip\noindent{\slsf Case 2. The orders of elliptic elements of
$\Gamma$ are not bounded.} Let $f\deff \{f_1,f_2,...\}\subset P$ be
the images of centers of all elliptic elements of $\Gamma$. This
sequence is infinite by assumption. Denote by  $F$ the set of all
elliptic elements of $\Gamma$ with fixed points which project to
$\{f_6,f_7,...\}$. Let $\Gamma_F$ be the  subgroup of $\Gamma$
normally generated by $T$ and $F$. Remark that $\psi : \Delta^* \to
P$ lifts to some $q:\Delta^*\to H/\Gamma_F$, \ie is a composition of
$q$ with the projection $\psi_1: H/\Gamma_F\to P$. But now
$H/\Gamma_F$ is hyperbolic by formula (\ref{euler}). Indeed,
$H/\Gamma_F\to P$ is ramified over $f_1,...,f_5$ and therefore the
orbifold Euler characteristic of $H/\Gamma_F$ is negative whatever
$\chi (P)$ is. Therefore $q$ extends to a puncture by Ohtsuka's
theorem (after, may be, completing $H/\Gamma_F$).  At the same time
remark that the group $\Gamma /\Gamma_F$ of the deck transformations
of the cover $\psi_1:H/\Gamma_F \to P$ has torsions only over the
centers $f_1,...,f_5$, \ie their orders are uniformly bounded. That
means that we are now under the Case 1 and this finishes the proof.

\smallskip\qed

\smallskip Now we are prepared to state the main result of this
Section, which is a precise version of Theorem \ref{imb-shel} from
the Introduction:

\begin{thm}
\label{virtual1} Let $(X,\call)$ be a disc-convex foliated manifold
and let $f:(H_{\eps},\call^{\v})\to \hat\call_D$ be a foliated
holomorphic imbedding of the standard  Hartogs  figure into the
holonomy covering cylinder $\hat\call_D$ for some transversal
$D\subset X^0$. Suppose that:

\smallskip 1) $h\deff p\circ f$ extends as a foliated meromorphic
immersion to a complement of a closed subset
$S\subset\Delta^{n+1}_{1+\eps}$ of the form $S=\cup_{z\in
S_1}S_{z}$, where $S_1$ is $(n-1)$-pluripolar in $\Delta^n_{1+\eps}$
and all $S_z$ are at most countable.

\smallskip 2) For some $z^0\in \Delta_{1+\eps}^n$ and the disc
$\Delta_{z^0}\deff \{z^0\}\times\Delta_{1+\eps}$ the cycle
$\hat\gamma_0\deff f|_{\Delta_{z^0}}(\d\Delta_{z^0})\subset
\hat\call_{z^0}$ is an imbedded essential vanishing cycle in the
holonomy covering leaf $\hat\call_{z^0}$.

\smallskip Then the leaf $\call_{z^0}$ itself contains an  imbedded
essential vanishing cycle $\gamma_0\subset \call_{z^0}^0$.
\end{thm}

\proof  Note that from (1) we get that $f$ itself extends as a
foliated imbedding of $\Delta^{n+1}_{1+\eps}\setminus S$ into
$\hat\call_D$. The condition that every point $s\in S$ is an
essential singularity of $h$ (and therefore also of $f$) means that
there exists no neighborhood $V\ni s$ such that $h$ (and $f$)
meromorphically (holomorphically)  extends to $V$. Note also that
$\Delta_{z^0}$ intersects $S$, otherwise
$f|_{\Delta_{z^0}}(\d\Delta)$ cannot be a vanishing cycle.

We shall work locally around point $z^0\in \Delta^n_{1+\eps}$ and
therefore we shall take coordinates in  which this point is the
origin $0$. $h_0\deff h|_{\Delta_0}:\Delta_0\setminus S_{0}\to
\call_{z^0}$ is a holomorphic mapping of a pluri-punctured disc
$\Delta_0\setminus S_0$ to the Riemann surface $\call_{z^0}$ which
factors as $h_0=p_0\circ f_0$ through the holomorphic imbedding
$f_0:\Delta_0\setminus S_0\to \hat\call_{z^0}$. Here $S_0\deff
S\cap\Delta_0$ is at most countable compact in $\Delta_0$. Take some
isolated point in $S_0$, suppose it is the origin and remark that
for a boundary of a small disc around the origin its image by $f$ is
an imbedded essential vanishing cycle in $\hat\call_{z^0}$.
Therefore we shrink our polydisc to $\Delta^{n+1}$ to be as small as
necessary to have  that $0$ is the only intersection point of
$\Delta_0$ with $S$, \ie  $\{0\}=S_0=S\cap \Delta_0$.

\smallskip Remark that $f_0:\Delta_0\setminus \{0\}\to \hat\call_{z^0}$
extends to a holomorphic imbedding
of the disc $\Delta_0$ into a Riemann surface $\hat R$ which is
obtained from $\hat\call_{z^0}$ by adding to it a point, \ie $\hat
R\setminus \hat\call_{z^0} = \{ a_0\}$ and $a_0$ is the image of $0$
under the extended map (which we still denote as $f_0:\Delta_0\to
\hat R$). This follows easily from the fact that
$f_0:\Delta_0\setminus \{0\}\to \hat\call_{z^0}$ is an imbedding.

\smallskip
If $\call_{z^0}$ is hyperbolic  the mapping $p$ restricted
to $\hat\call_{z^0}$ extends to a holomorphic map $p_0:\hat R\to R$
which is a ramified covering between neighborhoods of $a_0$ and
$b_0$.  Here $R\deff \call_{z^0}\cup \{b_0\}$ with
$b_0\deff h_0(0)=f(a_0)$. This results to an imbedded essential vanishing cycle in
$\call_{z^0}$.

\smallskip If $\call_{z^0}$ is orbifold-hyperbolic, \ie hyperbolic is
$\hat\call_{z^0}$, then Lemma \ref{ohtsuka} still does the job as above.

\medskip In all other cases $\call_{z^0}$ will be parabolic, \ie
torus, sphere, plane or punctured plane, as well as $\hat\call_{z^0}$.
These few cases can be listed explicitly with the help of \cite{Mi1,Mi2}.
Note that in all these cases we have both $\chi (\call_{z^0})\ge 0$ and
$\chi (\hat\call_{z^0})\ge 0$.
If $h_0$ extends to zero as a mapping from $\Delta_0$ to $R=\{b_0\}\cup
\call_{z^0}$  then everything goes as above. Therefore below we shall
be concerned with $h_0$ not extending to the origin in the sense just
described, \ie with $p_0$ having an essential singularity at ``added''
point $a_0$.

\smallskip

\smallskip\noindent{\sl Case 1. $\chi (\hat\call_{z^0}) > 0$ and
$\call_{z^0}$ is compact.} In that case it should be a sphere as
well as $\hat\call_{z^0}$, and the covering
$p:\hat\call_{z^0}\to\call_{z^0}$ should be finite, see Remark E.5
in \cite{Mi1}. This case is trivial, \ie a vanishing cycle doesn't
occur ($\hat\call_{z^0}$ cannot contain a curve, which is nonhomotopic
to zero).

\begin{rema} \rm
By Lemma \ref{ratgeo1}, which will be proved in the last Section in this
case $(X,\call)$ is a rational quasifibration.
\end{rema}

\smallskip\noindent{\sl Case 2. $\call_{z^0}$ is non-compact.}
I.e., is $\cc$ or $\cc^*$. Then the formula (\ref{euler}) tells us
that $\call_{z^0}$ can be either $\cc$ with one ramification point,
or $\cc$ with two of index two, or $\cc^*$ with no ramifications.
All these cases are trivial, \ie we always get an imbedded vanishing
cycle.

\smallskip Now we consider the cases when $\call_{z^0}$ is compact
and $\chi (\hat\call_{z^0}) = 0$.

\smallskip\noindent{\sl Case 3. $\call_{z^0}$ is a torus.}
In that case formula (\ref{euler}) tells us that
$p:\hat\call_{z^0}\to\call_{z^0}$ is an unramified covering. Now
every loop in $\ttt^2$ is homotopic to a multiply covered imbedded
one and therefore $h_0(\d\Delta_0)$ is homotopic to a multiply
covered imbedded loop $\gamma_0$ and this homotopy lifts again to
$\hat\call_{z^0}$. This again produces an imbedded essential
vanishing cycle.

\smallskip In the last two cases $\call_{z^0}$ is a sphere.
Then the formula (\ref{euler}) tells that $p:\hat\call_{z^0}\to\call_{z^0}$
is a ramified covering with either three or four ramification points $\{z_j\}$
with  multiplicity function $\nu$ satisfying

\begin{equation}
\label{nu}
\sum_j \left(1 - \frac{1}{\nu (z_j)}\right) = 2.
\end{equation}

There are only four integer solutions of (\ref{nu}), see \cite{Mi1}
Remark E.6 and \cite{Mi2} Corollary 4.5 for more details. Here we
only list them together with the needed facts.

\begin{itemize}

\item The (orbifold) universal covering $\tilde \call_{z^0}$ of $\call_{z^0}$
(\ie the usual universal covering of $\hat\call_{z^0}$) is $\cc$  in all these
cases and the group of deck transformations of the
covering $\tilde p_{z^0} : \tilde \call_{z^0}\to\call_{z^0}$ is the extension of
$\zz^2$ by a finite group
$\zz_n$ of $n$-roots of $1$ for $n=2,3,4,6$ (\ie one has four options). In another
words the group in
question is $\zz^2\rtimes \zz_n$ - the semidirect product of $\zz^2$ with $\zz_n$.

\item $\zz^2$ acts on $\cc$ by translations along some lattice $\Lambda$ and $\zz_n$
by rotations onto the angle $e^{\frac{2\pi i}{n}}$.

\item In the case $(2,2,2,2)$ the lattice $\Lambda$ is generated by $1$ and
$ \tau$ where $\tau$ is an arbitrary
complex number which belongs to the Siegel region $\cals \deff \{
\tau : |\tau |\ge 1, |\re (\tau) |\le 1/2, \im (\tau)>0 \text{ and }
\re (\tau)\ge 0 \text{ if } |\tau |=1 \text{ or } |\re
(\tau)|=1/2\}$ - the fundamental domain of $SL(2,\zz)$. The finite
group is $\zz_2=\{\pm 1\}$ in this case.

\item In all other cases the lattice is rigid, \ie unique, and is determined by the
condition to be invariant under the rotations from $\zz_n$ for $n=3,4,6$.

\end{itemize}

With this information at hand one should distinguish here two cases.

\smallskip\noindent{\sl Case 4. $\call_{z^0}$ is a sphere and the ramification
function is one of $(2,4,4)$, $(2,3,6)$, $(3,3,3)$.}

\smallskip Recall the following well known fact:

\smallskip\noindent{\it The group $G_n\deff\zz^2\rtimes\zz_n$ for $n=3,4,6$ has no
nontrivial normal subgroups of infinite index.}

\smallskip Really, let $N\triangleleft G_n$ be a nontrivial normal subgroup. We see
it as acting on $\cc$ as described. Suppose $N$ contains a rotation $\rho$.  Take any
translation $t\in G_n$. Then the commutator $t_1\deff[\rho ,t]$ is (obviously) a
translation and it belongs to $N$ because $[\rho ,t]=\rho
(t\rho^{-1}t^{-1})$ and $N$ is normal. But $t_2\deff\rho t\rho^{-1}$ is a translation
transversal to $t_1$
(if $n\not= 2$) and therefore $N\supset \zz\cdot t_1\times\zz\cdot t_2$ and we are done.

\smallskip Remark now that the group $N$ of deck transformations of the covering
$\tilde\call_{z^0}\to \hat\call_{z^0}$ should be a normal subgroup of the group $G_n$
of the deck transformations of the covering $\tilde\call_{z^0}\to \call_{z^0}$. By the
fact, just mentioned, $N$ is either trivial or of finite
index in $G_n$. In both cases there cannot be any vanishing cycles in $\hat\call_{z^0}$.

\smallskip\noindent{\sl Case 5. $\call_{z^0}$ is a sphere and the ramification
function is $(2,2,2,2)$.}

This case is not rigid in the sense that there is one conformal
parameter, namely the cross-ratio of four (ramification) points on
$\cc\pp^1$. But this doesn't matter. Again, if $N\triangleleft G_2$
is a nontrivial normal subgroup then it contains a translation $t_1$
as it was explained above. So $N\supset \zz\cdot t_1$. If $N$
contains also a rotation then it contains also an another
translation $t_2=[\rho ,t]$ transversal to $t_1$ if $t$ was taken
transversal to $t_1$, the proof goes exactly as above.

\smallskip Therefore we are left with the case  $N=\zz\cdot t$, \ie
the group $N\triangleleft G_2$ of the deck transformations of the
cover $\tilde\call_{z^0}\to \hat\call_{z^0}$ can be only $\zz\cdot
t$ in this case (other cases are trivial). Take $t=k$ for simplicity
(after and appropriate choice of a basis for $\Lambda$, \ie $1,
\tau$ as above). Then $\hat\call_{z^0}$ is a cylinder $\cc/\zz\cdot
t$, \ie $\hat\call_{z^0}= \cc/k\zz = [0,k]+\rr\tau$ with left and
right boundary lines identified by $z\to z+k$. Every imbedded loop
$\hat\gamma_0$ in this cylinder is homotopic to the interval
$[0,k]$. Covering $\hat\call_{z^0}\to \call_{z^0}$ is a composition
of a unramified $k$-sheeted covering $p_1:\cc/\zz\cdot k\to\cc/\zz$
and a ramified one $p_2:\cc/\zz\to \cc/\zz^2\times\zz_2$. Under the
first mapping $\hat\gamma_0$ maps to a $k$-times taken imbedded loop
$[0,1]$. This loop is homotopic to ($k$-times taken)
$[0,1]+\frac{i}{4}$ and the last lies entirely in the fundamental
domain of $G_2$ (only the ends are identified). Therefore it
projects to an imbedded loop $\gamma_0$ in the factor $\call_{z^0}$.
As a result we got an imbedded vanishing cycle. Figure 5 from
\cite{Mi2} might be helpful for better understanding the last few
lines above. Theorem is proved.

\smallskip\qed

\smallskip Theorem 3 from the Introduction follows now immediately
from  Theorem \ref{virtual1} and Proposition 4.2. More precisely, we
obtain the following result.

\begin{corol}
\label{metric-hart}
Let $(X,\call, \omega)$ be a disc-convex, singular holomorphic
foliation by curves such that the pluritaming form $\omega$ is
actually a metric form. If some leaf $\call_{z^0}$ of $\call$
contains an essential vanishing cycle then it contains also an
imbedded essential vanishing cycle.
\end{corol}

\begin{rema}\rm
The same is true for disc-convex foliated manifolds $(X,\call)$
provided $X$ is a total space of an elliptic fibration (with
possibly singular fibers) over a disc-convex K\"ahler manifold
(apply Proposition 4.3).

\end{rema}

\newprg[prgPCFS.dim-two]{Imbedded shells in dimension two}

It would be instructive to understand something to the very end.
Also it is a time to get more examples and see how restrictive is
the presence of a foliated shell in $(X,\call)$. That's why let us
look closely to foliations on compact complex surfaces. $X$ in this
subsection will denote a compact complex surface, \ie a complex
manifold of dimension two. $\call$ will be a singular holomorphic
foliation by curves on $X$. We will work only with
$(X,\call)\in\cals$ in this subsection.

\smallskip
As we know on a compact complex surface there always exists a
$dd^c$-closed metric form. This was for the first time observed by
Gauduchon in \cite{Ga}.  Moreover all compact complex surfaces are
almost Hartogs, this is explained in \cite{Iv1, Iv4}. Really, the
K\"ahler ones are simply Hartogs, elliptic ones are served by
Proposition 4.3 and that of class $VII$ by the Proposition 4.2.
Therefore results of this paper are applicable to compact complex
surface in their full scale. Our task here is simple: to get
consequences from the presence of shells. This can be done using the
following beautiful and extremely powerful idea (I call it a
``pseudoconvex surgery'') due to Kato, we shall step by step use his
results from \cite{K1,K2} adapting them to our ``foliated'' case.

\begin{figure}[h]
\centering
\includegraphics[width=3.7in]{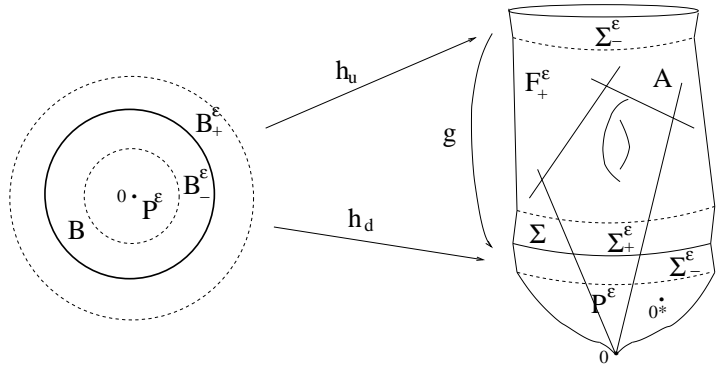}
\caption{Pseudoconvex surgery.} \label{surgery-fig}
\end{figure}

\smallskip\noindent{\sl Pseudoconvex surgery.}
Let $h:(P^{\eps}\setminus\{0\},\call^{\v})
\to (X,\call)$ be an imbedded foliated shell. We keep the notations
of the Introduction and of the proof of Theorem \ref{imbshell}.
Recall that $P=\cup_{z\in G}\Delta_z$ for a domain $G\ni 0$. In an $\eps$-
neighborhood $B^{\eps}$ of the boundary $B=\d P$ the mapping $h$ is
a foliated imbedding (but it is only immersion on the whole of
$P\setminus \{0\}$). The origin $\{0\}$ is the only essential
singular point of $h$. $\gamma_0\deff h(\d\Delta_0)$ is an essential
vanishing cycle. Set $\Sigma\deff h(B)$. Denote by $B^{\eps}_{\pm}$
one sided neighborhoods of $B$. Set $\Sigma^{\eps}_{\pm} =
h(B^{\eps}_{\pm})$ as on the Figure \ref{surgery-fig}. Cut $X$ along
$\Sigma$ to get a connected open set $E\deff X\setminus\Sigma$.

\smallskip\noindent Construct a pseudoconvex manifold $F^{\eps}_+$ by
gluing to $E$ the domain $P^{\eps}$ by  the biholomorphism
$h_d:B^{\eps}_+\to\Sigma^{\eps}_+$ - a copy of $h$ (in fact it may
have one cyclic singularity). Note that $F^{\eps}_+$ inherits the
foliation $\call$. Moreover, $F^{\eps}_+$ contains two copies of
$\Sigma^{\eps}_-$, one near its boundary - second in the interior
(see our Figure). There is a natural map $g$ between these two
copies of $\Sigma^{\eps}_-$, we refer to \cite{K1}, \S 1 for the
construction of $g$. For us it will be important that $g$ is a
foliated biholomorphism in its domain of definition. Really, $g$
comes in \cite{K1} and \cite{K2} as a part of a deck transformation
$\tilde g$ of a certain unramified covering $\tilde X$. The latter
inherits a foliation $\tilde\call$ which, of course, must be
preserved by $\tilde g$ and therefore is preserved by $g$. In fact,
one can see $g$ in our Figure: in coordinates on both copies of
$\Sigma^{\eps}_-$ in $F^{\eps}_+$ given by $h_u$ and $h_d$ the
mapping $g$ is the identity. $h_u$ is an ``upper'' copy of $h$.
Anyway, by the Hartogs extension theorem for holomorphic functions
$g$ extends onto the whole $F^{\eps}_+$ as a foliated holomorphic
map $g:(F^{\eps}_+,\call)\to (F^{\eps}_+,\call)$.

\smallskip $F^{\eps}_+$ also contains a point - the ``origin'' -
it comes from  the origin $0$ in $P^{\eps}$ when attaching it to
$E$. We keep noting this point as $0$.

\smallskip\noindent{\sl Claim 1. (M. Kato, Lemma 1 in \cite{K1},
Lemma 2 in \cite{K2}.) There exists a point $0^*\in F^{\eps}_+$ such
that
\begin{equation}
\bigcap_{n\ge 1} g^n(F^{\eps}_+)=\{0^*\}.
\end{equation}
}

\smallskip Let $A$ be the maximal compact subvariety of
$F^{\eps}_+$. Note that $A$ is contracted by $g$ to points. The set
$A$ on our Figure \ref{surgery-fig} is drawn as a chain of four
segments. We need to distinguish two cases.

\smallskip\noindent{\sl Case 1. $0^*\not\in A$.}

 $0^*$ may coincide with $0$ or not, we treat both
cases simultaneously. Remark that due to the fact that $g$ is
foliated and $0^*\not\in A$ it is a biholomorphism in a neighborhood
of $0^*$. Take a cyclic quotient $\bb_{l,d}$ of the standard ball
$\bb = \{z\in\cc^2:\norm{z}<1\}$ centered at $0^*$ and contained in
some $g^{n_0}(F^{\eps}_+)$ if $0^*=0$ and our shell was
$(l,d)$-cyclic shell. If $0^*\not= 0$ or $d=1$ then it is just the
ball.

\smallskip Now as in \cite{K2}, Lemma 5 one proves that, if $d=1$,
then:

\smallskip

\sli $g$ is a contracting biholomorphism in a neighborhood of $0^{*}$.

\slii Moreover there exists a strongly plurisubharmonic function
$\phi$ near $0^{*}$ such that for every $c>0$ small enough $P_c\deff
\{z:\phi (z)<c \}$ is biholomorphic to $\bb$ and $ g$ contracts each
$P_c$, \ie $ g(P_c)\comp P_c$.

\smallskip The proof if entirely local. Therefore if $d>1$ one lifts
$g$ from $\bb_{l,d}$ to $\bb$ and has the same properties for the
lifted local biholomorphism. In the first case one gets a primary
Hopf surface in the second - non-primary. Namely, one masters from
the shell between $\d P_c$ and $g^n(\d P_c)$ (for appropriate $n$
and $c>0$) a Hopf surface and proves that our surface $X$ blows down
to this one, call it $Y$ and the foliation obtained in $Y$ denote by
$\calf$. Bimeromorphic transformations/unramified coverings do not
effect essential (!) vanishing cycles. In particular our leaf
$\call_0$ with a vanishing cycle $\gamma_0$ descends to the same in
$Y$ with foliation $\calf$. Now let us see what happens in $Y$. Our
foliation is vertical in an appropriate coordinates near $0^*$.
Therefore, after appropriate change in $z_1$-coordinate can write
the contracting map $g$ (or its lift) in $\bb$ in the following
form:
\begin{equation}
\label{g}
g(z_1,z_2) =  (\alpha_1z_1, \alpha_2z_2+z_1g_1(z_1,z_2)),
\end{equation}
where $0<|\alpha_1|,|\alpha_2|<1$. Now it is obvious  that it is the
central fiber $\calf_0$ of $\calf$ which carries an essential
vanishing cycle and this fiber is a torus.

\smallskip\noindent{\sl Case 2. $0^*\in A$.}

Take a connected component of $A$ containing $0^*$ and from this
moment denote it as $A$. Let $\lambda :\tilde F^{\eps}_+\to
F^{\eps}_+$ be the minimal resolution of singularities of
$F^{\eps}_+$ and let $B$ be the proper preimage of $0$. Remark that
in the case of a cyclic quotient singularity all components of $B$
are rational curves, see pp. \cite{BHPV} 107-110. Consider an, a
priori meromorphic mapping $\tilde g\deff\lambda^{-1}\circ g\circ
\lambda : \tilde F^{\eps}_+\to \tilde F^{\eps}_+$.  Kato proved in
\cite{K2} that:

\smallskip
\sli $\tilde g$ is holomorphic (and foliated in our case) and there
exists $n$ such that $\tilde g^n(A\cup B)= \{\text{point}\}$.

\slii There exists $0^{**}\in B$ such that $\bigcap_{n\ge 1}\tilde
g^n(\tilde F^{\eps}_+)=\{0^{**}\}$.

\sliii There exists a strongly plurisubharmonic function $\phi$ near
$0^{**}$ such that for every $c$ small enough $P_c\deff \{z:\phi
(z)<c \}$ is biholomorphic to $\bb$ and $\tilde g$ contracts each
$P_c$, \ie $\tilde g(P_c)\comp P_c$.

\sliv $\tilde g: \tilde g^{-1}(P_c\setminus \{0^{**}\} \to P_c\setminus
0^{**}$ is a biholomorphism.

\smallskip\noindent See again Lemma 5 in \cite{K2}. Kato then masters
from these data (in a clear way) a surface $Y$ with Global Spherical
Shell in the terminology of Kato, or a Kato surface (his shell is
clearly foliated in our sense) and proves that our $X$ blows down to
$Y$ (as well as foliation $\call$ goes down to some $\calf$). On $Y$
one gets a divisor $C$ as factor of $(A\cup B)\setminus \tilde
g^n(F^{\eps}_+)$ by $\tilde g^n$ for an appropriate $n$. $C$ is
proved to be a chain (or two chains) of rational curves. Again the
foliation in a neighborhood of $0^{**}$ is vertical in an
appropriate coordinates. The image of the leaf $\call_0$ which
supports a vanishing cycle under $\tilde g^n$ cannot miss the set
$\tilde g^n(F^{\eps}_+)\cap A\cup B$, otherwise the corresponding
$\calf_0$ would not contain a vanishing cycle - this was already
once explained. Therefore $\call_0\subset C$. I.e. it is contained
in a rational curve and we are done.

\begin{rema}\rm
Corollary \ref{dim-2} from the Introduction is proved.
\end{rema}

It would be instructive to see clearly an example of the Case 2.
Let's take the simplest one.

\begin{figure}[h]
\centering
\includegraphics[width=2.7in]{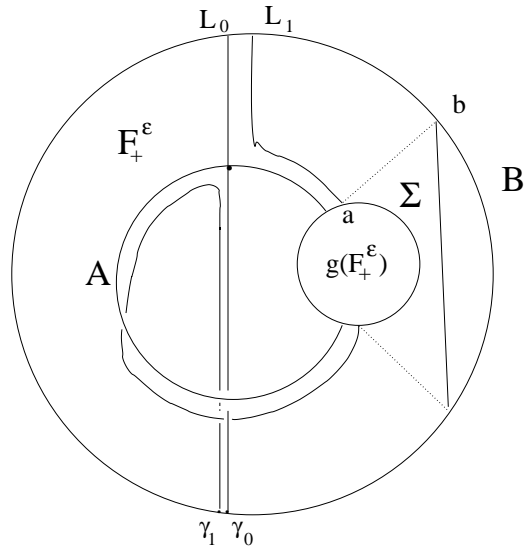}
\caption{Example to the case $0^*\in A$. $B$ is the standard sphere.
$F^{\eps}_+$ is the one time blown up unit ball. $\tilde g$ is given
by $(z_1,z_2)\to (\frac{1}{2}z_1,\frac{1}{2}z_2)$ and  $\Sigma$ is
the image of $B$ under $\tilde g$. The image $\tilde g(F^{\eps}_+)$
of  the blown up ball is removed and $X$ is obtained by identifying
$B$ with $\Sigma$. $\call_0$ lands to $A$ which is a rational curve
with one point of selfintersection. $\gamma_0$ is a circle (a point
on this Figure) and $\gamma_1$ on a nearby leaf $\call_1$
(complicated curve on the Figure) bounds a disc. To understand this
note that circles $a$ and  $b$ should be identified.}
\label{shel_kato-fig}
\end{figure}

\newsect[sectPLEX]{Pluriexact Manifolds and Foliations}

\newprg[prgPLEX.plu-exa]{Characterization of pluriexact foliations}

This Section is entirely devoted to pluriexact foliations. We start with the
following:

\begin{lem}
If $X$ is a compact complex manifold then $dd^c:\cale^{\rr}_{2,2}\to
\cale^{\rr}_{1,1}$ has closed range.
\end{lem}

\proof Observe the following resolution of the sheaf $\calh_{\cc}$
of complex valued pluriharmonic functions

\begin{equation}\label{exact-seq}
0\longrightarrow \calh_{\cc} \buildrel (-\d , \id) \over
\longrightarrow \Omega^1\oplus \left[ \calh_{\rr}
+i\cale_{\rr}\right] \buildrel (\id\oplus \d )\over \longrightarrow
\cale^{1,0} \buildrel (\d \oplus \bar\d) \over \longrightarrow
\cale^{1,1}_{\rr}\buildrel dd^c \over \longrightarrow
\cale_{\rr}^{2,2}\buildrel d \over \longrightarrow ...
\end{equation}

Here $\calh_{\rr}$ is the sheaf of real valued pluriharmonic
functions, $\Omega^1$ the sheaf of holomorphic $1$-forms,
$\cale_{\rr}$ the sheaf of smooth real valued functions. This
resolution tells that

\begin{equation}
\ker\{d:\cale_{\rr}^{2,2}\to\cale_{\rr}^3\}/\im\{dd^c:\cale_{\rr}^{1,1}
\to\cale_{\rr}^{2,2}\}\equiv H^4(X,\calh_{\cc}).
\end{equation}
Therefore $dd^c:\cale_{\rr}^{1,1}\to\cale_{\rr}^{2,2}$ has closed
range (in fact of finite codimension). By duality
$dd^c:\cale^{\rr}_{2,2}\to\cale^{\rr}_{1,1}$ has also closed range.

\qed

Fix some strictly positive $(1,1)$-form $\Omega$ on $X$. Let $\call$
be a holomorphic foliation by curves on $X$. Denote by
$K_{1,1}(\call)$ the compact set in $\cale_{1,1}^{\rr}$ which
consists from positive $(1,1)$-currents $T$ tangent to $\call$ such
that $(\Omega ,T)=1$, \ie the compact base of currents directed by
$\call$. Let us prove now the Proposition \ref{harv-laws} from
the Introduction, it is analogous to Theorem 3.18 from
\cite{Go}, a non-foliated version for $dd^c$-closed metric forms was
given in \cite{Iv1}.

\smallskip\noindent{\slsf Proof of Proposition \ref{harv-laws}.} Let
$\omega$ be a pluriclosed taming form for $\call$. If $dd^cS=T\in K_{1,1}(\call)$ for
some $S\in \cale_{2,2}^{\rr}$ then $0<(\omega ,T)=(\omega
,dd^cS)=(dd^c\omega ,S)=0$ - a contradiction. Vice verse, if
$K_{1,1}(\call)\cap dd^c\cale_{2,2}^{\rr}=\emptyset$ then, since
$dd^c\cale_{2,2}^{\rr}$ is closed, by Hanh-Banach theorem there
exists $\omega$ such that $\omega|_{K_{1,1}(\call)}>0$ and
$\omega|_{dd^c\cale_{2,2}^{\rr}}=0$.

\qed

\newprg[prgPLEX.forms]{Plurinegative metric and taming forms}

A form $\omega\in \cale^{p,p}$ is positive if its restriction onto
any germ of $p$-dimensional complex submanifold is positive (meaning
$\ge 0$). This is equivalent to the positivity of $(n,n)$-forms
\[
\omega\wedge i\theta_1\wedge\bar\theta_1\wedge ...\wedge
i\theta_{n-p}\wedge\bar\theta_{n-p}
\]
for all $(1,0)$-forms $\theta_1,...,\theta_{n-p}$. A current $T\in
\cale_{p,p}^{\rr}$ is positive if $\big< T,\omega\big>\ge 0$ for
every positive $w\in\cale^{p,p}$.  Denote by $K_{p,p}$ the compact
set in $\cale_{p,p}^{\rr}$ which consists from such strictly
positive $(p,p)$-currents $T$ on $X$ that $<T,\Omega^{p}>=1$.

\begin{prop}
\label{harv-laws-neg}
Let $X$ be a compact complex manifold. Then the following alternative
holds true:

\sli either $X$ admits a plurinegative metric form,

\slii or, there exists a sequence $S_n\in K_{2,2}$ and  increasing
sequence of real numbers $t_n$ such that $t_ndd^cS_n$ converges to
some $T\in K_{1,1}$.
\end{prop}
\proof (\sli $\Rightarrow$ (\slii Let $\omega$ be a plurinegative
metric form. Then we need to see that such sequence cannot exist.
Suppose it does. Then
\[
0 \ge <t_nS_n,dd^c\omega> = <t_ndd^cS_n, \omega> \to <T,\omega> > 0.
\]
Therefore $<T,\omega> = 0$, contradiction.

\smallskip\noindent (\slii $\Rightarrow$ (\sli   Set
$B=dd^c(K_{2,2}$. This is a convex compact in
$\cale_{1,1}^{\rr}(X)$. Let $K_B$ be the cone generated by $B$.
Nonexistence of a sequence as in (\slii means exactly that $K_B\cap
K_{1,1}=\emptyset$. Hanh-Banach theorem gives us a continuous linear
form $\omega$ on $\cale_{1,1}^{\rr}(X)$ such that $\omega
|_{K_{1,1}} >0$, \ie $\omega$ is a metric form, and such that
$\omega|_B\le 0$, \ie $\omega$ is plurinegative.

\smallskip\qed

\smallskip

Let $\call$ be a holomorphic foliation by curves on $X$. Denote by
$K_{1,1}(\call)$ the compact set in $\cale_{1,1}^{\rr}$ which
consists from positive $(1,1)$-currents $T$ tangent to $\call$ such
that $(\Omega ,T)=1$, \ie the compact base of currents directed by
$\call$. The following Proposition is analogous to the previous one and
the proof is identic to that already given.

\begin{prop}
\label{harv-laws-fol}
Let $(X,\call)$ be a compact foliated pair. Then the following alternative
holds true:

\sli either $(X,\call)$ admits a plurinegative taming form,

\slii or, there exists a sequence $S_n\in K_{2,2}$ and  increasing
sequence of real numbers $t_n$ such that $t_ndd^cS_n$ converges to
some $T\in K_{1,1}(\call)$.
\end{prop}

\newprg[prgPLEX.divis]{Subdivision of pluriexact manifolds and foliations}

\smallskip
Based on what was said in this section we divide the class $\cale$
of compact pluriexact (foliated) manifolds into three classes:

\begin{itemize}
 \item Class $\cale_{-}$ of (foliated) manifolds from $\cale$ admitting a
 plurinegative (taming) form.

 \item Class $\cale_{+}$ of (foliated) manifolds from $\cale$ admitting a
 (non-trivial) positive (directed) $(1,1)$-current $T$ such that $T=dd^cS$
 for some positive $(2,2)$-current $S$.

 \item Class $\cale_0\deff\cale\setminus\big(\cale_{-}\sqcup\cale_{+}\big)$.
\end{itemize}

\smallskip Let us give an example of a foliated manifold from class
$\cale_0$. The example is the classical one due to Hironaka. In
nonhomogeneous coordinates $(x,y,z)$ in $\cc\pp^3$ consider the
rational curve $C$ with exactly one transverse self-intersection,
which is defined by the equation
\begin{equation}
\eqqno(curve)
y^2 = x^2 + x^3, \qquad z=0.
\end{equation}
Manifold $X$ of the example is a proper modification of $\cc\pp^3$ along
$C$. In an neighborhood of the origin one blows-up first one local
irreducible branch of $C$ and then an another. Let $p:\cc\pp^1\to
C\subset\cc\pp^3$ be a holomorphic parameterizing map. Let
$D_{\infty}$ be a disc in $\cc\pp^1$ near $\infty$ and let
$p(D_{\infty})$ be its image. We perform a local blow-up of
$\cc\pp^3$ with the center $p(D_{\infty})$ and denote by
$C^{'}_{\infty}$ the strict transform of the origin. Then we blow-up
again along $p(D_{0})$, where $D_0$ is a disc around $0$ n
$\cc\pp^1$. By $C_{\infty}$ we denote the strict transform of
$C^{'}_{\infty}$ under this second blow-up, by $C_{0}$ the strict
transform of the point of intersection of $p(D_{0})$ with the
exceptional divisor of the first blow-up. We perform the blow-up
along the remaining part of $C$ and denote by $E$ the exceptional
divisor of the resulting blow-down map $\pi :X\to\cc\pp^3$. As it is
well known (and obvious), see \cite{S} for example, $[C_0]$ is
homologous to  $[C_0]+[C_{\infty}]$ and therefore $C_{\infty}$ is
homologous to zero.

\medskip\noindent{\slsf Note. In fact $[C_{\infty}]$ is $dd^c$-exact
as a current.} To see this take as $S=\frac{1}{\pi}\pi^*\ln
|t|^2\cdot [E]$ where $t$ is a parameter on the parameterizing curve
$\cc\pp^1$. This means that we consider $\frac{1}{\pi}\pi^*\ln
|t|^2$ as a function on $E$ and the action of $S$ on $(2,2)$-form
$\phi$ is given by
\begin{equation}
\eqqno(currS)
<S,\phi> = -\frac{1}{\pi}\int_E\pi^*\ln |t|^2\phi .
\end{equation}
Since $dd^c\ln |t|^2 = \pi\delta_{\{0\}}$ we get that for any $(1,1)$-form
$\psi$
\[
<dd^cS,\psi> = <S, dd^c\psi> = -\frac{1}{\pi}\int_E\pi^*\ln |t|^2dd^c\psi =
-\frac{1}{\pi}\int_Edd^c\big(\pi^*\ln |t|^2\big)\psi =
\]
\begin{equation}
\eqqno(ddc-exact)
= -\frac{1}{\pi}\int_E\pi^*dd^c\big(\ln |t|^2\big)\psi = \int_{C_{\infty}}\psi
= <[C_{\infty}],\psi>.
\end{equation}
Let us see that $[C_{\infty}]$ is not a $dd^c$ of a {\slsf positive}
$(2,2)$-current.

\begin{lem}
\label{non-exact}
There exists no {\slsf positive} $(2,2)$-current $S$ such that $dd^cS$ is
positive and non-zero.
\end{lem}
\proof Suppose the opposite, \ie that there exists positive
$(2,2)$-current $S$ on $X$ such that $T\deff dd^cS$ is a non-zero
positive $(1,1)$-current. Then for positive currents $\tilde
S\deff\pi_*S$ and $\tilde T\deff\pi_*T$ one has $dd^c\tilde S =
\tilde T$.  Since $\cc\pp^3$ is K\"ahler this implies that $\tilde
T=0$. Since $\pi : X\setminus E\to \cc\pp^3\setminus C$ is a
biholomorphism the current $T$ is supported on $E$.

\smallskip By the Cut-off Theorem of Bassanelli $\chi_E\cdot S$ is a
$dd^c$-positive current supported on $E$,  \ie $\chi_E\cdot S =
h[E]$ where $h$ is a plurisubharmonic function on $E$. But then $h$
is constant and $dd^c\chi_E\cdot S = 0$.  Now we can write
\begin{equation}
\eqqno(line1)
T=dd^cS = dd^c(\chi_E\cdot S) + dd^c(\chi_{X\setminus E}\cdot S) =
dd^c(\chi_{X\setminus E}\cdot S).
\end{equation}
Let us see that $\chi_E\cdot dd^c(\chi_{X\setminus E}\cdot S)=0$.
Since $\chi_E\cdot T = T$, from this will follow that $T=0$, which
is a contradiction. By one more theorem of Bassanelli, see Theorem
3.5 in \cite{Bs} $dd^c(\chi_{X\setminus E}\cdot S) -
\chi_{X\setminus E}dd^cS$ is a negative current supported on $E$.
Therefore $\chi_Edd^c(\chi_{X\setminus E}\cdot S)$ is negative and
supported on $E$. But from \eqqref(line1) we see that this current
is positive. So it is zero.

\smallskip\qed

\begin{lem}
\label{sequence} There exists a sequence of positive  bidimension
$(2,2)$ currents $\{ S_n\}$ of mass one and an increasing sequence
of positive real numbers $t_n\nearrow +\infty$ such that
\[
 [C_{\infty}] = \lim_{n\to\infty}t_ndd^cS_n.
\]
\end{lem}
\proof For $p>0$ consider the following function on the Riemann
sphere $\cc\pp^1$
\begin{equation}
\eqqno(phi)
\phi_p (z) =
\begin{cases}
|z|^{2p} & \text{ if } |z|\le 1\cr
2 - \frac{1}{|z|^{2p}} & \text{ if } |z|\ge 1.
\end{cases}
\end{equation}
We have that
\begin{equation}
\eqqno(ddcphi1)
dd^c\phi_p = \frac{i}{2}\d\dbar \phi_p = \frac{i}{2}\d (p|z|^{2p-2}zd\bar z)
= p^2|z|^{2p-2}\frac{i}{2}dz\wedge d\bar z
\quad \text{ for } |z|\le 1.
\end{equation}
Analogously in coordinate $w\deff 1/z$ we have
\begin{equation}
\eqqno(ddcphi2)
dd^c\phi_p = dd^c(2-|w|^{2p}) = -p^2|w|^{2p-2}\frac{i}{2}dw\wedge d\bar w
\quad \text{ for } |w|\le 1.
\end{equation}
For $p>0$ consider $(2,2)$-currents $S_p$ which are define as
\begin{equation}
\eqqno(currSp)
<S_p,\psi> = \frac{1}{\pi}\int_E\pi^*\phi_p\cdot\psi|_E
\end{equation}
for a $(2,2)$-form $\psi$ on $X$. Note that $S_p$ are positive. Set
$t_p=\frac{1}{\pi p}$. We want to prove that $t_pdd^cS_p\to
[C_{\infty}]$ as  $p\searrow 0$.

\smallskip First remark that $\big< t_pdd^cS_p, \psi\big> = \big<t_pdd^c\phi_p,\psi \big>$
for a test $(1,1)$-form $\psi $. Here $dd^c\phi_p$ should be
understood in the sense of distributions. From \eqqref(ddcphi1)  we
see that
\[
t_p\int\limits_{\Delta}dd^c\phi_p =1, \quad dd^c\phi_p\ge 0 \text{ and } t_pdd^c\phi_p =
\frac{p}{\pi}|z|^{2p-2}\frac{i}{2}\to 0
\]
as $p\searrow 0$ on compacts of $\Delta^*\deff\Delta\setminus
\{0\}$. This implies that $t_pdd^c\phi_p\to \delta_{\{0\}}$ in
$\Delta$. In the same manner from \eqqref(ddcphi2) we get that
$t_pdd^c\phi_p\to - \delta_{\{\infty\}}$ in $\cc\pp^1\setminus
\bar\Delta$.

\smallskip\noindent{\slsf Claim. $t_pdd^c\phi_p\to 0$ on the annulus $A\deff\Delta
(2)\setminus \bar\Delta (1/2)$.}

\smallskip Take a test function $\psi$ with support in $A$, set $\tilde\psi (r) =
\int_0^{2\pi}\psi (r,\theta)d\theta$ and
write
\[
4t_p\int\limits_Add^c\phi_p\cdot \psi \deff
t_p\int_{1/2}^2\int_0^{2\pi}\phi_p(r) \Delta\psi rd\theta dr =
t_p\int_{1/2}^2\int_0^{2\pi}\phi_p(r)\big(r\frac{\d^2 \psi}{\d r^2}
+ \frac{1}{r}\frac{\d^2 \psi}{\d\theta^2} + \frac{\d \psi}{\d r}
\big) d\theta dr =
\]
\[
= t_p\int_{1/2}^2\phi_p(r)\big(r\tilde\psi^{''} + \tilde\psi^{'}\big)dr =
t_p\big[\int_{1/2}^1r^{2p+1}\tilde\psi^{''}dr  +
\int_1^2\big(2 - \frac{1}{r^{2p}}\big)r\tilde\psi^{''}dr\big] +
\]
\[
t_p\big[\int_{1/2}^1r^{2p}\tilde\psi^{'}dr  + \int_1^2\big(2 - \frac{1}{r^{2p}}\big)
\tilde\psi^{'}dr\big] =: I^{''}_p + I^{'}_p.
\]
Next we compute these integrals separately.
\[
I^{'}_p = t_p\big[\int_{1/2}^1\big(r^{2p}\tilde\psi\big)^{'}dr  +
\int_1^2\big(\big(2 -
\frac{1}{r^{2p}}\big)\tilde\psi\big)^{'}dr\big] - \frac{2p}{\pi
p}\int_{1/2}^1r^{2p-1}\tilde\psi dr - \frac{2p}{\pi
p}\int_{1}^2r^{-2p-1}\tilde\psi dr =
\]
\begin{equation}
\eqqno(one) = \frac{\tilde\psi (1)}{\pi p} - \frac{\tilde\psi
(1)}{\pi p} - \frac{2}{\pi } \big(\int_{1/2}^1r^{2p-1}\tilde\psi dr
+ \int_{1}^2r^{-2p-1}\tilde\psi dr\big) \to -
\frac{2}{\pi}\int_{1/2}^2\frac{\tilde\psi}{r}dr
\end{equation}
as $p\searrow 0$. At the same time
\[
I^{''}_p = t_p\int_{1/2}^1\big(r^{2p+1}\tilde\psi^{'}\big)'dr  +
t_p\int_1^2 \big(\big(2 -
\frac{1}{r^{2p}}\big)r\tilde\psi^{'}\big)'dr - \frac{2p+1}{\pi
p}\int_{1/2}^1r^{2p}\tilde\psi'dr -
\]
\[
- \frac{1}{\pi
p}\int_1^2\big(2+\frac{2p-1}{r^{2p}}\big)\tilde\psi'dr =
\frac{\tilde\psi' (1)}{\pi p} - \frac{\tilde\psi' (1)}{\pi p} -
\frac{2p+1}{\pi p}\int_{1/2}^1\big(r^{2p}\tilde\psi\big)'dr -
\frac{1}{\pi
p}\int_1^2\big(\big(2+\frac{2p-1}{r^{2p}}\big)\tilde\psi\big)'dr
\]
\[
+ \frac{2(2p+1)}{\pi}\int_{1/2}^1r^{2p-1}\tilde\psi dr -
\frac{2p(2p-1)}{\pi p}\int_1^2\frac{\tilde\psi}{r^{2p+1}}dr \to
-\frac{2p+1}{\pi p} \tilde\psi (1) + \frac{2p+1}{\pi p}\tilde\psi
(1) +
\]
\begin{equation}
\eqqno(two)
+ \frac{2}{\pi}\big(\int_{1/2}^1\frac{\tilde\psi}{r}dr  + \int_1^2
\frac{\tilde \psi}{r}dr\big) = \frac{2}{\pi}\int_{1/2}^2
\frac{\tilde\psi}{r}dr.
\end{equation}
\eqqref(one) and \eqqref(two) cancel each other and the Claim is
proved. We conclude that
\begin{equation}
\eqqno(ddcphi3) t_pdd^c\phi_p\to \delta_{\{0\}} -
\delta_{\{\infty\}}
\end{equation}
on $\cc\pp^1$. From \eqqref(currSp) we conclude immediately that
$t_pdd^cS_p\to [C_{\infty}]$.

\smallskip\qed

\begin{rema}\rm
Lemmas \ref{non-exact} and \ref{sequence} show that the manifold in
the example of Hironaka belongs to class $\cale_0$, \ie that
$\cale_0\not=\emptyset$. In order to get a holomorphic foliation  by
curves  $\call$ on $X$ such that  the foliated manifold $(X,\call)$
belongs to $\cale_0$ it is sufficient to take $\call$ in such a way
that the current $T\deff[C_{\infty}]$ is directed by $\call$. The
rational fibration $\{ x=y=\const\}$ will do the job.
\end{rema}

\newprg[prgPLEX.pl-neg]{Foliations with Plurinegative Taming Forms}

In this section we consider the class $\cale_{-}$ of foliated
manifolds $(X,\call)$ which possess a plurinegative taming form
$\omega$. Our goal is to understand the vanishing cycles in $\call$.
We start with the following:

\begin{prop}
\label{exact-curr} If a leaf $\call_z$ of a disc-convex foliated
manifold $(X,\call,\omega)\in\cale_-$ contains an essential vanishing cycle
then there exists a nontrivial positive $d$-exact $(1,1)$-current
$T$ tangent to $\call$. The support of $T$ is contained in
$\overline{\call_z}$.
\end{prop}
\proof Let $\gamma\subset \call_z^0$ be (an essential) vanishing
cycle. From the proofs of Theorems \ref{reparalldim1} and
\ref{immshell} we see that there exists a foliated meromorphic
immersion $h:(\Delta^{n+1}, \call^{\v})\setminus S \to (X,\call)$
such that $h(\d \Delta_0) = \gamma$ and the essential singularity
set $S$ of $h$ intersects $\Delta_0$. From (\ref{areainfty}) we see
that there exist $q_n\to 0$ in $\Delta^n$ such that
$h|_{\d\Delta_{q^n}}$ converges but $\area (h(\Delta_{q^n}))$
diverges to infinity. Therefore by standard (and obvious) reasoning
currents
\begin{equation}
T_n = \frac{[h(\Delta_{q^n})]}{\area (h(\Delta_{q^n}))}
\end{equation}
converge to a closed, positive current $T$ of mass one  tangent to
$\call$.

\smallskip To prove that $T$ is, in fact, exact, observe that by Lemma
2.2 from \cite{Iv3} one have that $H^2_{DR}(\Delta^2\setminus S) =
0$, where $H_{DR}$ denotes the de Rham cohomology. Let $\phi$ be a
$d$-closed $2$-form on $X$. Then $h^*\phi = d\psi$ for some $1$-form
in $\Delta^2\setminus S$. Therefore
\[
<T,\phi> = \lim_{n\to \infty}<T_n, \phi> = \lim_{n\to
\infty}\frac{1}{\area (h(\Delta_{q^n}))}
\int\limits_{\Delta_{q^n}}h^*\phi = \lim_{n\to \infty}\frac{1}{\area
(h(\Delta_{q^n}))} \int\limits_{\d\Delta_{q^n}}h^*\psi =
\]
\[
= \int\limits_{\d\Delta_{0}}h^*\psi\lim_{n\to \infty}\frac{1}{\area
(h(\Delta_{q^n}))} = 0.
\]

\smallskip\qed

\smallskip Now we turn to the

\smallskip\noindent{\slsf Proof of Theorem \ref{plurineg}.} Let $z_0\in X^0$
be a point
such that the leaf $\call^0_{z^0}$ of the foliation $\call$ through
$z^0$ contains a vanishing cycle $\gamma_0$. According to subsection
3.2 the situation can be reduced to the following. Take a transverse
$n$-disc $D$ through $z^0$, $n+1$ being the complex dimension of
$X$, in such a way that after identification of $D$ with
$\Delta^{n}$ the point $z^0$ corresponds to $0$. Let $\hat\call_D$
be  the completed holonomy covering cylinder of $\call$ over $D$.
Then there exists an  imbedded vanishing cycle $\hat\gamma_0 \subset
\hat\call_{z^0}$ which projects to (some another) vanishing cycle
$\gamma_0\subset \call_{z^0}$ under the natural meromorphic
projection $p:\hat\call_D\to X$. This $\gamma_0$ may not be an
imbedded loop and may be different from that one which was taken at
the beginning.

\smallskip Shrinking $D$, if necessary, we can construct a
generalized Hartogs figure $(W,\pi,U,D)$ over $D$, where $W$ is an
open subset of $\hat\call_D$, $U$ is a neighborhood of some point
$z\in D$, $\pi :W\to D$ is the restriction of the natural projection
$\pi : \hat\call_D\to D$ to $W$.

\smallskip According to  Theorem \ref{reparalldim1} we get that
our  meromorphic foliated immersion $p:(W,\pi,U,D)\to
(X,\call,\omega)$ extends, after a reparametrization, to a
meromorphic foliated immersion $\tilde p :\widetilde{W}\setminus
S\to X$, where::

\begin{itemize}

\item
$(\widetilde{W},  \pi , D) $ is a {\slsf complete} Hartogs
figure over $D$. That means the complex manifold $\widetilde{W}$ is
holomorphically foliated by discs by the holomorphic submersion
$\pi : \widetilde{W}\to D$.

\smallskip
\item $S$ is a closed subset of $\widetilde{W}$ which has the
following structure: there exists a closed, complete $(n-1)$-polar
subset $S_1\subset\Delta^{n}$ such that $S=\bigcup_{s\in S_z}S_z$
where $S_z$ is a compact subset of $\widetilde{W}_z \deff
\pi^{-1}(z)$.

\smallskip
\item $(n-1)$-polarity of $S_1$ means that $\Delta^{n}$ in its turn
(after shrinking, if necessary) can be decomposed as
$\Delta^{n-1}\times\Delta$ in such a way that for every
$\lambda\in\Delta^{n-1}$ the intersection $S_{1,\lambda}\deff
\big(\{\lambda\}\times\Delta\big)\cap S_1$ is a complete polar (and
compact in our case) subset of $\Delta$.
\end{itemize}

\medskip Now we have the following two possible cases.

\medskip\noindent{\slsf Case 1. For every $\lambda\in\Delta^{n-1}$
the set $S_{1,\lambda}$ is not empty.}

\smallskip Take a point $s_{\lambda}\in S_{1,\lambda}$ for every
$\lambda\in\Delta^{n-2}$. From Theorem \ref{repar2dim} it readily
follows that for every $\lambda$ there exists a sequence
$(\lambda_n,z_{1,n}) \to s_{1,\lambda}$ such that
\begin{equation}
\eqqno(to-infty)
\area \big[\tilde p(\Delta_{\lambda_n,z_{1,n}})\big]\to\infty .
\end{equation}
As it was explained in Proposition \ref{exact-curr} such sequence accumulates
to an exact, positive $(1,1)$-current $T_{\lambda}$ directed by $\call$.
Remark that the support of each $T_{\lambda}$ belongs to
$\call_{\Delta^2_{\lambda}}$ and therefore they are distinct for
different $\lambda$-s.

\medskip\noindent{\slsf Case 2. For some $\lambda_0\in\Delta^{n-1}$
the set $S_{1,\lambda_0}$ is empty.}

\smallskip Then $S_{1,\lambda} =\emptyset$ for $\lambda$ in a neighborhood
of $\lambda_0$. We can assume that $\lambda_0$ is as close to $0$ in
$\Delta^{n-1}$ as we wish. Otherwise for some neighborhood of $0$
the Case 1 will occur. Restricting $\widetilde{W}$ to a smaller
polydisc, if needed, we find ourselves in the conditions where
$(\widetilde{W},  \pi , D)$ is isomorphic to $(\Delta^{n+1}, \pi ,
\Delta^{n})$ where $\pi:\Delta^{n+1}\to \Delta^{n}$ is the canonical
projection. We find ourselves in the assumptions of \cite{IS}. More
precisely, we have a meromorphic mapping $\tilde p: H^n_2(\eps)\to
X$ where:

\smallskip
\sli $H^n_2(\eps) = \big[\Delta^{n-1}(\lambda_0,\eps)\times
\Delta^2\big]\cup \big[\Delta^{n-1} \times A^2(1-\eps , 1)\big]$ -
the Hartogs figure of bidimension $(n-1,2)$;

\slii $\tilde p$ is holomorphic on $\Delta^{n-1}\times A^2(1-\eps ,
1)$;

\sliii the image manifold $X$ admits a $dd^c$-closed positive $(2,2)$-form
$\omega_2$.

\smallskip By the result of \cite{IS} the mapping $\tilde p$
meromorphically extends to $\Delta^n\setminus R$, where the
singularity set $R$ is either empty, or has the following structure:

\smallskip
\sli $R=\cup_{\lambda\in R_1}R_{\lambda}$, where $R_1$ is complete
$(n-3)$-polar closed subset of $\Delta^{n-2}$ and each $R_{\lambda}$
is complete polar compact of $\Delta^2_{\lambda}$ of Hausdorff
dimension zero.

\slii Again, the $(n-3)$-polarity of $R_1$ means that in a
neighborhood of zero we can decompose
$\Delta^{n-1}=\Delta^{n-2}\times \Delta$ in such a way that for
every $\lambda' = (\lambda_1,...,\lambda_{n-2})$ the set
$R_{\lambda'}\deff R\cap \Delta^3_{\lambda'}$ is a compact,
pluripolar subset of Hausdorff dimension zero.

\sliii Moreover $\tilde p(\d\Delta^3_{\lambda'})$ is not homologous
to zero in $X$.

\medskip The last item means that $\tilde p(\d\Delta^3_{\lambda'})$
is a foliated three-dimensional shell in $X$ for every $\lambda'\in
\Delta^{n-2}$.

\smallskip Theorem \ref{plurineg} is proved.

\newprg[prgPLEX.exact]{Foliations with Pluriexact Directed Currents}

This subsection is devoted to the class $\cale_{+}$. By $\call^s \deff
\call^{\sing}$ we denote the singular set of $\call$. $\call^s$ has
complex codimension at least two.

\smallskip\noindent{\slsf Proof of Theorem \ref{pluripos}.}
Let $R$ be a positive $(2,2)$-current such that $T\deff dd^cR$ is
non-trivial, positive and directed by $\call$.  Exactly as in the
proof of Lemma \ref{non-exact} (taking $\call^s$ instead of $E$ and
$R$ as $S$ in the notations of the proof) we get that
$\chi_{\call^s}T = \chi_{\call^s}dd^c(\chi_{X\setminus\call^s}R)$ is
a negative current supported on $\call^s$. But $T$ is positive and
therefore $\chi_{\call^s}T = 0$. This proves the part (\sli of our
theorem.

\smallskip Now we shall prove the part (\slii. We fix a strictly positive
$(1,1)$-form $\Omega$ on $X$.
Masses of positive $(p,p)$-currents $R$ on $X$ of order zero will be
measured as $\norm{R}
\deff \big<R,\Omega^p\big>$. First remark that for every $(1,1)$-form
$\omega$ on $X$ such that $dd^c\omega + \Omega^2\ge 0$ one has
$\big< T,\omega \big> \ge - \norm{R}$. This is immediate, just write
\begin{equation}
\eqqno(normS) \big< T,\omega \big> = \big< dd^cR,\omega \big> =
\big< R ,dd^c\omega\big> \ge \big< R,-\Omega^2\big> = -\norm{R}.
\end{equation}

\smallskip\noindent Second, for $n\ge 3$ consider the following
$(n-2,n-2)$-form in $\cc^n$
\begin{equation}
\eqqno(omega0)
\omega_0 = -\frac{i}{2}\frac{\big(dz,dz\big)}{\norm{z}^2}\wedge
\big(dd^c\norm{z}^2\big)^{n-3}.
\end{equation}
One checks easily that:

\smallskip

\sli $dd^c\omega_0$ is positive;

\smallskip

\slii $dd^c\omega_0\in L^p_{loc}$ for every $p<\frac{n}{2}$.

\smallskip\noindent Both properties follow immediately from the following expression:
\[
dd^c\frac{i}{2}\frac{(dz,dz)}{\parallel z\parallel^2} =
\frac{i^2}{4\pi}\frac{(dz,dz)}{\parallel z\parallel^6}\wedge
(2(dz,z)\wedge(z,dz) - \parallel z \parallel^2(dz,dz)),
\]
where the form on the right hand side is non-positively definite.
Indeed, for a vector  $v=(v_1,...,v_n)\in\cc^n$ one has
\[
 (\parallel z \parallel^2(dz,dz)-2(dz,z)\wedge(z,dz))\wedge
 v\wedge\bar v = \parallel z \parallel^2\parallel
 v\parallel^2 - 2(v,z)(z,v)
 \]
 \[
 =\parallel z \parallel^2\parallel
 v\parallel^2 - 2|(v,z)|^2\geq 0
 \]
by Cauchy-Schwarz inequality.

\smallskip Taking $\rho\omega$ for an appropriate cut-off function
$\rho$ with support in the unit polydisc we can extend $\omega_0$ to
$\cc^n$ by zero and with $dd^c\omega_0\ge -\eta$ for $\eta > 0$ as
small as we need. Smoothing by convolution we get $\omega_{\eps}$
converging to $\omega_0$ with support in the unit polydisc and
having $dd^c\omega_{\eps}\ge -2\eta$. Now we can push
$\omega_{\eps}(z-\zeta)$ to $X$ verifying $dd^c\omega_{\eps}(z -
\zeta) + \Omega^2\ge 0$ for all $\zeta$ in a neighborhood of a fixed
center $p_0$ of some foliated coordinate chart.

\smallskip Let $B=\Delta^{n}\times \Delta$ be a foliated chart for $\call$
and $\mu$ the induced by $T$ Radon measure on $\Delta^{n}$. Coordinates in $B$
denote as $(z',z_{n+1})=(z_1,...,z_{n},z_{n+1})$. Write
\[
\big< T,\omega_{\eps}\big> \sim -\frac{i}{2}\int\limits_{\Delta^{n}}
\int\limits_{\Delta} \frac{dz_{n+1}\wedge d\bar z_{n+1}}{\norm{z'-\zeta'}^2 +
|z_{n+1}|^2}d\mu (z') = - 2\pi \int\limits_{\Delta^{n-1}}d\mu (z')
\int_0^1\frac{rdr}{r^2 + \norm{z'-\zeta'}^2} \sim
\]
\[
\sim \pi \int\limits_{\Delta^{n-1}}\ln \norm{z'-\zeta'}^2d\mu (z') \ge -\norm{S}
\]
by \eqqref(normS) and independently of $\zeta'\in \Delta^{n}$.
Theorem is proved.

\smallskip\qed

\begin{rema} \rm
\label{sing-size}
The singularity set $S$ repeatedly comes out in the proofs of this paper.
Let us make a few remarks about this issue.

\smallskip\noindent{\bf (a)} An appearance of $S$ is a highly {\slsf non-algebraic}
phenomenon. It is sufficient to say that due to the Hartogs type extension
theorem of \cite{Iv3} the set $S$ is always empty if $X$ is K\"ahler (or if
$(X,\call)$ admits a $d$-closed taming form).

\smallskip\noindent{\bf (b)} In the case of a pluriclosed taming form
$S$ is {\slsf proper} over a subspace of complex codimension two, this is due to the
homological nature of shells, see Lemma \ref{ext-shel} and Theorem \ref{reparalldim2}.

\smallskip\noindent{\bf (c)} If the taming form is only plurinegative the size
of $S$ can fall down, and this phenomena is responsible for the
appearance of three-dimensional shells, see Theorem \ref{plurineg}
and Example \ref{hopf-three}.
\end{rema}

\newsect[sect.EXOQ]{Other Results, Examples and Open Questions}

We still owe the proofs of some statements used in the text of this
paper and of some propositions from the Introduction. Moreover it is
the time to give more interesting examples (from the point of view
of this text) then just foliations on complex surfaces or on
K\"ahler manifolds. Looking on each example in this Section we shall
be rather attentive to its Hartogs properties because, as it should
be clear from the proofs of this paper, the failure of a foliated
manifold $(X,\call)$ to be Hartogs is "almost equivalent" to the
presence of essential vanishing cycles/foliated shells in $(X,\call
)$.

\newprg[prgEXOQ.har-fol]{Hartogs foliation on a compact non-Hartogs threefold}

The following example is due to Nakamura, see \cite{Na}. We only
interpret it according to our needs adding a foliation to it.

\begin{exmp}\rm
Take any matrix $A\in SL(2,\zz)$ with real eigenvalues $\alpha <1$
and $1/\alpha$. For example the following one:
\begin{equation}
A =
\left(%
\begin{array}{cc}
  1 & 1 \\
  1 & 2 \\
\end{array}%
\right).
\end{equation}
Here $\alpha = 3/2-\sqrt{5}/2$. Consider the standard integer
lattice $\Lambda_0\deff\zz^4$ in $\cc^2$. $A$ preserves $\Lambda_0$
and therefore defines a holomorphic automorphism $A$ of the torus
$\ttt_0\deff\cc^2/\Lambda_0$. Therefore we can construct a compact
complex threefold $X_0\deff \cc^*\times\ttt_0/<g>$ where $g(z,Z)=
(\alpha z, AZ)$. $X_0$ is a complex $2$-torus bundle over a complex
$1$-torus $\cc^*/<\alpha>$. We fix the coordinate $z$ for $\cc^*$.

\smallskip Let $\v$ be the eigenvector of $A$ with eigenvalue
$\alpha$ and $\w$ be that with $1/\alpha$. It will be appropriate
for the forthcoming construction to take $\v , \w$ as the basis in
$\cc^2$, where $A$ acts, and to introduce coordinates $Z=(z_1,z_2)$
in this basis, \ie now we have: $\v =(1,0)$ and $\w=(0,1)$. In these
coordinates $A$ acts as $AZ = (\alpha z_1, 1/\alpha z_2)$. Observe
that our lattice $\Lambda_0$ is irrational in these coordinates. A
foliation on $X_0$ we construct as follows. Take first the
"vertical" foliation $\{z_1=\const \}$ in $\cc^2$, factor it by
$\Lambda_0$. Due to the irrationality of $\Lambda$ in the new basis
it will have dense leaves. Now we observe that this foliation is
obviously invariant under the action of $A$, which is simply
multiplication by $1/\alpha$ on the leaves. Therefore the "vertical"
foliation $\call^{\v} = \{z=\const , z_1=\const_1\}$ descends from
$\cc^*\times\cc^2$  to $X_0$ and we denote it as $\call_0$.

\smallskip Now, following \cite{Na}, we shall deform
$(X_0,\call_0)$. In the subspace
$\cc^2_{z,z_1}\deff\cc_z\times\cc_{z_1}$ of our coordinate space
$\cc^3_{z,Z}\deff\cc_z\times\cc^2_{z_1,z_2}$ we take a real subspace
$\rr^2_{\tau}$ - a deformation of $\{0\}_z\times\cc_{z_1}$.
Parameter $\tau$ here runs in $Gr_{\rr}(2,4)$. This subspace
$\rr^2_{\tau}\subset\cc^2_{z,z_1}$ we see as the graph of the
uniquely defined $\rr$-linear map $L_{\tau}: \cc^2_{z_1}\to \cc_z$
and therefore the subspace $\rr^4_{\tau}\deff
\rr^2_{\tau}\times\cc_{z_2}$ is a graph of $(L_{\tau}, \id ):
\cc^2_{z_1,z_2}\to \cc_z$. By $\Lambda_{\tau}$ we denote the image
of the lattice $\Lambda_0$ under $(L_{\tau}, \id )$ - a deformation
of $\Lambda_0$. Denote by $\ttt_{\tau}$ the torus
$\rr^4_{\tau}/\Lambda_{\tau}$. Remark that $A$ still preserves
$\Lambda_{\tau}$ and therefore $\cc^3_{z,Z}\setminus \{0\}_z\times
\rr^4_{\tau}$ factors first by $\Lambda_{\tau}$ and then by $(\alpha
, A)$ to a compact complex threefold $X_{\tau}$ which is a real
$4$-torus bundle over a complex $1$-torus $\cc^*/<\alpha>$. Our
"vertical" foliation $\call^{\v}$ descends again to $X_{\tau}$ and
we denote the result as $\call_{\tau}$. The construction of
$(X_{\tau},\call_{\tau})$ is finished.
\end{exmp}

In the following Proposition $\calv$ denotes a sufficiently small
neighborhood of $\{0\}\times\cc_{z_1}$ in
$Gr_{\rr}(2,\cc^2_{z,z_1})$.

\begin{prop}
The family of foliated $3$-folds $\{(X_{\tau},\call_{\tau}):\tau\in
\calv\}$, constructed above, possesses the following properties:

\sli Manifolds $X_{\tau}$ do not admit a $dd^c$-closed (even $dd^c$-negative)
metric form for all $\tau\in \calv\setminus\cc\pp^1$ and $X_{\tau}$ is
not even almost Hartogs.

\slii At the same time all $(X_{\tau},\call_{\tau})$ are Hartogs.
\end{prop}
\proof (\sli The fact that $X_{\tau}$ are not K\"ahler is explained
in \cite{BK}, see pp. 82-84. For $\tau\in V\setminus \cc\pp^1$ our
$X_{\tau}$ has $\cc^3\setminus \rr^4_{\tau}$ as an unramified
covering. For this reason it is also not almost Hartogs. Really, the
covering map is singular along $\rr^4_{\tau}$ which is much more
massive then just a countable union of complex curves. But it is
also to massive as a singularity set for the covering map in the
event that $X_{\tau}$ would admit a plurinegative metric form, see
the Main Theorem from \cite{Iv6}.

\smallskip\noindent (\slii Let $h:(W,\pi , U,V)\to (X_{\tau},\call_{\tau})$
be a holomorphic foliated generic injection of a three dimensional
generalized Hartogs figure into the foliated manifold
$(X_{\tau},\call_{\tau})$. Without a loss of generality assume that
$U\subset V$ are bidiscs, so $W$ is simply connected. Lift $h$ to a
foliated generic injection $\tilde h$ of $(W,\pi , U,V)$ into
$(\cc^3\setminus \rr^4_{\tau},\call^{\v})$. Then it extends after a
reparametrization as a map with values in $(\cc^3,\call^{\v})$. But
the fiber of $\call^{\v}$ which touches $\rr^4_{\tau}$ is entirely
contained in $\rr^4_{\tau}$ and therefore the extended map never
hits $\rr^4_{\tau}$. After that we can descend the extended map back
to $(X_{\tau},\call_{\tau})$.

\smallskip\qed

\newprg[prgEXOQ.rat]{Rationality}

First of all we shall prove the following

\begin{lem}
\label{superharm1} Let $\omega$ be a plurinegative taming form for
the vertical foliation on the product $\Delta\times \cc\pp^1$. Then
the volume function
\[
\v_{\omega} (z_1) = \int\limits_{\{z_1\}\times \cc\pp^1} \omega
\]
is superharmonic. If, moreover, $\omega$ is pluriclosed then
$\v_{\omega}$ is harmonic.
\end{lem}
\proof For any test function $\psi$ in $\Delta$ we have
\[
<\psi ,\Delta \v_{\omega}> = \frac{i}{2}\int_{\Delta}\Delta\psi
\left(\int_{\{z_1\}\times \cc\pp^1} \omega \right)d\zeta\land
d\bar\zeta = \int_{\Delta\times\cc\pp^1}dd^c(\pi^*\psi) \land \omega
=
\]
\begin{equation}
\eqqno(areafunction2)
=  \int_{\Delta\times\cc\pp^1}\pi^*\psi
\land dd^c\omega \leq 0,
\end{equation}
\ie $\v_{\omega}$ is superharmonic. If $\omega$ was pluriclosed then
\eqqref(areafunction2) becomes an equality and consequently in this
case $\v_{\omega}$ is harmonic.

\smallskip\qed

\smallskip
Recall that by $\calr_{\call}$ we denoted the analytic space of
rational cycles on $X$ tangent to $\call$. Fix a plurinegative
taming form $\omega$ and consider the area function
$\v_{\omega}:\calr_{\call}\to \rr^+$ defined by
(\ref{areafunction1}).

\begin{corol}
\label{ratgeo}
Suppose that $\call$ is tamed by a plurinegative form
$\omega$. Then every irreducible component of $\calr_{\call}$ is
compact and every connected component consists of finitely many
irreducible ones.  The volume function $\v_{\omega} $ is constant on
every connected component of $\calr_{\call}$.
\end{corol}
\proof Let first that  $\calk$ is an irreducible component of
$\calr_{\call}$. Denote by $\calc_{\calk}$ the universal family over
$\calk$. $\calc_{\calk}$ comes with two natural mappings:
projection $\pi:\calc_{\calk}\to \calk$ and inclusion
$p:\calc_{\calk}\to X$. Take an analytic disc $\phi :\Delta\to
\calk$ and the restriction $\calc_{\Delta}$. Lemma \ref{superharm1}
shows that $\v_{\omega}|_{\phi (\Delta)}$ is superharmonic.
Therefore $\v_{\omega}$ is plurisuperharmonic on $\calk$. Suppose
$\calk$ is not compact. Take a divergent sequence of points
$\{k_n\}\subset \calc$. Now two cases could occur:

\smallskip\noindent{\sl Case 1. $\v_{\omega}(k_n)$ stays bounded
(may be on some subsequence).}

In that case we can subtract a converging subsequence of rational
cycles $C_{k_n}$. The limit is again a rational cycle $C_0$ which
obviously should belong to our irreducible component $\calk$.
Contradiction.

\smallskip\noindent{\sl Case 2. $\v_{\omega}(k_n)\to \infty$.} So $v(k)$
increases when $k$ goes to infinity in $\calk$, \ie leaves every
compact. But this contradicts to the minimum principle for
(pluri)-harmonic functions.

\smallskip Therefore $\calk$ is compact and $v_{\omega}$ is constant on
$\calk$. This implies that $\v$ is constant on every connected
component of $\calr_{\call}$. Suppose there exist a sequence
$\calk_n$ of irreducible components of some connected component
$\caln$. Take $k_n\in \calk_n$. Then $\v_{\omega}(C_{k_n})$ is
constant and therefore some subsequence $C_{k_n}$ converges to some
rational cycle $C_0$ which corresponds to a point $k_0$ in
$\calr_{\call}$. But in this case $\calr_{\call}$ contains a
sequence of compact irreducible components having an accumulation
point $k_0$. This contradicts to the fact that $\calr_{\call}$ is a
complex space.

Therefore each connected component of $\calr_{\call}$ consists from
a finite number of compact irreducible ones.

\qed

The following is immediate:

\begin{corol}
A foliated manifold $(X,\call)$ which admits a plurinegative taming
form has bounded rational cycle geometry.
\end{corol}

Let us turn to the proof Corollary 4 from the Introduction. It is
sufficient to establish the following:

\begin{lem}
\label{ratgeo1} Let $(\tilde\call_D,\pi)$ be a covering cylinder of
holomorphic foliation by curves $\call$ on a compact complex
manifold $X$ which admits $dd^c$-negative metric form. Suppose that
$D$ is biholomorphic to the polydisc and that there exists $z\in D$
such that the fiber $\tilde\call_z=\pi^{-1}(z)$ is isomorphic to
$\cc\pp^1$. Then $\pi^{-1}(D)\sim D\times \cc\pp^1$.
\end{lem}
\proof The set $U$ of $z\in D$ such that $\tilde\call_z\sim
\cc\pp^1$ is clearly open. Each connected component $U'$ of $U$
naturally is included in some irreducible component $\calk$ of
$\calr_{\call}$. Therefore the area function $\v_{\omega} (z) =
\area_{\omega}(\tilde\call_z)$ is constant on $U'$. But this implies
that for any boundary point $z^0\in \d U'\cap D$ the fiber
$\tilde\call_{z^0}$ is again rational. Therefore $U'=D$ and
$\tilde\call_D=D\times\cc\pp^1$.

\qed

\begin{rema} \rm
\label{ratgeo2} Precisely the same argument as in
\eqqref(areafunction2) gives that every irreducible component of a
(not necessarily rational) cycle space of compact curves tangent to
$(X,\call)\in \cale_-$ is compact and the area is constant. This
implies, for example, that if a holonomy covering of some leaf is
compact then the same is true for all generic leaves.
\end{rema}

\begin{exmp}
\label{kato-exmp}
As it is shown in \cite{K3} there exists a compact
complex manifold $X$ of dimension $5$  and a smooth  holomorphic
foliation by curves $\call$ on $X$ such that there exists a
non-empty domain $W\subset X$ with $X\setminus \bar W\not=\emptyset$
having the following properties:

\sli If $z^0\in W$ then $\call_{z^0}\subset W$ and $\call_{z^0}\equiv
\cc\pp^1$.

\slii There exists thin subset $S$ of $X\setminus \bar W$ such all
compact leaves in $X\setminus \bar W$ are contained in $S$.
\end{exmp}

Lemma \ref{ratgeo1} implies now that this $(X,\call)$ doesn't admit
a plurinegative taming form.

\newprg[prgEXOQ.pres]{Preservation of cycles}

 Let $\call$ be
a foliation by curves on a disc-convex complex manifold $X$ and $D$
be a transversal smooth hypersurface. We shall work on the holonomy
covering cylinder $\hat\call_D$. If $\call$ is smooth the same works
also for $\call_D$. Take a point $z\in D$ and a loop $\gamma  \in
\pi_1(\hat\call_z)$. Reference point for $\pi_1(\hat\call_z)$ will
be always $z$.

\begin{defi}
The domain of preservation of the homotopy class $[\gamma ]$ is a
topological space $\Omega_{\gamma , D}$ defined as follows:

\smallskip\noindent
1) the points of  $\Omega_{\gamma , D}$ are homotopy classes
$[\gamma^{'}]\in \pi_1(\hat\call_{z^{'}})$ (where $z^{'}$ is any
point of $D$) such that some representative $\gamma^{'}$ of
$[\gamma^{'}]$ can be joined by a homotopy $\gamma_t$ of loops in
$\hat\call_{z(t)}$ with some representative $\gamma$ of $[\gamma]$.
Here $z(t)$ is a path in $D$ from $z^{'}$ to $z$.

\smallskip\noindent
2) the topology on $\Omega_{\gamma , D} $ is defined in a natural
way saying that $[\gamma_n]$ converge to $[\gamma ]$ if some
representatives converge uniformly.
\end{defi}

Let $\Omega_{\gamma}$ be  the domain of preservation of the
(homotopy class  $[\gamma ]$ in fact) of our loop $\gamma$. There is
a natural projection $p:\Omega_{\gamma}\to D$ sending
$[\gamma^{'}]\in \pi_1(z^{'})$ to $z^{'}$.

\medskip\noindent{\sl Proof of Proposition 5}. Suppose that for some
loop $\gamma \subset \hat\call_z$ the space $\Omega_{\gamma , D}$ is
not Hausdorff. That means that there exists $z^0\in D$, two loops
$\gamma ,\beta \subset \hat\call_{z^0}$ representing different
homotopy classes in $\pi_1(\hat\call_{z^0})$ and two sequences of
loops $\gamma_n,\beta_n\subset \hat\call_{z_n}$, homotopic to each
other in $\hat\call_{z_n}$, converging to $\gamma$ and $\beta$
respectively.

Taking $\alpha_n:=\gamma_n\circ\beta_n^{-1}$ we obtain a sequence of
loops, homotopic to zero and converging to a loop $\alpha\subset
\hat\call_{z^0}$ which is not homotopic to zero.

We are exactly in the situation of the proof of the Theorem 1 and
therefore deduce the existence of a foliated shell in $\call$.

The local biholomorphicity of the projection $p$ is obvious.

\medskip\qed

The phenomena of preservation of cycles to our knowledge was first
studied by Landis-Petrovsky in \cite{LP}, see also \cite{Iy1}.

\newprg[prgEXOQ.comp]{Foliations with compact fibers}

Let $\call $ be a smooth holomorphic foliation of dimension one
on an $n$-dimensional complex manifold $X$. We suppose that all
leaves of $\call$ are compact.

\smallskip\noindent{\sl Proof of Proposition 4.} (i) We denote by
$\omega$ an adapted to $\call$ plurinegative $(1,1)$-form. Let $\call_z$
be a leaf of $\call$ through the
point $z\in X$. If $\call_z$ is compact with finite holonomy we
denote by $n(\call_z)$ the cardinality of the holonomy group of
$\call_z$ and set

\begin{equation}
\v (z) = \vol (\call_z) = n(\call_z)\int_{\call_z}\omega .
\end{equation}

Denote by $\Omega$ the connected component of the set of $z\in X$
such that the leaf $\call_z$ of $\call$ through $z$ is compact and
has finite holonomy which contains our compact leaf. By the Reeb
local stability theorem $\Omega$ is an open set in $X$.

\smallskip\noindent{\sl Case 1.} {\it There exists $z^0\in \d\Omega$
which is a limit of $z_n\in \Omega$ with $\v (z_n)$ uniformly
bounded.}

\smallskip For any transversal $D$ to $\call_{z^0}$ the intersection
$D\cap \Omega$ is open in $D$ and every $\call_{z_n}$ cuts $D$ by a
bounded number of points, say $N$. This readily follows from the
boundedness of volumes of $\call_{z_n}$. Therefore for every $h\in
\hol (\call_{z^0})$ its order is at most $N!$, \ie $h^{N1}=\id $.
Therefore the holonomy group $\hol (\call_{z^0})$ has finite
exponent and therefore it is finite itself, see Lemma 2 from
\cite{P}. Therefore $z^0$ is an interior point of $\Omega$.
Contradiction.

We are left with the following possibility:

\smallskip\noindent{\sl Case 2.} {\it $\v (z)\to\infty$ when $z\to\d\Omega$.}

This case is excluded by Remark \ref{ratgeo2}. All is left to remark
that if $\d\Omega\not=\emptyset$ we obtain a contradiction with the
minimum principle for plurisuperharmonic functions.

\medskip Therefore $\Omega = X$ and (i) is proved.

\medskip\noindent (ii) By the standard observation in foliation
theory, see ex. \cite{Go} the set of leaves without holonomy is not
thin in $X$. Therefore we are done by (i).

\medskip\qed

\begin{rema}\rm
Without any changes this proof applies to smooth holomorphic
$q$-dimensional foliations on compact complex manifolds admitting
$dd^c$-negative taming $(q,q)$-forms.
\end{rema}

\newprg[prgEXOP.comp-fibr]{Holomorphicity of Complex Fibrations}

Let us give one more corollary of the compactness of the spaces of
cycles tangent to a foliation admitting a plurinegative adapted
form. Denote by $\calb_{\call}$ the space of $q$-cycles tangent to $\call$.
By $\calz $ we denote
the corresponding universal family and let $\ev :\calz \to X$ be the
evaluation map.

\begin{thm}
Let $\call$ be a smooth real foliation on compact complex manifold $X$
with all leaves being compact complex manifolds of complex dimension $q$.
Suppose that $\call$ admits a plurinegative adapted $(q,q)$-form $\omega$
and that $q\ge \frac{1}{2}\dim_{\cc}X$. Then $\call$ is holomorphic.
\end{thm}
\proof
It is easy to see that $\calb_{\call}$ has at most countably many
irreducible components, they are all compact and the volume function is
constant on each of them. The argument is similar to that of
Lemma \ref{superharm1}. That immediately gives that some irreducible
component $\calk$ of $\calb_{\call}$ such that $\ev (\calz_{\calk})$ has
positive measure and therefore $\ev (\calz_{\calk}) = X$.

\smallskip Let us prove that $\calk$ contains all leaves of $\call$.
Suppose that there is a leaf $\call_z$ which is not in $\calk$. Find
$k\in\calk$ such that $\calz_k\cap\call_z\not=\emptyset$. This
intersection represents a nontrivial element in corresponding
homology group. And this homological intersection doesn't depend on
the choice of $k\calk$ and $\call_z$. This is a contradiction,
because for some $k$ the cycle $\calz_k$ coincides with a leaf of
$\call$.

Therefore all leaves of $\call$ belong to $\calk$. In the same manner one establishes
that all cycles $\calz_k$ from $\calk$ are leaves of $\call$.

\medskip\qed

\begin{rema}\rm
In dimension two, \ie when $X$ is a compact complex surface ,this result
is due to J. Winkelmann, see \cite{W}.
\end{rema}

\newprg[prgEXOP.iwasawa]{Foliations on Iwasawa manifold}

Example \ref{kato-exmp} of Kato already provided us a foliation
without a plurinegative taming form. However it is very inexplicit.
Let us give a very simple one.
\begin{exmp}\rm
Let $H(3)$ be the group of matrices of the form

\begin{equation}
A =
\left(%
\begin{array}{ccc}
  1 & z_1 & z_3 \\
  0 & 1 & z_2 \\
  0 & 0 & 1 \\
\end{array}%
\right)
\end{equation}
with complex $z_1,z_2,z_3$. Denote by $\zz (3)$ the subgroup of
$H(3)$ which consists from $z_1,z_2,z_3\in \zz + i\zz$. The quotient
$H(3)/\zz (3)$ is a compact, complex three-dimensional manifold
$\cali $ - Iwasawa manifold. The holomorphic forms $\omega_1=dz_1$,
$\omega_2=dz_2$ and $\omega_3=dz_3-z_1dz_2$ are left invariant with
respect to the action of $\zz (3)$ and therefore project to
holomorphic forms on $\cali$. Define a holomorphic foliation by
curves $\call_1$ on $\cali $ by $\omega_1=\omega_2=0$.
\end{exmp}

\begin{prop}
Foliated manifold $(\cali , \call_1)$ possesses the following
properties:

\sli It is Hartogs.

\slii It doesn't admit a plurinegative taming form.
\end{prop}
\proof (\sli Hartogs property is invariant with respect to
unramified coverings. Since the universal covering of $\cali$ is
$H(3)\equiv \cc^3$ we are done.

\smallskip\noindent (\slii Consider $S\deff \frac{i}{2}\omega_3\wedge
\bar\omega_3$ as a positive $(2,2)$-current on $\cali$. A simple
calculation
\[
dd^cS=i\d\bar\d S = \d\omega_3\wedge\bar\d\omega_3 =
\frac{i^2}{2}\omega_1\wedge\bar\omega_1\wedge\omega_2\wedge\bar\omega_2=:T.
\]
And $T$ is a positive current directed by $\call_1$. A positive
current $S$ such that $dd^cS$ is also positive and directed by
$\call$ is a clear obstruction to the existence of a plurinegative
taming form for $\call_1$.

\smallskip\qed

\begin{rema}\rm
Iwasawa manifold carries also other foliations. For example
$\call_2\deff \{\omega_1=\omega_3=0\}$. This $\call_2$ is tamed by
the closed form $\frac{i}{2}\omega_2\wedge\bar\omega_2$ and
therefore the Proposition 1 from the Introduction applies to
$\call_2$ - it belongs to $\calu$.
\end{rema}

\newprg[prgEXOQ.quest]{Open questions}

In this subsection we shall formulate some open questions which are
important to complete our knowledge about holomorphic foliation by
curves on compact complex manifolds and which naturally come out from
the discussions in this paper.

\medskip\noindent{\bf Class $\cals$.}

\smallskip Recall that foliations of class $\cals$ are pluritamed
foliations containing foliated shells.

\begin{quest}
In the conditions of the Theorem 3.1 let $\call_z$ be the leaf which
contains an essential vanishing cycle. Is it true that its closure
$\bar\call_z$ is a compact complex curve?
\end{quest}
Recall that this is like the proof of Novikov's theorem works, see
\cite{Go} for example. If $\dim X=2$ we proved it in Corollary \ref{dim-2}.

\begin{quest}
Prove that a foliation with shells is parabolic.
\end{quest}

\begin{quest}
Let $(X,\call)\in \cals$ be compact and $\dim X\geq 3$. Does $(X,\call)$
contains  a total space of a deformation $(\calx,\pi, Y)$ of foliated compact
Hopf or Kato surfaces $\calx_{\lambda}, \lambda\in Y$ with compact $Y$?
\end{quest}

Most probably the answer to this question is ``no'' if stated in its full
generality. But even a partial positive result (or a counterexample) would
be important. In this concern let us give an example shoving that the total
space of deformation may not sweep the whole of $X$, \ie that a foliated
shell may "disappear in the limit".

\begin{exmp}\rm Let $E'$ be a holomorphic rank two bundle over a Hopf
surface $H^2=\cc\setminus\{0\}/z\sim 2z$ which admits a holomorphic
section $\sigma$ vanishing exactly at one point $z_0\in H^2$ with
multiplicity one, see \cite{GH}, p.726. Denote by $E$ the bundle
dual to $E'$. Let $\tau_0$ be the zero section of $E$. The quotient
of $E\setminus \tau_0$ by the action $(z,\v)\to(z,\frac{1}{2}\v)$ is
a compact complex $4$-manifold which we denote as $X$. It is fibered
over $H^2$ and the fiber over $z\in H^2$ we denote as $X_z$.

\smallskip $E\setminus \tau_0$ carries a singular holomorphic
foliation by curves defined as follows: its leaves in each fiber
$E_z\setminus\{0\}, z\not=z_0$ are $\{ x\in E_z:
\sigma_z(x)=const\}$. Actually on each $E_z\setminus\{0\}$ it is
again our "vertical" foliation. It factors under the chosen action
to a foliation $\call$ on $X$. The singularity set of $\call$ is
$E_{z_0}$. $(X,\call)$ carries an obvious family of foliated shells
over $H^2\setminus \{z_0\}$, and this family extends over $z_0$ (!)
as a family of shells. But $\call$ itself is singular over $z_0$ and
therefore the shell in $X_{z_0}$ is not a foliated one.
\end{exmp}

\begin{quest}
Is it true that immersed foliated shells could be always made spherical? The
same question about imbedded ones. In that case one expects them to be
holomorphic foliated images of quotients of the standard sphere in $\cc^2$
with the standard vertical foliation.
\end{quest}

The problem here lies in reducing of the size of the singularity set
$S$, see Subsection 4.1 for a more detailed discussion.

\begin{quest}
Let $(X,\call)$ be pluritamed by a $dd^c$-closed {\slsf metric} form
$\omega$. Is it true that the singularity set $S$ of meromorphic
foliated immersions appearing in Theorem \ref{reparalldim2} is at
most countable union of analytic subsets of pure  codimension two?
\end{quest}

\medskip\noindent{\bf Class $\calu$.}

\begin{quest}
Let $D$ be a transversal polydisc. Suppose that the skew cylinder
$\tilde\call_D$ exists (and $\call$ admits a plurinegative adapted
form). Prove that $\tilde\call_D$ is disc-convex.
\end{quest}

This is known for Stein $X$, \cite{Iy1}, in that case
$\tilde\call_D$ is Stein. It is also known for algebraic $X$,
\cite{Br3}.

\begin{quest}
Prove that the set of $z\in D$ such that $\tilde\call_z =\cc$ is pluripolar in
$D$ or is the full $D$.
\end{quest}

Algebraic case is treated in \cite{Br3}.

\begin{quest}
Suppose that the  domain $\Omega_{\gamma}$  of preservation of
cycle $\gamma$ (as in Definition 1.4) exists.
Prove that $\Omega_{\gamma }$ is good in the sense of
Landis-Petrovsky, \ie that for a natural projection
$p:\Omega_{\gamma } \to D$ the set $\Sigma :=D\setminus
p(\Omega_{\gamma })$  doesn't separate $D$.
\end{quest}

This question for holomorphic foliations by curves on arbitrary
compact complex manifolds Il'yashenko calls the generalized
Landis-Petrovsky conjecture. The answer is positive for algebraic
foliations, but it is wrong for holomorphic ones on Stein manifolds,
see \cite{Iy2}. Example of Kato in \cite{K3} leaves little hope for
the positive answer in general, but it is not a direct
counterexample. In Question 7
We propose to solve the Landis-Petrovsky conjecture when $\call$
admits an adapted pluriclosed (or plurinegative) taming form.

\newprg[prg5.8]{Class $\cale$}

\begin{quest}
Suppose $T$ is a nontrivial $dd^c$-exact foliated cycle for $\call$
such that $T=dd^cS$ where $S$ is also positive. Can one provide more
restrictions on the support of $T$ then it is done in Theorem
\ref{pluripos}? Can one say something about the structure os $S$?
\end{quest}

Not that $S$ has a well defined Lelong numbers, see \cite{Sk}.

\begin{quest}
Is it true that every pluriexact foliation contains a compact curve tangent to the leaves?
\end{quest}

\ifx\undefined\bysame
\newcommand{\bysame}{\leavevmode\hbox to3em{\hrulefill}\,}
\fi

\def\entry#1#2#3#4\par{\bibitem[#1]{#1}
{\textsc{#2 }}{\sl{#3} }#4\par\vskip2pt}

\end{document}